\documentclass[11pt]{elsart-maktor}
\renewcommand{\and}{et}
\usepackage{amssymb,
amsfonts,amssymb,euscript,mathrsfs,pifont
}
\usepackage[english,francais]{babel}
\usepackage{amsxtra}
\usepackage{amscd}
\usepackage[T1]{fontenc}                
\usepackage[all]{xy}

 \DeclareMathAlphabet{\got}{U}{euf}{m}{n}     
\DeclareMathAlphabet{\mat}{U}{msb}{m}{n}     
\DeclareMathAlphabet{\mathbold}{OML}{cmm}{bx}{it} 

\makeatletter

\renewcommand{\section}{\@startsection{section}{1}%
\z@{-3.5ex \@plus -1ex \@minus -.2ex}
  {2.3ex \@plus.2ex}
  {\bfseries\large}}           

\renewcommand{\subsection}{\@startsection{subsection}{2}%
\z@{.5ex \@plus.7ex}{-.5em}%
{\normalfont\bfseries}}

\renewcommand{\subsubsection}{\@startsection{subsubsection}{3}%
\z@{.5ex\@plus.7ex}{-.5em}%
  {\normalfont\itshape}}
\makeatother
\numberwithin{equation}{section}

\newtheorem{theo}{Th\'eor\`eme}[subsection]
\newtheorem{pr}{Proposition}[subsection]
\newtheorem{lem}{Lemme}[subsection]
\theoremstyle{definition}
\newtheorem{defi}{D\'efinition}[subsection]
\newtheorem{co}{Corollaire}[subsection]

\newcommand{\transposee}[1]{{\vphantom{#1}}^{\mathit t}{#1}}
\newcommand{\goth}{\mathfrak}

\newenvironment{rem}{\vskip 1em\bf Remarque.\rm}{\par\rm}
\newenvironment{rems}{\vskip 1em\bf Remarques.\rm}{\par\rm}
\newenvironment{theo*}{\vskip 1em\bf Th\'eor\`eme.\it}{\par\rm}
\newenvironment{dem}{\vskip 1em{\it D\'emonstration} :}%
{\unskip\hfill\null\nobreak\hfill\carre\vskip1em\par}
\newcommand{\carre}{\rule{1ex}{1ex}} 

\renewcommand{\thefootnote}{\fnsymbol{footnote}}

\makeatletter
\providecommand*{\diff}%
{\@ifnextchar^{\DIfF}{\DIfF^{}}}
\def\DIfF^#1{%
\mathop{\mathrm{\mathstrut d}}%
\nolimits^{#1}\gobblespace}
\def\gobblespace{%
\futurelet\diffarg\opspace}
\def\opspace{%
\let\DiffSpace\!%
\ifx\diffarg(%
\let\DiffSpace\relax
\else
\ifx\diffarg[%
\let\DiffSpace\relax
\else
\ifx\diffarg\{%
\let\DiffSpace\relax
\fi\fi\fi\DiffSpace}

\providecommand*{\eu}%
{\ensuremath{\mathrm{e}}}

\providecommand*{\iu}%
{\ensuremath{\mathrm{i}}}

\newcommand{\im}{\mathop{\mathrm{Im}}}

\newcommand{\Supp}{\mathop{\mathrm{Supp}}}

\newcommand{\Ind}{\mathop{\mathrm{Ind}}\nolimits}

\newcommand{\Card}{\mathop{\mathrm{Card}}}

\newcommand{\tr}{\mathop{\mathrm{tr}}}
\newcommand{\norm}{\mathop{\mathrm{N}}}

\def \Z {\Bbb{Z} }

          \def\var{\varphi}

            \def \O { \mathcal O}

\def \O{\mathcal O}      
\def\var{\varphi}      
    \def \a{\goth a}

     \def \e { \varepsilon }

\begin{document}

\begin{frontmatter}
\title{Restriction de la représentation de Weil à un sous-groupe
  compact maximal ou à un tore maximal elliptique.}

\author{Khemais Maktouf}
\address{Université de Monastir, Faculté des Sciences de Monastir,
Département de Mathématiques,
5019 Monastir, Tunisie}
\ead{khemais.maktouf@fsm.rnu.tn}

\author{Pierre Torasso} \address{UMR 6086 CNRS, Université de Poitiers,
  Laboratoire de Mathématiques et Applications, Boulevard Marie et
  Pierre Curie, BP~30179,
  86962 Chasseneuil Cedex, France}
\ead{pierre.torasso@math.univ-poitiers.fr}

\selectlanguage{french}
\begin{abstractf}
Dans cet article nous démontrons que la représentation de Weil sur un
corps $p$-adique avec $p\neq2$ restreinte à un sous-groupe compact
maximal ou un tore maximal elliptique du groupe métaplectique se
décompose sans multiplicité et décrivons explicitement les
représentations irréductibles ou les caractères qui interviennent.
\end{abstractf}
\begin{keywordf}
représentation de Weil, groupe métaplectique, sous-groupe compact
maximal, tore maximal elliptique \MSC 22E50
\end{keywordf}
\selectlanguage{english}
\begin{abstract}
In this article, we prove that the restriction of the Weil
representation over a $p$-adic field, $p\neq2$, to maximal compact
subgroups or elliptic maximal tori of the metaplectic group is
multiplicity free and give an explicit description of the irreducible
representations or characters occurring.
\end{abstract}
\begin{keyword}
Weil representation, metaplectic group, maximal compact subgroup,
elliptic maximal torus \MSC 22E50
\end{keyword}
\end{frontmatter}
\selectlanguage{french}

\section{Introduction}\label{1}

La représentation de Weil intervient dans de nombreux domaines de la
théorie des représentations des groupes presque algébriques généraux
ou des groupes réductifs sur un corps local. Citons la construction,
dans le cadre de la méthode des orbites, des représentations unitaires
irréductibles dans \cite{duflo-1982}, la correspondance de Howe (voir
\cite{prasad-ip-1997}, \cite[Chapitre
  2]{moeglin-vigneras-waldspurger-1987}, \cite{waldspurger-ip-1990}),
les séries theta \cite{lion-vergne-1980}. Il est donc important d'en
comprendre la structure et notamment comment sa restriction à un
sous-groupe compact maximal ou un tore anisotrope se décompose en
irréductibles. Dans cet article, nous nous plaçons sur un corps
$p$-adique de caractéristique résiduelle différente de $2$ et nous
démontrons que les représentations des sous-groupes compacts maximaux
et les caractères des tores maximaux elliptiques qui apparaissent
effectivement le font avec la multiplicité $1$ et les décrivons
explicitement.

Soit $k$ un corps local non archimédien de caractéristique nulle, $\O$
l'anneau des entiers de $k$, $\goth p$ son idéal maximal et $\varpi$
une uniformisante de $\O$.  Le corps résiduel $\O/\goth p$ (not\'e
$\Bbb{F}_{q}$) est fini de cardinal $q$ et de caractéristique $p$.
Nous supposons que $p\neq 2$.  Nous fixons un caractère unitaire
$\psi$ de conducteur $\lambda_\psi$ de $k$.

On se donne un $k$-espace symplectique $(W ,\beta)$ de dimension
$2r$. On note $Sp(W)$ le groupe symplectique associé à $(W,\beta)$ et
$Mp(W)$ le groupe m\'etaplectique correspondant, revêtement à deux
feuillets non trivial de $Sp(W)$. On a une suite exacte
courte~:
\begin{eqnarray}\label{eq0.0.0}
1 \longrightarrow \Z/2\Z \longrightarrow
Mp(W) \longrightarrow Sp(W)
\longrightarrow 1.
\end{eqnarray}
En fait, le groupe de Heisenberg $H(W)$ construit sur l'espace
symplectique $W$ admet une unique (à équivalence près) représentation
unitaire irréductible $(\rho_{\psi},\mathcal{H})$ de caractère central
$\psi$~; on l'appelle la représentation de
Schrödinger de caractère central $\psi$ de $H(W)$. Le groupe $Sp(W)$
opère dans $H(W)$ par authomorphismes agissant trivialement sur le
centre. Soit $U(\mathcal{H})$  le groupe des transformations
unitaires de l'espace de Hilbert $\mathcal{H}$, muni de la topologie
de la convergence forte.  L'ensemble
$\widehat{Sp(W)}_{\psi}$ des couples $(g,U)\in Sp(W)\times
U(\mathcal{H})$ vérifiant~:
\begin{equation*}
U\rho_{\psi}(h)U^{-1}=\rho_{\psi}(g.h)\mbox{, }g\in Sp(W)\mbox{,
}h\in H(W),
\end{equation*}
est un sous-groupe fermé localement compact de $Sp(W)\times
U(\mathcal{H})$, extension centrale de $Sp(W)$ par le groupe $\Bbb{U}$
des nombres complexes de module $1$. On a donc une suite exacte courte
\begin{eqnarray}\label{eq0.0.1}
  1 \longrightarrow \Bbb{U} \longrightarrow \widehat{Sp(W)}_{\psi}
\longrightarrow Sp(W) \longrightarrow 1
\end{eqnarray}
et l'application $S_{\psi}:(g,U)\mapsto U$ est une
représentation de $\widehat{Sp(W)}_{\psi}$, appelée représentation
métaplectique. Il existe un unique homomorphisme de groupes 
$Mp(W)\longrightarrow\widehat{Sp(W)}_{\psi}$ rendant commutatif le
diagramme suivant
\begin{equation*}
\xymatrix{ Mp(W) \ar[rr] \ar[rd] && \widehat{Sp(W)}_{\psi}  \ar[ld]
  \\ & Sp(W) } 
\end{equation*}
Ce morphisme est injectif et permet d'identifier $Mp(W)$ à un
sous-groupe de $\widehat{Sp(W)}_{\psi}$. Alors la restriction à
$Mp(W)$ de la représentation métaplectique est une représentation
fidèle encore notée $S_{\psi}$ et appelée représentation de Weil de
type $\psi$ de $Mp(W)$.

Nous montrons que la suite exacte \ref{eq0.0.0} est scindée au dessus
de chaque sous-groupe compact maximal de $Sp(W)$. 
On sait qu'un sous-groupe compact
maximal de $Sp(W)$ est le stabilisateur $K_{B}$ d'un bon réseau $B$ de
$W$, i. e. $B$ est un sous-$\mathcal{O}$-module ouvert et compact de
$W$ vérifiant
\begin{equation*}
\varpi B^* \subset B \subset B^*,
\end{equation*}
où $B^*$ est le réseau dual de $B$ relativement à $\beta$ et
$\psi$. Dans \cite[II.3]{waldspurger-ip-1990}, Waldspurger a associé à
un tel réseau une réalisation
$(\rho^{B}_{\psi},\mathcal{H}_{\psi}^{B})$ de la représentation de
Schrödinger, appelée modèle latticiel généralisé, et une
représentation $S_{\psi}^{B}$ de $K_{B}$ dans $\mathcal{H}_{\psi}^{B}$
telle que l'application $g\mapsto(g,S_{\psi}^{B}(g))$ soit une section
de la suite exacte courte \ref{eq0.0.1}. Lorsque $B$ est un réseau
autodual, on parle de modèle latticiel~: dans ce cas, Moeglin a montré
que cette application est une section canonique, et unique lorsque
$q\geq 4$, de la suite exacte \ref{eq0.0.0} (voir \cite[Chapitre 2, II
  10, Lemme]{moeglin-vigneras-waldspurger-1987}). Nous montrons que ce
résultat reste vrai dans le cas général (voir les théorèmes
\ref{theo4.2.1} et \ref{theo4.6.1}). Pour ce faire, nous comparons le
modèle latticiel généralisé associé au bon réseau $B$ avec le modèle
latticiel associé à un réseau autodual $A$ tel que $B\subset A\subset
B^{*}$ (voir le théorème \ref{theo4.5.1}). Nous obtenons en
particulier les formules explicites pour la réalisation, dans ce
modèle latticiel, de la représentation $S_{\psi}^{B}$ de Waldspurger
pour les éléments d'un système de générateurs du groupe $K_{B}$
constitué du sous-groupe \flqq parabolique\frqq~ $P_{B}$ stabilisateur
de $A$ et d'un élément particulier $\varsigma_{B}\in K_{B}\backslash
P_{B}$ (voir le corollaire \ref{co4.5.1}). Il est à remarquer que
l'énoncé du théorème \ref{theo4.5.1} se trouve déjà dans
\cite[4.3.e]{pan-2001}.

Dans la suite de cette introduction, nous noterons $S_{\psi}^{B}$ la
réalisation de la représentation de Weil dans le modèle latticiel
associé au bon réseau $B$. On voit donc que la représentation
$S_{\psi}^{B}$ de Waldspurger est la restriction de la représentation
de Weil $S_{\psi}^{B}$ à l'image de $K_{B}$ par sa section canonique
dans le groupe métaplectique. De même, les formules que nous avons
obtenues pour la réalisation de la représentation de Waldspurger de
$K_{B}$ dans le modèle latticiel décrivent la restriction à cette
section de la représentation de Weil $S_{\psi}^{A}$.

Ces dernières formules, nous permettent d'étudier la décomposition en
irréductibles de la restriction de la représentation de Weil à $P_{B}$
et à $K_{B}$. Dans les deux cas, nous montrons que cette décomposition
est sans multiplicité~; nous montrons également que les
représentations de $P_{B}$ qui apparaissent sont monomiales et en
donnons une description explicite (voir le lemme \ref{lem5.1.2} et le
théorème \ref{theo5.2.1}). Lorsque le réseau $B$ est autodual, les
représentations de $K_{B}$ qui apparaissent sont également
monomiales~; dans ce cas le résultat a été obtenu par Prasad (voir
\cite{prasad-ip-1997}).

Si $n$ est un entier naturel, on désigne par $O_{n}$ l'anneau
$\mathcal{O}/\varpi^{n+1}\mathcal{O}$ et on pose
$\mathrm{b}_{n}=B/\varpi^{n+1}B^{*}$ et
$\mathrm{b}_{n}^{*}=B^{*}/\varpi^{n}B$. Alors, $O_{n}$ est un anneau
local fini et principal de caractéristique différente de $2$,
$\mathrm{b}_{n}$ (resp. $\mathrm{b}_{n}^{*}$) est un $O_{n}$-module de
type fini muni naturellement d'une structure symplectique et l'action
de $K_{B}$ dans $B$ passe au quotient à $\mathrm{b}_{n}$
(resp. $\mathrm{b}_{n}^{*}$), induisant un morphisme surjectif de
$K_{B}$ sur le groupe symplectique $Sp(\mathrm{b}_{n})$
(resp. $Sp(\mathrm{b}_{n}^{*})$). On remarque que $O_{0}=\Bbb{F}_{q}$
et que $\mathrm{b}_{0}$ (resp. $\mathrm{b}_{0}^{*}$), que nous notons
plus simplement $\mathrm{b}$ (resp. $\mathrm{b}^{*}$) est un espace
symplectique sur $\Bbb{F}_{q}$~; on note $2l$ la dimension de
$\mathrm{b}^{*}$. Alors le nombre $l=l(B)$, qui prend toutes les
valeurs entières entre $0$ et $\frac{1}{2}\dim W$, détermine la classe
de conjugaison du bon réseau $B$ et donc celle du sous-groupe compact
maximal $K_{B}$. On montre que $K_{B}$ est naturellement isomorphe à
la limite projective des groupes $Sp(\mathrm{b}_{n})$
(resp. $Sp(\mathrm{b}_{n}^{*})$) (voir le lemme \ref{lem2.5.2}).

Dans \cite{cliff-mcneilly-szechtman-2000} et
\cite{cliff-mcneilly-szechtman-2003}, Cliff, McNeilly et Szechtman ont
construit les représentations de Weil pour le groupe symplectique d'un
module symplectique sur un anneau local fini de caractéristique
différente de $2$ et décrit leur décomposition en irréductibles dans
le cas où l'anneau est principal et le module symplectique libre (les
différentes représentations de Weil associées à un même caractère
primitif de l'anneau local diffèrent d'un caractère du groupe
symplectique). Il s'avère que pour $n\geq1$, le $O_{n}$-module
$\mathrm{b}_{n}$ (resp. $\mathrm{b}_{n}^{*}$) est libre si et
seulement si le réseau $B$ est autodual ($l(B)=0$) ou vérifie
$B=\varpi B^{*}$ ($l(B)=r$). Nous étendons les résultats de
\cite{cliff-mcneilly-szechtman-2003} au cas des groupes
$Sp(\mathrm{b}_{n})$ (resp. $Sp(\mathrm{b}_{n}^{*})$) avec $B$
quelconque (voir le théorème \ref{theo3.5.1}). En fait la restriction
de la représentation de Weil $S_{\psi}^{B}$ à $K_{B}$ est la limite
inductive de représentations de Weil des groupes
$Sp(\mathrm{b}_{2n+1})$. Utilisant alors nos résultats, nous en
déduisons une description des représentations irréductibles
apparaissant dans la décomposition de la restriction de la
représentation de Weil à $K_{B}$ comme représentations induites à
partir de stabilisateurs génériques de $K_{B}$ dans les différents
$\mathrm{b}_{2n+1}$ (voir le théorème \ref{theo5.3.1} et le lemme
\ref{lem5.4.1}). Dans \cite{cliff-mcneilly-szechtman-2003} Cliff,
McNeilly et Szechtman ont remarqué que leurs résultats permettaient
d'obtenir la décomposition en irréductible de la restriction de la
représentation de Weil à $K_{B}$ lorsque $B$ est autodual ou vérifie
$B=\varpi B^{*}$.

D'autre part, dans \cite{dutta-prasad-2010} Dutta et Prasad ont
démontré que la représentation de Weil du groupe symplectique
construit sur un groupe abélien fini se décompose sans
multiplicité. Ils ont décrit cette décomposition en terme de la
combinatoire des orbites du groupe abélien sous l'action de son groupe
d'automorphismes et remarqué que l'on peut en déduire les mêmes
résultats pour la représentation du groupe symplectique d'un module
symplectique de type fini sur un anneau fini local et principal. Leurs
résultats contiennent ceux du théorème \ref{theo3.5.1}, sauf que leur
description des sous-modules irréductibles ne les fait pas
explicitement apparaître comme modules induits.

Un tore $T$ d'un groupe réductif est dit elliptique s'il ne contient
pas de sous-tore déployé non central. Un tore elliptique de $Sp(W)$
est donc anisotrope et compact. Par suite, le revêtement métaplectique
est scindé au dessus d'un tel tore. Nous étudions la restriction de la
représentation de Weil à un tore maximal elliptique de $Sp(W)$,
identifié à un sous-groupe de $Mp(W)$ par le choix d'une section.

Tout d'abord, nous démontrons, comme conséquence de la conjecture de
Howe, prouvée par Waldspurger dans \cite{waldspurger-ip-1990}, que la
restriction de la représentation de Weil à un tel tore est sans
multiplicité (théorème \ref{theo6.3.1}). Cette démonstration, plus
simple que celle que nous avions proposée, nous a été communiquée par
le {\it referee}. Nous l'en remercions vivement.

Cependant, ce résultat ne permet pas de décrire complètement cette
restriction. Pour ce faire, nous avons besoin de connaître la
structure des tores maximaux elliptiques de $Sp(W)$.

Soit $T$ un tore maximal de $Sp(W)$. Comme $T$-module, $W$ se
décompose de manière unique en somme directe de sous-modules
irréductibles sur $k$, $W=W_{1}\oplus\cdot\oplus W_{n}$. D'après
\cite{morris-1991}, $T$ est elliptique si et seulement si cette
décomposition est une somme directe orthogonale de sous-espaces
symplectiques. Dans ce cas, $T=T_{1}\times\cdot\times T_{n}$ où chaque
$T_{i}$, image de $T$ par l'application $x\mapsto x_{\vert W_{i}}$,
est un tore maximal elliptique de $Sp(W_{i})$ tel que $W_{i}$ soit un
$T_{i}$-module irréductible sur $k$. Utilisant ce fait, on montre que
l'étude de la restriction de la représentation de Weil aux tores
elliptiques maximaux de $Sp(W)$ se ramène au même problème pour ceux
de ses tores maximaux dont l'action dans $W$ est irréductible sur $k$
(voir la proposition \ref{pr6.3.1} et le théorème
\ref{theo7.4.1}). Pour simplifier, nous dirons qu'un tore $T$ de
$Sp(W)$ tel que le $T$-module soit irréductible sur $k$ est
irréductible.

Soit donc $T$ un tore maximal irréductible de $Sp(W)$. D'après
\cite{morris-1991}, il existe une extension $k'$ de degré $r$ de $k$,
une extension quadratique $k''$ de $k'$, un élément $u\in
k''^{\times}$ de trace nulle relativement à $k'$ et un isomorphisme
$\theta$ d'espaces symplectiques de $(k'',\beta_{u})$ sur $W$ tels que
$T$ soit l'image de $T_{k'',k'}$ par l'isomorphisme
$\tilde{\theta}:t\mapsto \theta\circ t\circ\theta^{-1}$, où la forme
$\beta_{u}$ est définie par
\begin{equation*} \beta_{u}(x,y)=\frac{1}{2}\sideset{}{_{k''/k}}\tr
  ux^{\tau}y\mbox{, }x,y\in k'',
\end{equation*} 
$\tau$ désignant l'élément non trivial du groupe de Galois de $k''$
sur $k'$, et $T_{k'',k'}$ est le sous-groupe multiplicatif de
$k''^{\times}$ constitué des éléments de norme $1$ relativement à
$k'$. On note $\mathcal{O}'$ (resp.  $\mathcal{O}''$) l'anneau des
entiers, $\varpi'$ (resp. $\varpi''$) une uniformisante, $v'$
(resp. $v''$) la valuation normalisée et $q'$ (resp.  $q''$) le
cardinal du corps résiduel de $k'$ (resp.  $k''$).  On prend
$\varpi''$ tel que $\varpi''^{2}=\varpi'$ lorsque $k''$ est ramifié
sur $k'$, et $\varpi''=\varpi'$ dans le cas contraire.  On montre que
l'on peut prendre $u=\varpi''$ lorsque $k''$ est ramifié sur $k'$ et
$v''(u)\in\{0,1\}$ dans le cas contraire (voir le théorème
\ref{theo6.1.1}).

Si $j\in\mathbb{N}$, on désigne par $T_{j}$ le $j$-ième sous-groupe de
congruence de $T$~:
\begin{equation*}
T_{j}=\{g\in T\vert g-1\in\varpi''^{j}\mathcal{O}''\}.
\end{equation*}
On a $T_{0}=T$ et, pour $j\geq1$, $T_{j}$ est un sous-groupe strict de
$T$. Lorsque $k''$ n'est pas ramifié sur $k'$, la suite $T_{j}$ est
strictement décroissante. Dans le cas contraire, la suite $T_{2j+1}$
est strictement décroissante et on a $T_{2j}=T_{2j+1}$ pour $j>0$. 

Soit $j\in\mathbb{N}$ et $\chi$ un caractère de $T_{j}$. On appelle
conducteur de $\chi$, le plus petit entier $\lambda$ tel que
$T_{\lambda}$ soit contenu dans le noyau de $\chi$. 

Lorsque $k''$ n'est pas ramifié sur $k'$, le groupe quotient $T/T_{1}$
est canoniquement isomorphe au groupe $\mu_{q'+1}$ des racines
$(q'+1)$-ièmes de l'unité dans $k''$. L'ensemble des caractères de
conducteur au plus $1$ de $T$ est donc égal à l'ensemble des
caractères de $\mu_{q'+1}$. Soit $(\, ,)_{k'}$ le symbole de Hilbert
du corps $k'$. D'après le théorème 90 de Hilbert, tout élément $g$ de
$T$ s'écrit $g=z(z^{-1})^{\tau}$ avec $z\in\mathcal{O}''^{\times}$. La
formule
\begin{equation}\label{eq0.0.2}
\eta_{0}(g)=(\varpi',zz^{\tau})_{k'}=
\left(\frac{p_{\mathbb{F}_{q'}}(zz^{\tau})}{q'}\right)
\end{equation}
définit un caractère de conducteur 1 de $T$ (ici,
$\left(\frac{\,}{q'}\right)$ désigne le symbole de Legendre relatif au
corps $\Bbb{F}_{q'})$.

Lorsque $k''$ est ramifié sur $k'$, le groupe quotient $T/T_{1}$ est
le groupe à deux éléments. Par suite, tous les caractères de $T$ sont
de conducteur pair à l'exception du caractère non trivial de $T/T_{1}$
qui est de conducteur $1$. Soit $\chi$ un caractère de $T$ de
conducteur $2\lambda\geq 2$ et $j$ un entier tel que
$j<\lambda\leq3j+1$~; alors, il existe $b\in\mathcal{O}'^{\times}$,
uniquement déterminé modulo le sous-groupe multiplicatif
$1+\varpi'^{\lambda-j}\mathcal{O}'$, tel que la restriction de $\chi$
au sous-groupe de congruence $T_{2j+1}$ soit égal au caractère
$\chi_{b,\lambda,j}$ défini par la formule \ref{eq7.1.4} du lemme
\ref{lem7.1.2}.

On note $e$ l'indice de ramification de $k'$ sur $k$ et $\delta$
l'entier tel que l'idéal $\varpi'^{\delta}\mathcal{O}'$ soit la
différente de $k'$ sur $k$. On pose
$\mu=e\lambda_{\psi}-\delta-v''(u)$.

On définit le réseau $B$ en posant 
\begin{equation*}
B=
\begin{cases}
\varpi''^{\mu}\mathcal{O}''\mbox{ si $k''$ est ramifié sur
  $k'$,}\\ \varpi'^{[\frac{\mu+1}{2}]}\mathcal{O}''\mbox{ si $k''$ est
  non ramifié sur $k'$.}
\end{cases}
\end{equation*}
Alors $B$ est un bon réseau $T$-invariant, autodual si et seulement si
$k''$ est ramifié sur $k'$ ou si $\mu$ est pair. Lorsque $k''$ est non
ramifié sur $k'$ et $\mu$ est impair, $B^{*}=\varpi'^{-1}B$ et
$\mathrm{b}^{*}=\mathcal{O}''/\varpi'\mathcal{O}''$ est un
$\mathbb{F}_{q}$ espace vectoriel symplectique de dimension
$\frac{2r}{e}$ (voir le corollaire \ref{co6.2.1}). La donnée du réseau
$B$ détermine une section du revêtement métaplectique au dessus de
$T$, restriction de la section canonique au dessus du sous-groupe
compact maximal $K_{B}$.

Supposons que le réseau $B$ est autodual. Alors, nous montrons que le
caractère trivial intervient dans la représentation de Weil et qu'un
caractère $\chi$ non trivial de $T$ intervient si et seulement si
\begin{itemize}
\item $\chi$ est de conducteur pair, lorsque $k''$ est non ramifié sur
  $k'$,
\item $\chi$ est de conducteur pair $2j>0$ et il existe
  $a\in\mathcal{O}''^{\times}$ tel que
  $\chi_{\vert T_{j}}=\chi_{aa^{\tau},j,[\frac{j}{2}]}$, lorsque $k''$ est
  ramifié sur $k'$.
\end{itemize}

Supposons que le réseau $B$ n'est pas autodual, de sorte que $k''$ est non
ramifié sur $k'$. Nous montrons qu'un caractère $\chi$ de $T$
intervient si et seulement si l'une des deux assertions suivantes est
vérifiée 
\begin{itemize}
\item $\chi$ est de conducteur au plus $1$ et $\chi\neq\eta_{0}$,
\item $\chi$ est de conducteur impair, strictement supérieur à $1$. 
\end{itemize}
A propos de la première assertion, il est à noter que la
représentation de Weil du groupe symplectique $Sp(\mathrm{b}^{*})$
apparaît comme sous-représentation de $K_{B}$ dans la
représentation de Weil de $Sp(W)$ et que les caractères de $T$ de
conducteur au plus $1$ interviennent dans le sous-$K_{B}$-module
correspondant. Or la représentation de Weil de $Sp(\mathrm{b}^{*})$
est de dimension $q'$, tandis que les caractères de $T$ de conducteur
au plus $1$ sont au nombre de $q'+1$.

Dans tous les cas, nous décrivons une base de l'espace des vecteurs
propres de poids $\chi$. Ces derniers résultats sont énoncés dans les
théorèmes \ref{theo7.2.1} et \ref{theo7.3.1} et généralisent ceux
démontrés par Yang dans le cas où $W$ est de dimension $2$ (voir
\cite{yang-1998}).

Nous remercions Paul Broussous, François Courtès et Claude Quitté pour
d'utiles conversations.

\section{Sous-groupes compacts maximaux du groupe symplectique}\label{2}

\subsection{}\label{2.1} 
Dans la suite, $k$ désigne un corps local non archimédien de
caractéristique nulle, $\O$ l'anneau des entiers de $k$ et $\goth p$
son idéal maximal. Nous fixons $\varpi$ une uniformisante de $\O$,
c'est-à-dire $\varpi \in \O$ tel que $\goth p = \varpi\O$.  Le corps
résiduel $\O/\goth p$ est fini de cardinal $q$ et de caractéristique
$p$, nous le notons $\Bbb{F}_{q}$.  L'entier $p$ est appelé la
caractéristique résiduelle de $k$. Nous supposons que $p\neq 2$.

Nous fixons aussi un caractère
unitaire $\psi$ non trivial de $k$.  Son conducteur $\lambda_\psi$ est
l'unique entier relatif vérifiant
\begin{equation*}
\psi_{|\goth p^{\lambda_\psi}}\equiv 1 \mbox{ et } \psi_{|\goth
   p^{\lambda_\psi-1}}\not \equiv 1 .
\end{equation*}

Le caractère $\psi$ induit un caractère $\overline{\psi}$ de
$\Bbb{F}_{q}$ tel que
\begin{equation*}
  \overline{\psi}(p_{\Bbb{F}_{q}}(x))=\psi(\varpi^{\lambda_{\psi}-1}x)\mbox{,
  }x\in\O,
\end{equation*}
$p_{\Bbb{F}_{q}}$ désignant la projection naturelle de  $\O$ sur $\Bbb{F}_{q}$.

\subsection{}\label{2.2} 
Dans la suite, $F$ désigne soit un anneau local fini de
caractéristique différente de $2$, muni de la topologie discrète, soit
un corps local non archimédien de caractéristique résiduelle
différente de $2$, muni de sa topologie localement compacte.  On appelle
$F$-espace symplectique tout couple $(W , \beta)$ où $ W$ est un
$F$-module de type fini et $\beta$ une forme bilinéaire alternée non
dégénérée sur $W$. 

Tout $F$-module de type fini est naturellement muni d'une topologie
localement compacte : lorsque $F$ est un anneau local fini, il s'agit
de la topologie discrète et, lorsque $F$ est un corps local, il s'agit
de sa topologie d'espace vectoriel de dimension finie sur ce corps
local. On munit alors $GL(W)$, le groupe linéaire de $W$, de la
topologie induite par celle de $End_{F}(W)$, qui en fait un groupe
topologique localement compact.

On désigne par $J_{r}$ la matrice antisymétrique d'ordre $2r$~:
\begin{equation*}
J_{r}  =  \left(
\begin{array}{ll}
 \,  \, \,  \, \,0& \, \, \,  I_r \\
   \,    -I_r & \, \, \,  0
\end{array}
\right),
\end{equation*}
$I_r$ étant la matrice identité d'ordre $r$.

Si $W$ est un $F$-module libre, on appelle base symplectique
(resp. presque symplectique) de $W$, toute base $(e_1, \ldots, e_r, f_1,
\ldots, f_r)$ telle que la matrice de $\beta$ dans cette base soit
$J_{r}$ (resp. $tJ_{r}$, où $t\in F$ est non nul). Dans ce cas, de
telles bases existent (voir \cite{klingenberg-1963}).

On note $Sp(W,\beta)$ ou plus simplement $Sp(W)$ le groupe
symplectique associé à $(W,\beta)$~: 
\begin{equation*}
Sp(W) = \{g \in GL(W)\, / \, \beta(gv, gw)= \beta(v, w),\, \mbox{
  pour tout } v , w \in W\}.
\end{equation*}
C'est un sous-groupe fermé de $GL(W)$. Il contient, lorsque $W$ est
non nul, le groupe à deux éléments $\{\pm Id\}$ comme sous-groupe
central~; lorsque $F$ est un corps, ce sous-groupe est le centre de
$Sp(W)$. Lorsque $W$ est nul, $Sp(W)$ est le groupe trivial.

Soit $g \in GL(W)$ dont la matrice dans une base presque symplectique
s'écrit par blocs $g=\left( 
\begin{smallmatrix}
a & b\\
c & d
\end{smallmatrix}\right)
.$ 
Alors, on a~: $g \in Sp(W)$ si et seulement si l'une des
conditions suivantes est satisfaite~:
\begin{equation*}
 {}^ta\, c = {}^tc\, a , \, \, {}^tb\, d = {}^td\, b \quad \mbox{ et }
 \quad {}^t a\, d - {}^tc\, b =  I_r,
\end{equation*}
\begin{equation*}
 a \,  {}^tb = b \,   {}^ta, \, \, c \,  {}^td = d \,  {}^tc \quad \mbox{ et } \quad a \,  {}^t d - b \,  {}^t c  = I_r.
\end{equation*}

Si $U$ est une partie de $W$, on note $U^{\perp}$ son orthogonal
relativement à $\beta$. Un sous-module $U$ de $W$ est dit totalement
isotrope, si $U\subset U^{\perp}$. Si $U=U^{\perp}$, on dit que c'est
un sous-espace lagrangien ou un lagrangien de $W$.

Lorsque $F$ est un corps, un sous-espace vectoriel totalement isotrope
est de dimension maximale si et seulement si c'est un lagrangien de
$W$. Dans cette situation, le groupe symplectique $Sp(W)$ opère
transitivement sur les lagrangiens de $(W, \beta)$.

\subsection{}\label{2.3}
Soit $(W , \beta)$ un $k$-espace symplectique.
\begin{defi}\label{def2.3.1}
Une base $(e_1, \ldots, e_r, f_1, \ldots, f_r)$ de $W$ est dite
autoduale relativement à $\beta$ et $\psi$ si elle vérifie les
relations~:
\begin{equation}\label{eq2.3.1}
\beta(e_i, e_j) = \beta(f_i, f_j) = 0 \mbox{ et } \beta(e_i,f_j)
=\varpi^{\lambda_\psi} \delta_{ij}, \, \,  1\leq i,
j \leq r 
\end{equation}
ou, de manière équivalente, si la matrice de $\beta$ dans cette base
est $\varpi^{\lambda_{\psi}}J_{r}$.
\end{defi}
 
Un réseau $B$ de $W$ est un $\mathcal O $-module, compact et ouvert.
On note
\begin{equation}\label{eq2.3.2}
\begin{split}
B^* &= \{v \in W , \quad \beta(v,b) \in  \varpi^{\lambda_\psi} \O,
\mbox{ pour tout } b \in B \}\\
&= \{v \in W , \quad \psi(\beta(v,b))=1,
\mbox{ pour tout } b \in B \}.
\end{split}
\end{equation}

Alors, $B^*$ est aussi un réseau de $W$ ; on l'appelle le réseau dual
de $B$ relativement à $\beta$ et $\psi$.

\begin{defi}\label{def2.3.2} 
On dit que $B$ est un bon réseau si 
\begin{equation*}
\varpi B^* \subset B \subset B^*.
\end{equation*}
Il est dit autodual si $ B = B^*$.
\end{defi}

Soit $B\subset W$ un bon réseau. On note ${\mathrm b}^*$ le quotient $ B^*
/B$. Il est naturellement muni d'une structure d'espace vectoriel sur
$\Bbb{F}_{q}$ et de la forme symplectique $\beta_{{\mathrm{b}
}^*}$ définie par~:
\begin{equation}\label{eq2.3.3}
  \beta_{{\mathrm{b}
}^*}(p_{{\mathrm{b}
}^{*}}(w),p_{{\mathrm{b}
}^{*}}(w')) =
p_{\Bbb{F}_{q}}(\varpi^{1-\lambda_{\psi}}\beta(w, w')),
\end{equation}
où $p_{{\mathrm{b}
}^{*}}$ désigne la projection naturelle de $B^{*}$ sur ${\mathrm{b}
}^{*}$.
On pose $l(B)=\frac{1}{2}\dim_{\Bbb{F}_{q}}b^{*}$. Alors $l(B)$ est un
entier et deux bons réseaux $B$ et $B'$ sont conjugués sous l'action
de $Sp(W)$ si et seulement si on a $l(B)=l(B')$.

Si $B$ est un bon réseau, on désigne par $K_{B}$ son stabilisateur
dans $Sp(W)$. Alors, les $K_{B}$, pour $B$ un bon réseau, forment
l'ensemble des sous-groupes compacts maximaux de $Sp(W)$ et deux tels
sous-groupes $K_{B}$ et $K_{B'}$ sont conjugués dans $Sp(W)$ si et
seulement si $B$ et $B'$ sont dans la même $Sp(W)$-orbite, autrement
dit si et seulement si $l(B)=l(B')$.

Soit $B\subset W$ un bon réseau et $l=l(B)$. Alors, il existe une base
autoduale $(e_1, \ldots, e_r,$ $f_1, \ldots, f_r)$ de $W$ telle que
\begin{equation}\label{eq2.3.4}
B = \varpi\O e_1 \oplus \cdots \oplus \varpi\O e_l\oplus\O e_{l+1}
\oplus \cdots \oplus \O e_r \oplus \O f_1 \oplus \cdots \oplus \O f_r.
\end{equation}

De plus, les éléments de $K_{B}$ sont les éléments $g$ de $Sp(W)$
dont la matrice dans la base $(e_1, \ldots,e_l, e_{l+1}, \ldots, e_r,
f_1, \ldots,f_l, f_{l+1}, \ldots, f_r)$ s'écrit par blocs~:
\begin{equation}\label{eq2.3.5}  
\begin{pmatrix}
a_{11}&\varpi a_{12}&\varpi b_{11}&\varpi b_{12}\\
a_{21}& a_{22}&\varpi b_{21}& b_{22}\\
\varpi^{-1}c_{11}&c_{12}&d_{11}&d_{12}\\
c_{21}&c_{22}&\varpi d_{21}&d_{22}
\end{pmatrix}, 
\end{equation}
les matrices $a_{i j }, b_{i j }, c_{i j }, d_{i j }$ étant à
coefficients dans $\O$, $a_{11}$ et $d_{11}$ (resp. $a_{22}$ et
$d_{22}$) étant carrées d'ordre $l$ (resp. $r-l$).

\subsection{}\label{2.4}
On garde les notations du paragraphe précédent. Soit $B\subset W$ un
bon réseau, $l=l(B)$ et $(e_{1},\ldots,e_{r},f_{1},\ldots,f_{r})$ une
base autoduale de $W$ vérifiant la relation \ref{eq2.3.4}. Alors, le
réseau dual de $B$ est
\begin{equation*}
B^* =\O e_1 \oplus \cdots \oplus \O e_r \oplus \varpi^{-1} \O f_1
\oplus \cdots \oplus \varpi^{-1}\O f_l \oplus \O f_{l+1} \oplus \cdots
\oplus \O f_r
\end{equation*}
et l'on a 
\begin{equation*}
B\subset A\subset B^{*},
\end{equation*}
où $A=\O e_{1}\oplus\cdots\oplus\O e_{r}\oplus\O
f_{1}\oplus\cdots\oplus\O f_{r}$ est un réseau autodual. On pose
$K=K_{A}$. Les éléments de $K$ sont les éléments de $Sp(W)$ dont la
matrice dans la base $(e_{1},\ldots,e_{r},f_{1},\ldots,f_{r})$ est à
coefficients dans $\O$.

Nous avons vu que ${\mathrm{b}
}^{*}=B^{*}/B$, muni de la forme
$\beta_{{\mathrm{b}
}^{*}}$ est un espace symplectique de dimension $2l$ sur
$\Bbb{F}_{q}$. On vérifie que la famille de vecteurs
\begin{equation*}
(p_{{\mathrm{b}
}^{*}}(e_1), \ldots, p_{{\mathrm{b}
}^{*}}(e_l),
p_{{\mathrm{b}
}^{*}}(\varpi^{-1} f_1), \ldots, p_{{\mathrm{b}
}^{*}}(\varpi^{-1} f_l))
\end{equation*}
est une base symplectique de $({\mathrm b}^*, \beta_{{\mathrm{b}
}^*})$.

De même, $\mathrm{b}=B/\varpi B^{*}$ est naturellement muni d'une
structure de $\Bbb{F}_{q}$-espace vectoriel et de la forme
symplectique $\beta_{\mathrm{b}}$ définie par
\begin{equation*}
  \beta_{\mathrm{b}}(p_{\mathrm{b}}(w),p_{\mathrm{b}}(w'))=
p_{\Bbb{F}_{q}}(\varpi^{-\lambda_{\psi}}\beta(w,w'))\mbox{, }w,w'\in B,
\end{equation*}
où $p_{\mathrm{b}}$ désigne la projection naturelle de $B$ sur
$\mathrm{b}$. La famille de vecteurs 
\begin{equation*}
  (p_{\mathrm{b}}(e_{l+1}),\ldots,p_{\mathrm{b}}(e_{r}),
p_{\mathrm{b}}(f_{l+1}),\ldots,p_{\mathrm{b}}(f_{r}))
\end{equation*}
est une base symplectique de $\mathrm{b}$ de sorte que la dimension de
$\mathrm{b}$ sur $\Bbb{F}_{q}$ est $2(r-l)$.

Le groupe $K_{B}$ laissant invariant $B$, laisse invariant $B^{*}$. Il
agit donc naturellement dans ${\mathrm{b}
}^{*}$ (resp. $\mathrm{b}$) et cette
action induit un morphisme de $K_{B}$ dans $Sp({\mathrm{b}
}^{*})$ (resp.
$Sp(\mathrm{b})$), le groupe symplectique associé à
$(\mathrm{b^{*}},\beta_{{\mathrm{b}
}^{*}})$ (resp. $(\mathrm{b},\beta_{\mathrm{b}})$)~; ce
morphisme est noté $p_{Sp(\mathrm{b^{*}})}$ (resp. $p_{Sp(\mathrm{b})}$). On a
le résultat suivant~:
\begin{lem}\label{lem2.4.1}
Le morphisme $p_{Sp({\mathrm{b}
}^{*})}\times p_{Sp(\mathrm{b})}$ est surjectif de
$K_{B}$ sur le groupe $Sp({\mathrm{b}
}^{*})\times Sp(\mathrm{b})$. 
\end{lem}
\begin{dem}
Rappelons que si $(W,\beta)$ est un espace symplectique sur le corps
$F$, le groupe symplectique $Sp(W)$ est engendré par les transvections
symplectiques de $(W,\beta)$ qui sont les transformations de $W$ de la
forme $\tau_{a,v}:w\mapsto w+a\beta(v,w)v$, $a\in F$, $v\in W$. Il est
immédiat que les transvections $\tau_{a,v}$ telles que
$a\in\varpi^{1-\lambda_{\psi}}\O$ et $v\in B^{*}$
(resp. $a\in\varpi^{-\lambda_{\psi}}\O$ et $v\in B$) sont dans $K_{B}$
et que leur image par $p_{Sp({\mathrm{b}
}^{*})}\times p_{Sp(\mathrm{b})}$
parcourt l'ensemble des éléments de la forme $(\tau,Id)$
(resp. $(Id,\tau)$) où $\tau$ est une transvection symplectique de
$({\mathrm{b}
}^{*},\beta_{{\mathrm{b}
}^{*}})$ (resp. $(\mathrm{b},\beta_{\mathrm{b}})$).
D'où le lemme.
\end{dem}

On désigne par $K'_{B}$ le noyau du morphisme
$p_{Sp({\mathrm{b}
}^{*})}\times p_{Sp(\mathrm{b})}$. On a 
\begin{equation*}
K'_{B}=\{g\in K_{B}\mbox{, } (g-1)B\subset\varpi B^{*}\mbox{ et
}(g-1)B^{*}\subset B\}. 
\end{equation*}
Les éléments de $K'_{B}$ sont les éléments de
$Sp(W)$ dont la matrice dans la base
$(e_{l},\ldots,e_{l},$ $e_{l+1},\ldots,e_{r},
f_{1},\ldots,f_{l},f_{l+1},\ldots,f_{r})$ s'écrit par blocs
\begin{equation*}
\begin{pmatrix}
I_{l}+\varpi a_{11}&\varpi a_{12} &\varpi^{2}b_{11}&\varpi b_{12}\\
a_{21}&I_{r-l}+\varpi a_{22}&\varpi b_{21}&\varpi b_{22}\\
c_{11}&c_{12}&I_{l}+\varpi d_{11}&d_{12}\\
c_{21}&\varpi c_{22}&\varpi d_{21}&I_{r-l}+\varpi d_{22}
\end{pmatrix}
\end{equation*}
les matrice $a_{ij}$, $b_{ij}$, $c_{ij}$ et $d_{ij}$ étant des matrices à
coefficients dans $\O$, $a_{11}$ et $d_{11}$ (resp. $a_{22}$ et
$d_{22}$) étant carrées d'ordre $l$ (resp. $r-l$).

\subsection{}\label{2.5}
On garde les notations du paragraphe précédent. Soit $n$ un entier
naturel. On désigne par $O_{n}$ l'anneau $\O/\varpi^{n+1}\O$ et par
$p_{O_{n}}$ la projection canonique de $\O$ sur $O_{n}$. Alors,
$O_{0}=\Bbb{F}_{q}$ et si $n\geq 1$, $O_{n}$ est un anneau local
principal fini dont les idéaux propres sont les puissances $k$-ièmes
de l'idéal maximal $\varpi\O/\varpi^{n+1}\O$, $1\leq k\leq n$. De plus
la projection $p_{O_{m}}$ passe au quotient, pour tout entier $m>n$,
en un morphisme d'anneau $p_{n,m}:O_{m}\longrightarrow O_{n}$~; on
note plus simplement $p_{n}=p_{n,n+1}$. La famille des $O_{n}$ munie
des morphismes $p_{n,m}$ constitue un système projectif d'anneaux dont
il est bien connu que l'anneau $\O$, muni de sa topologie, est la
limite projective~: $\O=\underset{\longleftarrow}{\lim}O_{n}$.

On pose ${\mathrm{b}
}_{n}=B/\varpi^{n+1}B^{*}$ et
${\mathrm{b}
}^{*}_{n}=B^{*}/\varpi^{n}B$. Ce sont des $O_{n}$-modules (et
donc des $O_{m}$-modules, pour $m\geq n$) de type fini. On a
${\mathrm{b}
}_{0}=\mathrm{b}$ et ${\mathrm{b}
}^{*}_{0}={\mathrm{b}
}^{*}$. Lorsque $n\geq1$,
${\mathrm{b}
}_{n}$ et ${\mathrm{b}
}^{*}_{n}$ admettent $2r$ générateurs et ils
sont libres si et seulement si $l=0$ ou $l=r$ (dans ce cas, on a
${\mathrm{b}
}_{n}={\mathrm{b}
}^{*}_{n+1}$ ou ${\mathrm{b}
}^{*}_{n}={\mathrm{b}
}_{n+1}$). On
désigne par $p_{{\mathrm{b}
}_{n}}$ (resp. $p_{{\mathrm{b}
}^{*}_{n}}$) la
projection canonique de $B$ (resp. $B^{*}$) sur ${\mathrm{b}
}_{n}$
(resp. ${\mathrm{b}
}^{*}_{n}$)~; elle induit une projection naturelle de
${\mathrm{b}
}_{n+1}$ (resp. ${\mathrm{b}
}^{*}_{n+1}$) sur ${\mathrm{b}
}_{n}$
(resp. ${\mathrm{b}
}^{*}_{n}$) notée $p_{n}$, laquelle est compatible avec
les structures de $O_{n+1}$-modules.

On désigne par $End_{B}(W)$ la sous-algèbre de $End_{k}(W)$ constituée
des endomorphismes linéaires laissant invariant $B$ et $B^{*}$ et par
$GL_{B}(W)$ le sous-groupe de $GL(W)$ constitué des éléments qui sont,
ainsi que leur inverse, dans $End_{B}(W)$. Tout élément $g$ de
$End_{B}(W)$ définit un élément $\dot{g}_{n}$
(resp. $\dot{g}^{*}_{n}$) de $End_{O_{n}}({\mathrm{b} }_{n})$
(resp. $End_{O_{n}}({\mathrm{b} }^{*}_{n})$) en posant
$\dot{g}_{n}p_{{\mathrm{b} }_{n}}(v)=p_{{\mathrm{b} }_{n}}(gv)$, $v\in
B$ (resp. $\dot{g}^{*}_{n}p_{{\mathrm{b} }^{*}_{n}}(v)=p_{{\mathrm{b}
  }^{*}_{n}}(gv)$, $v\in B^{*}$). On vérifie que l'application
$g\mapsto\dot{g}_{n}$ (resp. $g\mapsto\dot{g}^{*}_{n}$) est un
morphisme surjectif d'algèbres de $End_{B}(W)$ sur
$End_{O_{n}}({\mathrm{b} }_{n})$ (resp. $End_{O_{n}}({\mathrm{b}
}^{*}_{n})$) de noyau $I_{n}=\{g\in End_{k}(W) \vert
gB\subset\varpi^{n+1}B^{*}\}$ (resp. $I^{*}_{n}=\{g\in End_{k}(W)
\vert gB^{*}\subset\varpi^{n}B\}$). De plus, pour tout entier $m>n$,
le morphisme $g\mapsto\dot{g}_{m}$ (resp. $g\mapsto\dot{g}^{*}_{m}$)
passe au quotient en un morphisme d'algèbres
$r_{n,m}:End_{O_{m}}({\mathrm{b} }_{m})\longrightarrow
End_{O_{n}}({\mathrm{b} }_{n})$
(resp. $r^{*}_{n,m}:End_{O_{m}}({\mathrm{b} }^{*}_{m})\longrightarrow
End_{O_{n}}({\mathrm{b} }^{*}_{n})$)~; on note plus simplement
$r_{n}=r_{n,n+1}$ (resp. $r^{*}_{n}=r^{*}_{n,n+1}$).

\begin{lem}\label{lem2.5.1} Soit $n$ un entier naturel.

(i) L'application $g\mapsto\dot{g}_{n}$
  (resp. $g\mapsto\dot{g}^{*}_{n}$) définit un morphisme surjectif de
  groupes de $GL_{B}(W)$ sur le groupe linéaire $GL_{O_{n}}({\mathrm{b}
  }_{n})$ (resp. $GL_{O_{n}}({\mathrm{b} }^{*}_{n})$), de noyau $Id+I_{n}$
  (resp. $Id+I^{*}_{n}$).

(ii) Pour tout entier $m>n$, l'application $r_{n,m}$
  (resp. $r^{*}_{n,m}$) induit un morphisme surjectif de groupes de
  $GL_{O_{m}}({\mathrm{b} }_{m})$ sur $GL_{O_{n}}({\mathrm{b} }_{n})$
  (resp. $GL_{O_{m}}({\mathrm{b} }^{*}_{m})$ sur $GL_{O_{n}}({\mathrm{b}
  }^{*}_{n})$) de noyau $Id+I_{n}/I_{m}$
  (resp. $Id+I^{*}_{n}/I^{*}_{m}$).
\end{lem}
\begin{dem}
L'assertion (ii) est une conséquente immédiate de l'assertion
(i). Montrons l'assertion (i).

Il est clair que l'application $g\mapsto\dot{g}_{n}$
(resp. $g\mapsto\dot{g}^{*}_{n}$) est un morphisme de groupes de
$GL_{B}(W)$ dans $GL_{O_{n}}({\mathrm{b}
}_{n})$
(resp. $GL_{O_{n}}({\mathrm{b}
}^{*}_{n})$), de noyau $Id+I_{n}$
(resp. $Id+I^{*}_{n}$). Il reste à montrer la surjectivité.

Supposons que $n=0$. Soient $e_{1},\ldots,e_{r},f_{1},\ldots,f_{r}$
une base autoduale de $W$ telle que $B=\O\varpi
e_{1}\oplus\cdots\oplus\O\varpi e_{l}\oplus\O
e_{l+1}\oplus\cdots\oplus\O e_{r}\oplus\O f_{1}\oplus\cdots\oplus\O
f_{r} $ et $U$ (resp. $V$) le sous-espace vectoriel de $W$ engendré
par $e_{1},\ldots,e_{l},f_{1},\ldots,f_{l}$
(resp. $e_{l+1},\ldots,e_{r},f_{l+1},\ldots,f_{r}$). Alors $G=\{g\in
GL(W)\vert g_{\vert U}=Id\mbox{ et } g(B\cap V)=B\cap V\}$
(resp. $G^{*}=\{g\in GL(W)\vert g(B^{*}\cap U)=B^{*}\cap U\mbox{ et }
g_{\vert V}=Id\}$) est un sous-groupe de $GL_{B}(W)$ et son image par
l'application $g\mapsto\dot{g}_{0}$ (resp. $g\mapsto\dot{g}^{*}_{0}$)
est $GL_{\Bbb{F}_{q}}(\mathrm{b})$
(resp. $GL_{\Bbb{F}_{q}}({\mathrm{b} }^{*})$).

Supposons $n\geq1$. Soit $g\in End_{B}(W)$. Dire que $\dot{g}_{n}\in
GL_{O_{n}}({\mathrm{b} }_{n})$ (resp. $\dot{g}^{*}_{n}\in
GL_{O_{n}}({\mathrm{b} }^{*}_{n})$) revient à dire qu'il existe $g'\in
End_{B}(W)$ tel que $gg'\in Id+I_{n}$ (resp. $gg'\in
Id+I^{*}_{n}$). Or, 
 $Id+I_{n}$ et $Id+I^{*}_{n}$ sont des sous-groupes de
$GL_{B}(W)$. Il est alors clair que $g\in GL_{B}(W)$.
\end{dem}
Le $O_{n}$-module ${\mathrm{b} }_{n}$ (resp. ${\mathrm{b} }^{*}_{n}$) est
naturellement muni d'une forme symplectique, $\beta_{{\mathrm{b} }_{n}}$
(resp. $\beta_{{\mathrm{b} }^{*}_{n}}$) définie par
\begin{eqnarray}
\beta_{{\mathrm{b} }_{n}}(p_{{\mathrm{b} }_{n}}(v),p_{{\mathrm{b}
  }_{n}}(w)) &=& p_{O_{n}}(\varpi^{-\lambda_{\psi}}\beta(v,w))\mbox{,
} v,w\in B\label{eq2.5.1}\\ 
\beta_{{\mathrm{b} }^{*}_{n}}(p_{{\mathrm{b}
  }^{*}_{n}}(v),p_{{\mathrm{b} }^{*}_{n}}(w)) &=&
p_{O_{n}}(\varpi^{1-\lambda_{\psi}}\beta(v,w))\mbox{, }v,w\in B^{*}\label{eq2.5.2}
\end{eqnarray}
On désigne par $Sp({\mathrm{b}}_{n})$
(resp. $Sp({\mathrm{b}}^{*}_{n})$) le groupe symplectique correspondant. Il
est clair que, si $g\in K_{B}$, $\dot{g}_{n}\in Sp({\mathrm{b}}_{n})$
(resp. $\dot{g}^{*}_{n}\in Sp({\mathrm{b} }^{*}_{n})$) et que, pour tout
entier $m>n$, le morphisme $r_{n,m}$ (resp. $r^{*}_{n,m}$) envoie
$Sp({\mathrm{b} }_{m})$ dans $Sp({\mathrm{b} }_{n})$ (resp. $Sp({\mathrm{b}
}^{*}_{m})$ dans $Sp({\mathrm{b} }^{*}_{n})$). La famille des $Sp({\mathrm{b}
}_{n})$ (resp. $Sp({\mathrm{b} }^{*}_{n})$), munie des morphismes
$r_{n,m}$ (resp. $r^{*}_{n,m}$) constitue un système projectif de
groupes.
\begin{lem}\label{lem2.5.2}
(i) Pour tout  $n\in\mathbb{N}$, le morphisme de groupes
  $r_{n}:Sp({\mathrm{b}
}_{n+1})\longrightarrow Sp({\mathrm{b}
}_{n})$
  (resp. $r^{*}_{n}:Sp({\mathrm{b}
}^{*}_{n+1})\longrightarrow
  Sp({\mathrm{b}
}^{*}_{n})$) est surjectif.

(ii) L'application $g\longrightarrow(\dot{g}_{n})_{n\in\Bbb{N}}$
(resp. $g\longrightarrow(\dot{g}^{*}_{n})_{n\in\Bbb{N}}$) induit un
isomorphisme de groupes topologiques de $K_{B}$ sur
$\underset{\longleftarrow}{\lim}\,Sp({\mathrm{b}
}_{n})$
(resp. $\underset{\longleftarrow}{\lim}\,Sp({\mathrm{b}
}^{*}_{n})$).
\end{lem}
\begin{dem}
(i) On montre la surjectivité de $r_{n}$, la démonstration de celle de
  $r^{*}_{n}$ étant identique. 

Supposons dans un premier temps que $n\geq1$. Soit $g\in Sp({\mathrm{b}
}_{n})$. D'après le lemme \ref{lem2.5.1}, il existe $\tilde{g}\in
GL_{O_{n+1}}({\mathrm{b} }_{n+1})$ tel que $r_{n}(\tilde{g})=g$. Il suffit
de montrer qu'il existe $h\in I_{n}/I_{n+1}$ tel que
$(Id+h)\tilde{g}\in Sp({\mathrm{b} }_{n+1})$. Définissons l'application
$\gamma:b_{n+1}\times b_{n+1}\longrightarrow O_{n+1}$ en posant
$\gamma(v,w)= 
\beta_{{\mathrm{b}}_{n+1}}(\tilde{g}^{-1}v,\tilde{g}^{-1}w)-
\beta_{{\mathrm{b}}_{n+1}}(v,w)$~; c'est une application bilinéaire alternée sur
${\mathrm{b} }_{n+1}$ à valeurs dans l'idéal minimal
$\varpi^{n+1}\O/\varpi^{n+2}\O$. Soit $h\in I_{n}/I_{n+1}$. Dire que
$(Id+h)\tilde{g}\in Sp({\mathrm{b}}_{n+1})$ revient à dire que l'on a
\begin{equation}\label{eq2.5.3}
\beta_{{\mathrm{b}}_{n+1}}((Id+h)v,(Id+h)w)-\beta_{{\mathrm{b}}_{n+1}}(v,w)
=\gamma(v,w)\mbox{, }v,w\in{\mathrm{b}}_{n+1}
\end{equation}

Mais, $h$ étant un élément de $ I_{n}/I_{n+1}$, on a
$hb_{n+1}\subset\varpi^{n+1}B^{*}/\varpi^{n+2}B^{*}$. Comme $n\geq1$,
on en déduit que
\begin{eqnarray*}
\beta_{{\mathrm{b} }_{n+1}}(hb_{n+1},hb_{n+1})&\subset&
p_{O_{n+1}}(\varpi^{2(n+1)-\lambda_{\psi}}\beta(B^{*},B^{*}))\\ 
&\subset&p_{O_{n+1}}(\varpi^{2n+1}\O) \subset
p_{O_{n+1}}(\varpi^{n+2}\O)=\{0\}.
\end{eqnarray*}
La relation \ref{eq2.5.3} s'écrit alors
\begin{equation*}
\beta_{{\mathrm{b} }_{n+1}}(hv,w)+\beta_{{\mathrm{b}}_{n+1}}(v,hw)
=\gamma(v,w)\mbox{, }v,w\in{\mathrm{b} }_{n+1}.
\end{equation*}
Or, la forme $\beta_{{\mathrm{b}}_{n+1}}$ étant symplectique, il existe
$h'\in End_{O_{n+1}}({\mathrm{b}}_{n+1})$ tel que
\begin{equation*}
\gamma(v,w)=\beta_{{\mathrm{b}}_{n+1}}(h'v,w)=
\beta_{{\mathrm{b}}_{n+1}}(v,h'w)\mbox{, }v,w\in{\mathrm{b}}_{n+1}.
\end{equation*}
Il suffira de prendre $h=\frac{1}{2}h'$, dès que l'on aura montré que
$h'\in I_{n}/I_{n+1}$. Or, on a
$\gamma(v,w)\in\varpi^{n+1}\O/\varpi^{n+2}\O$, $v,w\in b_{n+1}$. On en
déduit que, pour tout $w\in\varpi B/\varpi^{n+2}B^{*}$,
$\beta_{{\mathrm{b}}_{n+1}}(h'{\mathrm{b}}_{n+1},w)=\{0\}$. Autrement dit,
$h'{\mathrm{b}}_{n+1}\subset(\varpi B/\varpi^{n+1}B^{*})^{\perp}=\varpi^{n+1}
B^{*}/\varpi^{n+2}B^{*}$, i.e. $h'\in I_{n}/I_{n+1}$ comme voulu.

Il reste à examiner le cas $n=0$. Mais, il suit du lemme
\ref{lem2.4.1} que l'application $g\mapsto\dot{g}_{0}$ est un
morphisme surjectif de $K_{B}$ sur $Sp(\mathrm{b})=Sp({\mathrm{b}
}_{0})$. Il
suffit alors de remarquer que $\dot{g}_{0}=r_{0}(\dot{g}_{1})$, $g\in
K_{B}$.

(ii) On donne la démonstration pour l'application
$g\longrightarrow(\dot{g}_{n})_{n\in\Bbb{N}}$, celle pour l'autre
application étant identique. Pour tout entier naturel $n$, la matrice
d'un élément de $I_{n}$ dans une base quelconque du $\O$-module $B$
est à coefficients dans $\varpi^{n}\O$. On en déduit que
$\cap_{n\in\Bbb{N}}I_{n}=\{0\}$ et que la famille
$(Id+I_{n})_{n\in\Bbb{N}}$ est une base de voisinages de l'élément
neutre dans $GL_{B}(W)$. Il est alors immédiat que l'application
$g\longrightarrow(\dot{g}_{n})_{n\in\Bbb{N}}$ est un morphisme
injectif de groupes topologiques de $K_{B}$ dans
$\underset{\longleftarrow}{\lim}\,Sp({\mathrm{b}
}_{n})$. Comme $K_{B}$ est
compact, il nous suffit de montrer que ce morphisme est surjectif.

Soit donc
$(g_{n})_{n\in\Bbb{N}}\in\underset{\longleftarrow}{\lim}\,Sp({\mathrm{b}
}_{n})$. Pour tout entier naturel $n$, soit $\tilde{g}_{n}\in
GL_{B}(W)$ tel que $\dot{\tilde{g}}_{n}=g_{n}$. Par construction, on a
$\tilde{g}_{n+1}-\tilde{g}_{n}\in I_{n}$. Il s'ensuit que la suite
$(\tilde{g}_{n})$ converge dans $End_{B}(W)$ vers un élément $g$
appartenant à $GL_{B}(W)$ et vérifiant $\dot{g}_{n}=g_{n}$,
$n\in\Bbb{N}$. Il reste à voir que $g\in Sp(W)$. Soit $v,w\in
B$. Alors, pour tout $n\in\Bbb{N}$, on a
$p_{O_{n}}(\varpi^{-\lambda_{\psi}}(\beta(gv,gw)-\beta(v,w)))=
\beta_{{\mathrm{b}
  }_{n}}(\dot{g}_{n}p_{{\mathrm{b}}_{n}}(v),\dot{g}_{n}p_{{\mathrm{b}}_{n}}(w))
-\beta_{{\mathrm{b}}_{n}}(p_{{\mathrm{b} }_{n}}(v),p_{{\mathrm{b}
  }_{n}}(w))=0$. On en déduit que
$\beta(gv,gw)-\beta(v,w)\in
\cap_{n\in\Bbb{N}}\varpi^{-\lambda_{\psi}+n}\O=\{0\}$. D'où
le lemme.
\end{dem}

\subsection{}\label{2.6}
On reprend les notations du paragraphe \ref{2.4}. On note ${\mathrm x} =
A/B$~; c'est un sous-espace totalement isotrope maximal de ${\mathrm
  b}^*$. On voit que ${\mathrm x}$ est le sous-espace vectoriel de ${\mathrm
  b}^*$ engendré par $(p_{{\mathrm{b}
}^{*}}(e_1), \ldots,
p_{{\mathrm{b}
}^{*}}(e_l))$.  Inversement~: si ${\mathrm y}$ est un lagrangien
de ${\mathrm b}^*$, alors $A=p_{{\mathrm{b}
}^{*}}^{-1}( {\mathrm y}) $ est un réseau
autodual de $W$ tel que $B\subset A\subset B^{*}$ et il existe une
base autoduale $(e_1, \ldots, e_r, f_1, \ldots, f_r)$ de $W$ vérifiant
la relation \ref{eq2.3.4} et engendrant $A$ comme $\O$-module.


%
D'autre part, $Sp({\mathrm b}^*)$ est engendré par $J_{l}$ et le
sous-groupe parabolique $P_{\mathrm{b}^{*}}= Sp(\mathrm{b}^{*})(\mathrm{x})$,
stabilisateur du lagrangien $\mathrm{x}$ dans $Sp({\mathrm b}^*)$.




\noindent Si on note $P_B = p_{Sp({\mathrm b}^*)}^{-1}(P_{\mathrm{b}^{*}})$, on a~:
\begin{equation*}
K_B' \subset K \mbox{ et } K \cap K_B = P_B.
\end{equation*}

On voit alors que $P_{B}$ est l'ensemble des
$g\in Sp(W)$ dont la matrice dans la base
\begin{equation*}
(e_1, \ldots,e_l, e_{l+1}, \ldots, e_r, f_1, \ldots,f_l, f_{l+1},
  \ldots, f_r)
\end{equation*}
s'écrit par blocs
\begin{equation}\label{eq2.4.1}
g=
\begin{pmatrix}
a_{11}&\varpi a_{12}&\varpi b_{11} &\varpi b_{12}\\
a_{21}& a_{22}& \varpi b_{21} & b_{22}\\
c_{11} & c_{12} & d_{11} & d_{12}\\
c_{21} & c_{22} & \varpi d_{21} & d_{22}
\end{pmatrix}
\end{equation}
les matrices $a_{ij}$, $b_{ij}$, $c_{ij}$ et $d_{ij}$ étant à
coefficients dans $\O$, $a_{11}$ et $d_{11}$ (resp. $a_{22}$ et
$d_{22}$) étant carrées d'ordre $l$ (resp. $r-l$).

On désigne par $\varsigma_B $ l'élément de $K_B$ dont la matrice dans
cette même base est
\begin{equation*}
\varsigma_B =\begin{pmatrix}
0&0&\varpi I_{l}&0\\
0&I_{r-l}&0&0\\
-\varpi^{-1}I_{l}&0&0&0\\
0&0&0&I_{r-l}
\end{pmatrix}
\end{equation*}

\begin{lem}\label{lem2.6.1} Avec les notations ci-dessus, $K_B$ est engendré par
$P_{B}$ et $\varsigma_B$.
\end{lem}
\begin{dem}
On a~:
\begin{itemize}
	\item $\varsigma_B(e_i)= -\varpi^{-1}f_i, \quad 1 \leq i \leq l$,
	\item $\varsigma_B(\varpi^{-1}f_i)= e_i, \quad 1 \leq i \leq l$.
\end{itemize}
Il s'en suit que $ p_{Sp({\mathrm b}^*)}(\varsigma_B)= J_{l}$.  Le
résultat voulu s'en déduit facilement.
\end{dem}

On désigne par $N_{B}$ (resp. $\overline{N}_{B}$) le sous-groupe de
$K_{B}$ constitué des éléments $x_{B}(a)$ (resp. $y_{B}(a)$), ayant
dans la base
\begin{equation*}
(e_1, \ldots,e_l, e_{l+1}, \ldots, e_r, f_1, \ldots,f_l, f_{l+1},
  \ldots, f_r)
\end{equation*}
une matrice de la forme
\begin{equation*}
\begin{pmatrix}
I_{l} & 0 & \varpi a & 0\\
0 & I_{r-l}& 0 & 0\\
0 & 0 & I_{l} & 0\\
0 & 0 & 0 & I_{r-l}
\end{pmatrix}
\qquad\mbox{ (resp.} 
\begin{pmatrix}
I_{l} & 0 & 0& 0\\
0 & I_{r-l}& 0 & 0\\
\varpi^{-1}a & 0 & I_{l}& 0\\
0 & 0 & 0 & I_{r-l}
\end{pmatrix}
\mbox{),}
\end{equation*}
où $a$ est une matrice symétrique d'ordre $l$ à coefficients dans
$\O$.
\begin{co}\label{co2.4.1}
Le groupe $K_{B}$ est engendré par $P_{B}$ et $\overline{N}_{B}$.
\end{co}
\begin{dem}
C'est une conséquence immédiate du lemme \ref{lem2.6.1} et de ce que
d'une part $N_{B}$ est un sous-groupe de $P_{B}$ et que d'autre part
on a  $\varsigma_{B}=x_{B}(I_{l})y_{B}(-I_{l})x_{B}(I_{l})$. 
\end{dem}
\begin{rem}
Il résulte de \cite{bruhat-tits-1987} que les sous-groupes $K_{B}$
sont les stabilisateurs des sommets de l'immeuble du groupe
$Sp(W)$. Pour être plus précis, si l'on choisit une base autoduale
$(e_{1},\ldots,e_{r},f_{1},\ldots,f_{r})$ de $W$ et si, pour $0\leq
l\leq r$, on pose $B_{l}=\O\varpi e_{1}\oplus\cdots\O\varpi
e_{l}\oplus\O e_{l+1}\oplus\cdots\oplus\O e_{r}\oplus\O
f_{1}\oplus\cdots\oplus\O f_{r}$, alors les $K_{B_{l}}$, $0\leq l\leq
r$, sont les stabilisateurs des sommets d'une chambre, les cas où
$l=0$ et $l=r$ correspondant aux deux sommets hyperspéciaux. De plus,
parmi les sous-groupes $K_{B_{l}}$, $0\leq l\leq r$, seuls $K_{B_{0}}$
et $K_{B_{r}}$ sont conjugués par un élément du groupe des similitudes
symplectiques~: si on désigne par $d_{r}$ la similitude symplectique
dont la matrice dans la base autoduale choisie est
\begin{equation*}
\begin{pmatrix}
\varpi I_{r} & 0\\
0 & I_{r}
\end{pmatrix}
\end{equation*}
on a $K_{B_{r}}= d_{r}K_{B_{0}}d_{r}^{-1}$.
\end{rem}

\subsection{}\label{2.7}
On garde les notations du paragraphe précédent. Il est clair que
$B^{*}\backslash\varpi B^{*}$ est une partie de $W$ qui est invariante
sous l'action de $K_{B}$.
\begin{lem}\label{lem2.7.1}
(i) Pour tout $b\in B^{*}\backslash B$ (resp. $b\in B\backslash\varpi
  B^{*}$), on a $K'_{B}b=b+B$ (resp. $K'_{B}b=b+\varpi B^{*}$).

(ii) Les orbites de $K_{B}$ dans $B^{*}\backslash\varpi B^{*}$ sont
  $B^{*}\backslash B$ et $B\backslash \varpi B^{*}$.
 
(iii) Les orbites de $P_{B}$ dans $B^{*}\backslash\varpi B^{*}$ sont
  $B^{*}\backslash A$, $A\backslash B$ et $B\backslash \varpi B^{*}$.
\end{lem}
\begin{dem}
Il est clair que $B^{*}\backslash B$ et $B\backslash \varpi B^{*}$
(resp. $B^{*}\backslash A$, $A\backslash B$ et $B\backslash \varpi
B^{*}$) sont invariants sous l'action de $K_{B}$ (resp. $P_{B}$). Il
reste à montrer (i) et le fait que chacun de ces ensembles est une
$K_{B}$-orbite (resp. $P_{B}$-orbite).

On vérifie facilement que l'image de $P_{B}$ par le morphisme
$p_{Sp({\mathrm{b}
}^{*})}\times p_{Sp(\mathrm{b})}$ est $P_{{\mathrm{b}
}^{*}}\times
Sp(\mathrm{b})$. Comme les orbites de $P_{{\mathrm{b}
}^{*}}$ dans ${\mathrm{b}
}^{*}$
sont $\{0\}$, $\mathrm{x}\backslash\{0\}$ et ${\mathrm{b}
}^{*}\backslash\mathrm{x}$
et comme $Sp(\mathrm{b})$ (resp. $Sp({\mathrm{b}
}^{*})$) agit transitivement
dans $\mathrm{b}\backslash\{0\}$ (resp. ${\mathrm{b}
}^{*}\backslash\{0\}$) il
suffit de montrer que, pour un élément particulier $b$ de
$B^{*}\backslash B$ (resp. $B\backslash\varpi B^{*}$), on a
$K'_{B}b=b+B$ (resp. $K'_{B}b=b+\varpi B^{*}$)~: on peut prendre
$b=e_{1}$ (resp $b=e_{l+1}$).

Dans la suite de la démonstration, on identifie chaque élément de $W$
(resp. $GL(W)$) avec le vecteur colonne de ses coordonnées (resp. sa
matrice) dans la base $(e_{1},\ldots,e_{l},$ $e_{l+1},\ldots,e_{r},
f_{1},\ldots,f_{l},f_{l+1},\ldots,f_{r})$.  

Supposons que $b=e_{1}\in B^{*}\backslash B$ et soit $u\in B$. Alors,
on a $u=\left(
\begin{smallmatrix}
\varpi\lambda_{1}\\
\lambda_{2}\\
\mu_{1}\\
\mu_{2}
\end{smallmatrix}\right)
$ avec $\lambda_{1},\mu_{1}\in\O^{l}$,
$\lambda_{2},\mu_{2}\in\O^{r-l}$, et l'élément $k=\left(
\begin{smallmatrix}
a^{-1} & 0\\
0 &\transposee{a}
\end{smallmatrix}\right)$
où 
\begin{equation*}
a  =
\begin{pmatrix}
I_{l}+\varpi\lambda_{1}\transposee{e}_{1}&0\\
\lambda_{2}\transposee{e}_{1}&I_{r-l}
\end{pmatrix}
\end{equation*}
est dans $K'_{B}$ et vérifie $k(e_{1}+u)=e_{1}+u'$ avec $u'=\left(
\begin{smallmatrix}
0\\
0\\
\mu'_{1}\\
\mu'_{2}
\end{smallmatrix}
\right)$, $\mu'_{i}$, $i=1,2$, étant à coefficients dans $\O$. Mais
alors l'élément $h=\left(
\begin{smallmatrix}
I_{r}&0\\
c&I_{r}
\end{smallmatrix}
\right)$
avec 
\begin{equation*}
c  =
-\begin{pmatrix}
\mu'_{1}\transposee{e}_{1}+e_{1}\transposee{\mu}'_{1}-
\transposee{\mu}'_{1}e_{1}E_{11}
&e_{1}\transposee{\mu}'_{2}\\
\mu'_{2}\transposee{e}_{1}&0
\end{pmatrix}
\end{equation*}
est dans $K'_{B}$ et vérifie $hk(e_{1}+u)=e_{1}$.

Plaçons nous dans le cas où $b=e_{l+1}\in B\backslash\varpi B^{*}$ et
soit $u\in\varpi B^{*}$. Alors, on peut écrire $b=\left(
\begin{smallmatrix}
0\\e_{1}\\0\\0
\end{smallmatrix}
\right)$ et $u=\left(
\begin{smallmatrix}
\varpi\lambda_{1}\\
\varpi\lambda_{2}\\
\mu_{1}\\
\varpi\mu_{2}
\end{smallmatrix}
\right) $ avec $\lambda_{1},\mu_{1}\in\O^{l}$,
$\lambda_{2},\mu_{2}\in\O^{r-l}$, et l'élément $k=\left(
\begin{smallmatrix}
a^{-1}&0\\
0&\transposee{a}
\end{smallmatrix}
\right)$ où
\begin{equation*}
a=
\begin{pmatrix}
I_{l}&\varpi\lambda_{1}\transposee{e}_{1}\\
0&I_{r-l}+\varpi\lambda_{2}\transposee{e}_{1}
\end{pmatrix}
\end{equation*}
est dans $K'_{B}$ et vérifie $k(e_{l+1}+u)=e_{l+1}+u'$ avec $u'=\left(
\begin{smallmatrix}
0\\
0\\
\mu'_{1}\\
\varpi\mu'_{2}
\end{smallmatrix}
\right)$, $\mu'_{i}$, $i=1,2$, étant à coefficients dans $\O$. Mais
alors l'élément $h=\left(
\begin{smallmatrix}
I_{r}&0\\
c&I_{r}
\end{smallmatrix}
\right)$
avec 
  
\begin{equation*}
 c=-
\begin{pmatrix}
0&\mu'_{1}\transposee{e}_{1}\\
e_{1}\transposee{\mu'}_{1}&\varpi(\mu'_{2}\transposee{e}_{1}
+e_{1}\transposee{\mu'}_{2})-\varpi\transposee{\mu'}_{2}e_{1}E_{11}
\end{pmatrix}
\end{equation*}
est dans $K'_{B}$ et vérifie $hk(e_{l+1}+u)=e_{l+1}$.
\end{dem}

\section{La représentation de Weil}\label{3}
Dans cette section, nous rappelons certains résultats sur la
représentation de Weil et le groupe métaplectique. On peut consulter
\cite{weil-1964}, \cite[Chapitre 2]{moeglin-vigneras-waldspurger-1987},
\cite{rao-1993}, \cite{cliff-mcneilly-szechtman-2000}, \cite{pan-2001}
et \cite{cliff-mcneilly-szechtman-2003}.

Dans la suite, les représentations induites que nous considérons sont
des représentations induites unitaires à droite. Plus précisément, si
$G$ est un groupe localement compact, $H$ un de ses sous-groupes fermés, tous
deux unimodulaires, et $\rho$ une représentation unitaire de $H$ dans
l'espace de Hilbert $\mathcal{H}$, $\Ind_{H}^{G}\rho$ désigne la
représentation de $G$ agissant par translations à droite dans l'espace
des fonctions $\varphi : G\longrightarrow\mathcal{H}$ qui vérifient
\begin{gather*}
\varphi(hg)=\rho(h)\varphi(g)\mbox{, }g\in G\mbox{, }h\in H,\\
\sideset{}{_{H\backslash G}}\int\Vert\varphi(g)\Vert^{2}\diff\dot{g}
<+\infty
\end{gather*}
où $\diff\dot{g}$ désigne une mesure de Radon $G$-invariante sur
$H\backslash G$.

\subsection{}\label{3.1}
Dans ce paragraphe, $F$ désigne de nouveau soit un anneau local fini
de caractéristique différente de $2$, soit un corps local non
archimédien de caractéristique résiduelle différente de $2$ et $\psi$
un caractère unitaire non trivial de $F$. Si $F$ est un anneau local fini, on
suppose que $\psi$ est primitif, en ce sens que son noyau ne contient
pas d'autre idéal que $\{0\}$.

Soit $(W,\beta)$ un $F$-espace symplectique~; si $F$ est
un corps, on note $2r$ la dimension de $W$ sur $F$. Le
groupe de Heisenberg associé est l'ensemble $H(W)= W\times F$, muni de
la multiplication
\begin{eqnarray*}
(w,t)(w',t') &=& (w+w', t+t' + \frac{1}{2}\beta(w,w')), w, w' \in W
  \mbox {, } t,t' \in F.
\end{eqnarray*}
Muni de la topologie produit, $H(W)$ est un groupe localement compact.

Le centre de $H(W)$ est $\{0\}\times F$ que l'on identifie à $F$, via
l'application~:
\begin{equation*}
t \in F \longmapsto (0, t) \in H(W).
\end{equation*}
On identifie également $W$ à un sous-ensemble de $H(W)$ au moyen de
l'application~:
\begin{equation*}
w\in W\longmapsto (w,0)\in H(W).
\end{equation*}
Dans ces conditions, on a $w^{-1}=-w$, pour tout $w\in W$.

Par le théorème de Stone Von-Neumann (voir
\cite[2]{cliff-mcneilly-szechtman-2000} pour le cas où $F$ n'est pas
un corps), nous savons qu'il existe une unique (à équivalence près)
représentation unitaire irréductible $(\rho_{\psi}, \mathcal{H})$ de
$H(W)$ de caractère central $\psi$. On l'appelle la représentation de
Schrödinger de caractère central $\psi$ de $H(W)$. On la note
$\rho_{W,\psi}$ lorsque l'on souhaite faire apparaître qu'il s'agit
d'une représentation du groupe de Heisenberg construit sur $W$.

Nous rappelons comment construire de telles représentations, selon
\cite[Chapitre 2]{moeglin-vigneras-waldspurger-1987} et
\cite{cliff-mcneilly-szechtman-2000}. Soit $A\subset W$ un sous-groupe
fermé. On pose
\begin{equation*}
A^{*}=\{v \in W , \quad \psi(\beta(v,a))=1,
\mbox{ pour tout } a \in A \}.
\end{equation*}
Alors $A^{*}$ est un sous-groupe fermé de $W$, appelé sous-groupe dual
: lorsque $F$ est un corps local et $A$ est un réseau, on retrouve la
notion de réseau dual (voir \ref{eq2.3.2}). On dit que $A$ est
autodual si $A=A^{*}$. Lorsque $F$ est un corps local, les
sous-espaces lagrangiens et les réseaux autoduaux sont des exemples de
sous-groupes autoduaux.

Lorsque $F$ est un anneau local fini, tout sous-espace lagrangien
est un sous-groupe autodual. On sait que si $F$ est un corps, il
existe des lagrangiens. Nous allons montrer qu'il en existe pour tout
anneau local fini.

Soit $\mathfrak{m}$ l'idéal maximal de $F$. Son annulateur
$\mathfrak{s}$ est l'unique idéal propre minimal de $F$. Fixons un
élément non nul $s_{0}\in\mathfrak{s}$. Alors l'application
$a\longmapsto as_{0}$ passe au quotient en un isomorphisme de
$F$-modules du corps résiduel $\overline{F}=F/\mathfrak{m}$ sur
$\mathfrak{s}$. On désigne par $\mu$ un tel isomorphisme.

\begin{lem}\label{lem3.1.1}
Soit $F$ un anneau local fini et $(W,\beta)$ un $F$-espace
symplectique. Soit $A\subset W$ un sous-module.

(i) On a $A^{\perp\perp}=A$.

(ii) On suppose que l'on a les inclusions $\mathfrak{m}A\subset
A^{\perp}\subset A$. Alors, $A/A^{\perp}$ est naturellement muni d'une
structure de $\overline{F}$-espace vectoriel, la restriction
$\beta_{\vert A}$ de $\beta$ à $A$ prend ses valeurs dans
$\mathfrak{s}$ et l'application $\mu^{-1}\circ\beta_{\vert A}$ passe au
quotient à $A/A^{\perp}$ en une forme symplectique.
\end{lem}
\begin{dem}
Pour le point (i), voir \cite[Lemma
  2.1]{cliff-mcneilly-szechtman-2003}.

Soit $A$ un sous-module de $W$ vérifiant les hypothèses du
(ii). L'unique chose à vérifier est que $\beta_{\vert A}$ est à
valeurs dans $\mathfrak{s}$. Or ceci résulte de ce que, compte tenu
des hypothèses,
$\mathfrak{m}\beta(A,A)=
\beta(\mathfrak{m}A,A)\subset\beta(A^{\perp},A)=\{0\}$.
\end{dem}
\begin{lem}\label{lem3.1.2}
Soit $F$ un anneau local fini et $W$ un $F$-espace symplectique. Alors
$W$ admet des sous-espaces lagrangiens.
\end{lem}
\begin{dem}
Soit $A\subset W$ un sous-$F$-module totalement isotrope et
maximal. Il suffit de montrer que $A=A^{\perp}$. Or, il suit de
\cite[Lemma 10.1]{cliff-mcneilly-szechtman-2003} que
$\mathfrak{m}A^{\perp}\subset A$ (le résultat cité est énoncé sous
l'hypothèse que $A$ est isotrope, $Sp(W)$-invariant et maximal, mais
il est clair que la démonstration reste valable si l'on y remplace
$Sp(W)$ par n'importe lequel de ses sous-groupes). Comme
$A^{\perp\perp}=A$, on voit alors que $A^{\perp}$ vérifie les
hypothèses du (ii) du lemme \ref{lem3.1.1}. Par suite,
$\mu^{-1}\circ\beta_{\vert A^{\perp}}$ passe au quotient en une forme
symplectique sur le $\overline{F}$-espace vectoriel $A^{\perp}/A$ et
l'image réciproque dans $A^{\perp}$ d'un lagrangien de cet espace
symplectique est un sous-groupe autodual de $W$ contenant $A$. On a
donc $A=A^{\perp}$. \end{dem} 

Soit $A$ un sous-groupe autodual de $W$. Alors, $A\times F$ est un
sous-groupe fermé de $H(W)$ et la formule
\begin{equation*}
\psi_{A}(w,t)=\psi(t)\mbox{, }(w,t)\in A\times F,
\end{equation*}
définit un caractère de $A\times F$. On note alors
$\rho_{\psi}^{A}=\Ind_{A\times F}^{H(W)}\psi_{A}$ la représentation
induite du caractère $\psi_{A}$ de $A\times F$ à $H(W)$. La
représentation $\rho_{\psi}^{A}$ se réalise dans l'espace
$\mathcal{H}_{\psi}^{A}$ des fonctions $\var:W\longrightarrow\Bbb{C}$ qui
vérifient
\begin{gather*}
\var(a+w)=\psi(\frac{1}{2}\beta(w,a))\var(w)\mbox{, }a\in A\mbox{,
}w\in W,\\
\int_{W/A}\vert\var(w)\vert^{2}\diff\dot w<+\infty,
\end{gather*}
où $\diff\dot w$ désigne une mesure de Haar sur le groupe additif $W/A$, par
la formule
\begin{equation}\label{eq3.1.1}
\rho_\psi^A(w,t)\var(w')= \psi( t + \frac{1}{2}\beta(w', w))\var(w'
+w)\mbox{, } w, w' \in W\mbox{, }t \in F.
\end{equation}

\begin{theo}\label{theo3.1.1}
  (i) Soit $A$ un sous-groupe autodual de $W$. Alors la représentation
  $\rho_{\psi}^{A}$ est irréductible et appartient à la classe de
  $\rho_{\psi}$.

(ii) Soit $A$ et $B$ deux sous-groupes autoduaux de $W$ et $\diff\dot
  b$ une mesure de Haar sur le groupe additif quotient $B/A\cap
  B$. Pour $\var\in\mathcal{H}_{\psi}^{A}$ continue à support compact
  modulo $A$, on pose
\begin{equation}\label{eq3.1.2}
\mathcal{I}_{B,A}\var(w)=\int_{B/A\cap
  B}\var(w+b)\psi(\frac{1}{2}\beta(b,w))\diff\dot b\mbox{, }w\in W.
\end{equation}
Alors l'intégrale \ref{eq3.1.2} converge pour tout $w$, la fonction
$\mathcal{I}_{B,A}\var$ est un élément de $\mathcal{H}_{\psi}^{B}$ et
$\mathcal{I}_{B,A}:\mathcal{H}_{\psi}^{A}\longrightarrow\mathcal{H}_{\psi}^{B}$ se
prolonge de manière unique en un opérateur continu, encore noté
$\mathcal{I}_{B,A}$, entrelaçant les représentations $\rho_{\psi}^{A}$
et $\rho_{\psi}^{B}$.
\end{theo}  

Supposons que $F$ soit un corps. Soit $X$ et $Y$ deux lagrangiens de
$W$ tels que $W=X\oplus Y$. Alors, pour un bon choix de la mesure de
Haar $\diff x$ sur $X$, l'application $\var\longmapsto\var_{|X}$ est
une isométrie de $\mathcal{H}_{\psi}^{Y}$ sur $L^{2}(X,\diff x)$
permettant d'identifier ces deux espaces.  La représentation
$\rho_{\psi}^{Y}$ réalisée dans $L^{2}(X,\diff x)$ est alors donnée
par les formules suivantes~:
\begin{eqnarray}
\rho_\psi^Y(x,0)\var(x')&=& \var(x'
+x)\mbox{, } x, x' \in X\label{eq3.1.3}\\
\rho_\psi^Y(y,0)\var(x')&=&\psi(\beta(x', y))\var(x')\mbox{,
}y\in Y\mbox{, }x'\in X.\label{eq3.1.4}
\end{eqnarray}

Lorsque $F$ est un corps local et $A$ est un lagrangien (resp. un
réseau) de $W$, on dit que $\rho_{\psi}^{A}$ est un modèle de
Schrödinger (resp. latticiel) de la représentation $\rho_{\psi}$.

\subsection{}\label{3.2}
On garde les notations du paragraphe précédent. Le groupe symplectique
$Sp(W)$ opère par automorphismes dans $H(W)$ par la formule
\begin{equation*}
x.(w,t) = (xw,t) , \, \mbox{ pour tout } w\in W\mbox{, } t \in F.
\end{equation*}

Cette action fixe $F$ point par point.

Soit $(\rho_\psi, \mathcal{H})$ une représentation unitaire
irréductible de $H(W)$ de caractère central $\psi$ et désignons par
$U(\mathcal{H})$ le groupe des transformations unitaires de l'espace
de Hilbert $\mathcal{H}$.  On note $\widehat{Sp(W)}_{\psi}$ le
sous-groupe topologique de $Sp(W) \times U(\mathcal{H})$ formé des
couples $(g, M(g))$ vérifiant~:
\begin{eqnarray}\label{eq3.2.1}
M (g)\rho_\psi(h)M(g)^{-1} = \rho_\psi(g.h), h \in H(W).
\end{eqnarray}
On voit que si $(g, M(g)) \in \widehat{Sp(W)}_{\psi}$ alors $(g,
zM(g)) \in \widehat{Sp(W)}_{\psi}$ pour tout $z \in \Bbb{U}$, le
groupe des nombres complexes de module $1$. Il
suit du théorème de Stone Von-Neumann que $\widehat{Sp(W)}_{\psi}$
est une extension centrale de $Sp(W)$ par les opérateurs scalaires de
norme $1$~; on a une suite exacte courte~:
\begin{eqnarray}\label{eq3.2.2}
1 \longrightarrow \Bbb{U} \stackrel{\alpha_1}{\longrightarrow}
\widehat{Sp(W)}_{\psi}\stackrel{\alpha_2} \longrightarrow Sp(W)
\longrightarrow 1,
\end{eqnarray}
où $\alpha_1: z \longmapsto (1, z\mathrm{Id})$ et $\alpha_2 : (g,
M(g))\longmapsto g$. La représentation métaplectique $S_\psi$ de
$\widehat{Sp(W)}_{\psi}$ agissant dans $\mathcal{H}$ est donnée par
\begin{eqnarray}\label{eq3.2.3}
S_\psi(g, M(g)) = M(g).
\end{eqnarray}
On la note $S_{W,\psi}$ si l'on souhaite faire apparaître qu'il s'agit
de la représentation métaplectique associée au groupe symplectique $Sp(W)$.

En fait, on peut se restreindre à un revêtement d'ordre au plus $2$ de
$Sp(W)$ grâce au résultat suivant~:
 
\begin{theo}\label{theo3.2.1}
(i) Si $F$ est un anneau local fini, il existe un morphisme de groupes
  de $Sp(W)$ dans $\widehat{Sp(W)}_{\psi}$ qui scinde la suite
  exacte \ref{eq3.2.2}. Lorsque $F$ est le corps fini $\Bbb{F}_{q}$,
  à l'exception du cas où $q=3$ et $\dim W=2$, cet homomorphisme est
  unique.

(ii) Si $F$ est un corps local non archimédien et $W$ est non nul, la
  suite exacte \ref{eq3.2.2} ne se scinde pas, mais il existe un
  unique sous-groupe $Mp(W)_{\psi}$ de $\widehat{Sp(W)}_{\psi}$ tel
  que la restriction de $\alpha_{2}$ à ce sous-groupe soit surjective
  et ait un noyau d'ordre $2$. Ce sous-groupe est fermé et la
  restriction $\alpha_{2}$ à ce sous-groupe en fait un revêtement
  d'ordre $2$ de $Sp(W)$.
\end{theo}

Voir \cite{steinberg-ip-1962} et \cite{cliff-mcneilly-szechtman-2000}
pour le point (i) et \cite{weil-1964} pour le point (ii).

Lorsque $F$ est un corps local non archimédien, on sait qu'à
isomorphisme canonique près, il n'existe qu'un revêtement d'ordre $2$
de $Sp(W)$, non trivial. Ce revêtement s'identifie donc de manière
canonique au sous-groupe $Mp(W)_{\psi}$ de $\widehat{Sp(W)}_{\psi}$~:
désormais on le note $Mp(W)$ et on l'appelle le groupe métaplectique
de l'espace symplectique $W$ ou le revêtement métaplectique du groupe
$Sp(W)$. Dans la suite, on appelle représentation de Weil de type
$\psi$ et on note $S_{\psi}$ ou $S_{W,\psi}$ si l'on tient à préciser,
la restriction de la représentation $S_{\psi}$ à $Mp(W)$.

Si $A$ est un sous-groupe autodual de $W$, on note
$\widehat{Sp(W)}_{\psi}^{A}$, $Mp(W)_{\psi}^{A}$ et
$S_{\psi}^{A}$ ou $S_{W,\psi}^{A}$ les objets précédents construits en
utilisant la réalisation $\rho_{\psi}^{A}$ dans l'espace
$\mathcal{H}_{\psi}^{A}$ de la représentation $\rho_{\psi}$. On parle
alors de modèle de Schrödinger ou de modèle latticiel de la
représentation de Weil, suivant que $A$ est un lagrangien ou bien
un réseau.

Lorsque $F$ est un anneau local fini, on appelle représentation de
Weil de type $\psi$ et on note $S_{\psi}$ ou $S_{W,\psi}$, la
représentation composée de $S_{\psi}$ avec un homomorphisme de
$Sp(W)$dans $ \widehat {Sp(W)}_{\psi}$ scindant la suite exacte
\ref{eq3.2.2}. \'Etant donnée une représentation de Weil, toutes les
autres s'écrivent com\-me produit de celle-ci par un caractère
unitaire de $Sp(W)$. Si $A$ est un sous-groupe autodual de $W$, on
désignera par $S_{\psi}^{A}$ ou $S_{W,\psi}^{A}$ une représentation de
Weil construite en utilisant la réalisation $\rho_{\psi}^{A}$ dans
l'espace $\mathcal{H}_{\psi}^{A}$ de la représentation $\rho_{\psi}$.

Plus généralement, $F$ étant toujours un anneau local fini, si $G$ est
un groupe fini et si $\gamma:G\longrightarrow Sp(W)$ est un morphisme
de groupes, on appelle représentation de Weil de type $\psi$ de $G$
relative au morphisme $\gamma$, toute représentation $S$ de $G$ dans
l'espace de la représentation $\rho_{\psi}$ vérifiant la relation
\begin{equation*}
  S(g)\rho_{\psi}(h)S(g)^{-1}=\rho_{\psi}(\gamma(g).h)\mbox{, }g\in
  G\mbox{, }h\in H(W).
\end{equation*}
Ici encore, étant donnée une représentation de Weil, toutes les autres
s'écrivent com\-me produit de celle-ci par un caractère unitaire de
$G$.

Si $F=\Bbb{F}_{q}$, $Sp(W)$ admet une unique représentation de Weil,
sauf dans le cas $q=3$ et $\dim W=2$ (à l'exception de ce cas, le
groupe $Sp(W)$ est parfait). Cela dit, il existe dans tous les cas un
choix canonique, que nous décrivons dans le paragraphe suivant~: dans
la suite, on convient que la représentation de Weil $S_{\psi}$ est
celle correspondant à ce choix.

\subsection{}\label{3.3}
Dans ce paragraphe, on garde les notations du paragraphe précédent et
on suppose que $F$ est un corps fini. On rappelle que si $\eta$  est
un caractère additif de $F$, l'indice de Weil relatif à $\eta$ est
l'application $\omega_{\eta}:F^{\times}\longrightarrow\Bbb{U}$ définie par
\begin{equation}\label{eq3.3.-2}
  \omega_{\eta}(a)=[F]^{-\frac{1}{2}}\sum_{t\in F}\eta(at^{2}).
\end{equation}
Si $a\in F^{\times}$, on désigne par $a\eta$ le caractère additif
$t\mapsto\eta(at)$ de $F$. On pose alors
\begin{equation}\label{eq3.3.-1}
\omega=\omega_{\frac{1}{2}\psi}.
\end{equation}

Soit $X$ un sous-espace
lagrangien de $W$ et soit $(e_{1},\ldots,e_{r},f_{1},\ldots,f_{r})$
une base symplectique de $W$ telle que $(e_{1},\ldots,e_{r})$ soit une
base de $X$. 

Si $g\in Sp(W)$, on pose
\begin{equation}\label{eq3.3.0}
j(g)=r-\dim X\cap gX
\end{equation}
Pour $0\leq j\leq r$, soit
\begin{equation*}
\Omega_j = \{g \in Sp(W)\vert j(g)=j\}.
\end{equation*}
Alors, $\Omega_{0}$ n'est autre que $P$, le sous-groupe parabolique de $Sp(W)$
stabilisateur de $X$, chaque $\Omega_{j}$ est une double classe sous
$P$ et $Sp(W)$ est la réunion disjointe des $\Omega_{j}$, $0\leq j\leq
r$. Pour $S\subset\{1,\ldots,r\}$, on désigne par $\tau_{S}$ l'élément
de $Sp(W)$ tel que
\begin{equation*}
\tau_{S}.e_{i}=
\begin{cases}
 f_i  & \mbox{ si }  i \in S\\
e_i & \mbox{ si } i \notin S 
\end{cases}
\qquad
\tau_{S}.f_{i}=
\begin{cases}
-e_i  & \mbox{ si }  i \in S\\
f_i & \mbox{ si } i \notin S 
\end{cases}
\end{equation*}
Alors, si $\mbox{ card} S = j$, on a~: $\Omega_j = P \tau_SP$.

D'après \cite[Lemma 5.1]{rao-1993}, il existe une application $\theta
: Sp(W)\longrightarrow \Bbb{F}^{\times} / (\Bbb{F}^{\times})^{2}$
telle que

\begin{itemize}
	\item $\theta(p_1 gp_2)= \theta(p_1 )\theta(
          g)\theta(p_2)\mbox{, } p_1, p_2 \in P\mbox{, }g \in Sp(W)$ ;
	\item $\theta(\tau_S) = 1\mbox{, } S\subset \{1, \ldots, r\}$ ; 
	\item $\theta(p)= \det_{X}(p) \mod
          (\Bbb{F}^\times)^2\mbox{, }p \in P$.
\end{itemize}
Si $g \in Sp(W)$, on pose
\begin{equation}\label{eq3.3.1}
m(g)= \omega (\theta(g))^{-1} \omega(1)^{1- j(g)},
\end{equation}
où $\omega$ est l'indice de Weil relatif au caractère $\frac{1}{2}\psi$
(voir les formules \ref{eq3.3.-2} et \ref{eq3.3.-1}).

Soit $\mu_{X}$ la mesure de Haar sur $X$ de masse totale
$1$. Si $g \in Sp(W)$ on pose
\begin{equation}\label{eq3.3.2}
  S_{\psi}^{X}(g)=m(g) q^{\frac{j(g)}{2}} M_{X}(g),
\end{equation}
où $M_{X}(g)$ est l'opérateur de l'espace $\mathcal{H}_{\psi}^{X}$ défini par
\begin{equation}\label{eq3.3.3}
M_{X}(g)\var(w)= \int_{X} \psi(\frac{1}{2}\beta(x,
w))\var(g^{-1}(w+x))\diff \mu_{X}(x)\mbox{,
}\var\in\mathcal{H}_{\psi}^{X}\mbox{, }w\in W.
\end{equation}
Le résultat suivant est alors conséquence de \cite[4.2,
  Lemma]{pan-2001}~:
\begin{lem}\label{lem3.3.1}
L'application 
\begin{equation*}
  g\longmapsto(g,S_{\psi}^{X}(g))
\end{equation*}
est un homomorphisme de groupes de $Sp(W)$ dans $\widehat{Sp(W)}_{\psi}^{X}$
scindant la suite exacte \ref{eq3.2.2}.
\end{lem}

\subsection{}\label{3.4} 
Dans ce paragraphe, on suppose que $F$ est un
anneau local fini et on se donne un sous-groupe $G$ de $Sp(W)$. Nous
rappelons un certain nombre des résultats de
\cite{cliff-mcneilly-szechtman-2000} et
\cite{cliff-mcneilly-szechtman-2003} concernant la décomposition des
représentations de Weil sur $F$.

On désigne par $\mathcal{F}(W)$ le $\mathbb{C}$-espace vectoriel des
fonctions à valeurs complexes sur $W$. On munit $\mathcal{F}(W)$ de la
structure de $Sp(W)$-module définie par
$g.\varphi(w)=\varphi(g^{-1}w)$, $w\in W$.  

Soit $(\rho_{\psi},\mathcal{H})$ une représentation de Schrödinger de
$H(W)$ et soit $S$ une représentation de Weil de type $\psi$ de
$Sp(W)$. La représentation $S$ induit une structure de $Sp(W)$-module
dans $End_{\mathbb{C}}(\mathcal{H})$ définie par
$g.A=S(g)AS(g^{-1})$. On a alors le résultat suivant (voir \cite[Theorem
  4.5]{cliff-mcneilly-szechtman-2000})~:

\begin{theo}\label{theo3.4.1}
L'application $\varphi\mapsto\sum_{w\in W}\varphi(w)\rho_{\psi}(w,0)$ de
$\mathcal{F}(W)$ dans $End_{\mathbb{C}}(\mathcal{H})$ est un
isomorphisme de $Sp(W)$-module.
\end{theo}

On munit le sous-groupe $G$ de $Sp(W)$ de la mesure de Haar
normalisée. Le théorème précédent a pour conséquence le résultat
suivant (voir \cite[Corollary 4.6]{cliff-mcneilly-szechtman-2000})~:
\begin{co}\label{co3.4.1}
Le carré de la norme, comme élément de $L^{2}(G)$, du caractère de la
restriction à $G$ d'une représentation de Weil de $Sp(W)$ est égal au
nombre d'orbites du groupe $G$ dans $W$.
\end{co}

Les résultats présentés dans la suite sont ceux de \cite[paragraphe
  4]{cliff-mcneilly-szechtman-2003}. Soit $U\subset W$ un sous-module
isotrope, non nul et $Sp(W)$-invariant. Alors,
$\overline{U^{\perp}}=U^{\perp}/U$ est un $F$-espace symplectique pour
la forme symplectique $\overline{\beta}$ définie par
$\overline{\beta}(v+U,w+U)=\beta(v,w)$. On pose
$H(U^{\perp})=U^{\perp}\times F$~: c'est un sous-groupe du groupe de
Heisenberg $H(W)$ et l'application $p_{U}:(v,t)\mapsto (v+U,t)$ est un
morphisme surjectif de $H(U^{\perp})$ sur
$H(\overline{U^{\perp}})$. De plus l'action du groupe $Sp(W)$ dans
$U^{\perp}$ passe au quotient en une action dans
$\overline{U^{\perp}}$ par automorphismes symplectiques, fournissant
ainsi un morphisme de groupes $r_{U}:Sp(W)\longrightarrow
Sp(\overline{U^{\perp}})$. Soit $\rho_{U^{\perp},\psi}$ l'inflation de
la représentation de Schrödinger de $H(\overline{U^{\perp}})$ via
$p_{U}$ à $H(U^{\perp})$ et soit $\sigma$ une représentation de Weil
de type $\psi$ de $Sp(W)$ relative à $r_{U}$, que l'on réalise dans le
même espace $\mathcal{H}_{\sigma}$~: $\sigma$ est le produit tensoriel
de l'inflation d'une représentation de Weil de
$Sp(\overline{U^{\perp}})$ via $r_{U}$ à $Sp(W)$ par un caractère de
$Sp(W)$. Comme nous l'avons vu, le groupe $Sp(W)$ agit par
automorphismes dans le groupe $H(W)$ et il laisse stable le
sous-groupe $H(U^{\perp})$. Nous pouvons donc considérer le groupe
produit semi-direct $Sp(W)\ltimes H(W)$ et son sous-groupe
$Sp(W)\ltimes H(U^{\perp})$. On définit alors la représentation
$\sigma\rho_{U^{\perp},\psi}$ de $Sp(W)\ltimes H(U^{\perp})$ dans
l'espace $\mathcal{H}_{\sigma}$ en posant
$\sigma\rho_{U^{\perp},\psi}(gh)= \sigma(g)\rho_{U^{\perp},\psi}(h)$,
$g\in Sp(W)$ et $h\in H(U^{\perp})$. Soit $R^{U,\sigma}$ la
représentation de $Sp(W)\ltimes H(W)$ définie par
\begin{equation*}
R^{U,\sigma}=\Ind_{Sp(W)\ltimes H(U^{\perp})}^{Sp(W)\ltimes H(W)}
\sigma\rho_{U^{\perp},\psi}
\end{equation*}
On note $S^{U,\sigma}$ la restriction de la représentation
$R^{U,\sigma}$ à $Sp(W)$. Si $\chi$ un caractère de $Sp(W)$, on a
$S^{U,\chi\otimes\sigma}= \chi\otimes S^{U,\sigma}$.

On vérifie que l'espace de la représentation $R^{U,\sigma}$
s'identifie, via l'application de restriction au sous-ensemble $W$ du
sous-groupe $H(W)$ du produit semi-direct $Sp(W)\ltimes H(W)$, à
l'espace $\mathcal{H}^{U}$ des fonctions $\varphi$ définies sur $W$, à
valeurs dans $\mathcal{H}_{\sigma}$,
qui vérifient la relation
\begin{equation}\label{eq3.4.1}
\varphi(x+u)=
\psi(\frac{1}{2}\beta(x,u))\rho_{U^{\perp},\psi}(u)
\varphi(x) \mbox{, }x\in W\mbox{, }u\in U^{\perp}.
\end{equation}

Si $x\in W$, on note $\hat{x}$ (resp. $\dot{x}$) sa classe dans le
module quotient $W/U^{\perp}$ (resp. $W/U$).  Le groupe $Sp(W)$ agit
naturellement dans $W/U^{\perp}$ et si $x\in W$, on désigne par $G(\hat{x})$ le
stabilisateur de $\hat{x}=x+U^{\perp}$ dans $G$. Rappelons que l'on identifie
$x$ à l'élément $(x,0)$ du groupe de Heisenberg $H(W)$. Alors
$xG(\hat{x})x^{-1}$ est un sous-groupe de $G\ltimes H(U^{\perp})$. Plus
précisément, si $g\in G(\hat{x})$, on a
\begin{equation*}
xgx^{-1}=g(g^{-1}x-x,\frac{1}{2}\beta(x,g^{-1}x)).
\end{equation*}
On définit une représentation, que nous noterons $\sigma_{G,\dot{x}}$, de
$G(\hat{x})$ dans $\mathcal{H}_{\sigma}$ en posant pour $g\in G(\hat{x})$,
\begin{equation*}
\begin{split}
\sigma_{G,\dot{x}}(g)&= \sigma\rho_{U^{\perp},\psi}(xgx^{-1})\\&=
\sigma(g)
\rho_{U^{\perp},\psi}(g^{-1}x-x,\frac{1}{2}\beta(x,g^{-1}x))
\end{split}
\end{equation*}
(on vérifie que la représentation $\sigma_{G,\dot{x}}$ ne dépend que
de la classe $\dot{x}$ de $x$ dans $W/U$) et on considère la
représentation $S_{G,\dot{x}}$ de $G$ définie par
\begin{equation*}
S_{G,\dot{x}}
=\Ind_{G(\hat{x})}^{G}\sigma_{G,\dot{x}}.
\end{equation*}
On désigne par $\mathcal{H}_{G,\dot{x}}$ l'espace de la représentation
$S_{G,\dot{x}}
$. Il est constitué des fonctions
$\varphi:G\mapsto\mathcal{H}_{\sigma}$ qui vérifient
\begin{equation*}
\varphi(hg)=\sigma_{G,\dot{x}}(h)\varphi(g)\mbox{, }g\in G\mbox{,
}h\in G(\hat{x}).
\end{equation*}
Si $\varphi\in\mathcal{H}_{G,\dot{x}}$,
$\mathcal{I}_{G,\dot{x}}\varphi$ désigne la fonction élément de
$\mathcal{H}^{U}$ à support dans l'orbite de $x$ modulo $U^{\perp}$
sous l'action de $G$ telle que
\begin{equation*}
\mathcal{I}_{G,\dot{x}}\varphi(g.x)=\sigma(g)\varphi(g^{-1})\mbox{, }g\in G.
\end{equation*}

Si $x=0$, on a $S_{G,0}=\sigma_{\vert G}$. Lorsque $G=Sp(W)$, on note
simplement $\sigma_{\dot{x}}$ (resp. $S_{\dot{x}}$) au lieu de
$\sigma_{G,\dot{x}}$ (resp. $S_{G,\dot{x}}$).

Si $\mathcal{H}$ est l'espace d'une représentation $\pi$ de $Sp(W)$ ,
on désigne par $\mathcal{H}^{\pm}$ le sous-espace propre de $-Id$ dans
$\mathcal{H}$ pour la valeur propre $\pm 1$. Comme $\{\pm Id\}$ est un
sous-groupe central de $Sp(W)$, le sous-espace $\mathcal{H}^{\pm}$ est
invariant sous la représentation $\pi$ et on désigne par $\pi^{\pm}$
la représentation de $Sp(W)$ dans ce sous-espace qui en résulte,
lorsque celui-ci est non nul. On a clairement
$\mathcal{H}=\mathcal{H}^{+}\oplus\mathcal{H}^{-}$ et donc
$\pi=\pi^{\pm}$ ou $\pi=\pi^{+}\oplus\pi^{-}$, suivant que l'un des
espaces $\mathcal{H}^{\pm}$ est nul ou non.

Dans ce qui suit, on suppose que $G=Sp(W)$. Soit $x\in W$. Si $x\in
U^{\perp}$, on a $G(\hat{x})=Sp(W)$ et deux cas se présentent alors~: si
$\overline{U^{\perp}}=\{0\}$, $S_{\dot{x}}=\sigma_{\dot{x}}$ est de dimension $1$
et irréductible~; dans le cas contraire, $S_{\dot{x}}=\sigma_{\dot{x}}$ est somme
directe des sous-représentations $S^{\pm}_{x}$. Par contre, si
$x\notin U^{\perp}$, $-Id$ n'est pas un élément de $G(\hat{x})$ et la
représentation $\sigma_{\dot{x}}$ s'étend au sous-groupe
$\widetilde{G(\hat{x})}=\{\pm Id\}G(\hat{x})$ de $Sp(W)$ en la représentation
$\sigma^{\pm}_{x}$ en décidant que $\sigma^{\pm}_{x}(-Id)=\pm Id$.
Alors, on a
\begin{equation*}
S^{\pm}_{\dot{x}}=\Ind_{\widetilde{G(\hat{x})}}^{Sp(W)}\sigma^{\pm}_{\dot{x}}
\end{equation*}
et donc 
\begin{equation*}
S_{\dot{x}}=S^{+}_{\dot{x}}\oplus S^{-}_{\dot{x}}.
\end{equation*}

\begin{lem}\label{lem3.4.1}
(i) La restriction à $H(W)$  de la représentation $R^{U,\sigma}$
est équivalente à la représentation de Schrödinger $\rho_{\psi}$.

(ii) La représentation $S^{U,\sigma}$ est une représentation de Weil
de type $\psi$ de $Sp(W)$ et toute représentation de Weil de type
$\psi$ de $Sp(W)$ est obtenue ainsi.

(iii) Soit $G$ un sous-groupe de $Sp(W)$. Si $x\in W$, le sous-espace
$\mathcal{H}^{U}_{G,\hat{x}}$ de $\mathcal{H}^{U}$ constitué des
fonctions à support dans l'orbite de $x$ modulo $U^{\perp}$ sous
l'action de $G$ est $(S^{U,\sigma})_{\vert G}$-invariant. De plus,
l'application $\mathcal{I}_{G,\dot{x}}$ est un isomorphisme de
$G$-modules de l'espace $\mathcal{H}_{G,\dot{x}}$ de la représentation
$S_{G,\dot{x}}$ sur $\mathcal{H}^{U}_{G,\hat{x}}$. Ainsi, la
  représentation de $G$ induite par $S^{U,\sigma}$ dans
  $\mathcal{H}^{U}_{G,\hat{x}}$ est équivalente à la représentation
  $S_{G,\dot{x}} $.

(iv) Si $X\subset W$ est un ensemble de représentants des $G$-orbites
dans $W/U^{\perp}\backslash\{0\}$, on a~:

\begin{equation}\label{eq3.4.2}
(S^{U,\sigma})_{\vert G}= \sigma_{\vert G}\oplus(\bigoplus_{x\in
    X}S_{G,\dot{x}} )
\end{equation}

(v) Si $X\subset W$ est un ensemble de représentants des
$Sp(W)$-orbites dans $W/U^{\perp}\backslash\{0\}$, on a~:
\begin{equation}\label{eq3.4.3}
S^{U,\sigma}
=
\begin{cases}
\sigma\oplus(\bigoplus_{x\in X}S^{+}_{\dot{x}}\oplus
S^{-}_{\dot{x}}) &\mbox{ si }
\overline{U^{\perp}}=\{0\}\\ 
(\sigma^{+}\oplus\sigma^{-})\oplus(\bigoplus_{x\in
  X}S^{+}_{\dot{x}}\oplus S^{-}_{\dot{x}}) &\mbox{ si }
\overline{U^{\perp}}\neq\{0\}
\end{cases}
\end{equation}
\end{lem}

Jusqu'à la fin de ce paragraphe, on suppose que l'anneau local fini $F$
n'est pas un corps. Soit $U\subset W$ un sous-module isotrope,
$Sp(W)$-invariant et maximal pour cette propriété. Rappelons que
$\mathfrak{m}$ désigne l'idéal maximal de $F$, $\mathfrak{s}$ son
idéal minimal et $\overline{F}$ le corps résiduel de $F$. D'après
\cite[Lemma 10.1]{cliff-mcneilly-szechtman-2003}, on a
$\mathfrak{m}U^{\perp}\subset U$. Alors, il suit du lemme
\ref{lem3.1.1} que $\overline{U^{\perp}}$ est naturellement muni d'une
structure de $\overline{F}$-espace vectoriel, que la forme
symplectique $\overline{\beta}$ prend ses valeurs dans $\mathfrak{s}$
et que, si $\mu$ est un isomorphisme de $F$-module de
$\overline{F}=F/\mathfrak{m}$ sur $\mathfrak{s}$, l'application $\mu
^{-1}\circ\overline{\beta}$ en fait un $\overline{F}$-espace vectoriel
symplectique. De plus, les groupes symplectiques pour
$\overline{\beta}$ et $\mu^{-1}\circ\overline{\beta}$ sont identiques
et notés $Sp(\overline{U^{\perp}})$.

Désignons par $H_{F}(\overline{U^{\perp}})$
(resp. $H_{\overline{F}}(\overline{U^{\perp}})$) le groupe de
Heisenberg construit sur $\overline{U^{\perp}}$ considéré comme un
$F$-module (resp. $\overline{F}$-espace vectoriel) symplectique. Alors
l'isomorphisme $\mu$ induit un plongement 
\begin{equation*}
\begin{split}
\mu_{*}:H_{\overline{F}}(\overline{U^{\perp}})&
\longrightarrow  H_{F}(\overline{U^{\perp}})\\
(u,t)&\longmapsto (u,\mu(t)).
\end{split}
\end{equation*}
L'application $\mu_{*}$ commute avec l'action de
$Sp(\overline{U^{\perp}})$ par automorphismes dans les deux groupes de
Heisenberg. 

La composée de la représentation de Schrödinger de
$H_{F}(\overline{U^{\perp}})$ de caractère central $\psi$ avec le
morphisme $\mu_{*}$ est clairement la représentation de Schrödinger de
$H_{\overline{F}}(\overline{U^{\perp}})$ de caractère central
$\psi\circ\mu$~: comme $\psi$ est un caractère primitif,
$\psi\circ\mu$ est un caractère non trivial de $\overline{F}$. Il est
alors immédiat qu'une représentation de Weil de type $\psi$ de
$Sp(\overline{U^{\perp}})$ est également une représentation de Weil de
type $\psi\circ\mu$ et réciproquement.  

Il suit de ces considérations que les représentations de Weil pour
$Sp(W)$ peuvent être construites à partir des représentations de Weil pour un
groupe symplectique sur le corps résiduel de $F$, de même que les
composantes d'une telle représentation qui apparaissent dans le lemme
\ref{lem3.4.1}.

\subsection{}\label{3.5} Soit $n$ un entier naturel et $B\subset W$ un bon
réseau. Nous allons appliquer les résultats du paragraphe précédent
pour déterminer la décomposition en irréductibles des représentations
de Weil des groupes $Sp({\mathrm{b}}_{n})$ et $Sp({\mathrm{b}}^{*}_{n})$
introduits au paragraphe \ref{2.5}.

Rappelons que ${\mathrm{b}}_{n}$ et ${\mathrm{b}}^{*}_{n}$ sont des modules
sur $O_{n}=\O/\varpi^{n+1}\O$. L'annulateur dans $O_{n}$ de
${\mathrm{b}}_{n}$ (resp. ${\mathrm{b}}^{*}_{n}$) est trivial si et seulement
si $\varpi B^{*}\subsetneq B$ (resp. $B\subsetneq B^{*}$), tandis que
${\mathrm{b}}_{n}={\mathrm{b}}^{*}_{n-1}$
  (resp. ${\mathrm{b}}^{*}_{n}={\mathrm{b}}_{n-1}$) dans le cas
  contraire. Nous pouvons donc nous placer dans la situation
  suivante~: $\psi$ est un caractère primitif de $O_{n}$ et on
  considère les représentations de Weil de type $\psi$ pour
  $Sp({\mathrm{b}}_{n})$ (resp. $Sp({\mathrm{b}}^{*}_{n})$) lorsque $\varpi
  B^{*}\subsetneq B$ (resp. $B\subsetneq B^{*}$).
\begin{lem}\label{lem3.5.1}
(i) L'espace symplectique ${\mathrm{b}}_{n}$ admet un unique sous-module
  isotrope invariant sous l'action de $Sp({\mathrm{b}}_{n})$ et maximal,
  noté $\mathrm{u}$. On a 

${\mathrm{u}}=\varpi^{m+1}B^{*}/\varpi^{n+1}B^{*}$,
${\mathrm{u}}^{\perp}=\varpi^{m}B/\varpi^{n+1}B^{*}$ et
$\overline{\mathrm{u}^{\perp}}=\mathrm{b}$, si $n=2m$,

${\mathrm{u}}=\varpi^{m+1}B/\varpi^{n+1}B^{*}$,
  ${\mathrm{u}}^{\perp}=\varpi^{m+1}B^{*}/\varpi^{n+1}B^{*}$ et
  $\overline{\mathrm{u}^{\perp}}= \mathrm{b}^{*} $, si $n=2m+1$.

(ii) L'espace symplectique ${\mathrm{b}}^{*}_{n}$ admet un unique
  sous-module isotrope invariant sous l'action de
  $Sp({\mathrm{b}}^{*}_{n})$ et maximal, noté $\mathrm{u}^{*}$. On a

${\mathrm{u}^{*}}=\varpi^{m}B/\varpi^{n}B$, 
${\mathrm{u}}^{*\perp}=\varpi^{m}B^{*}/\varpi^{n}B$ et 
$\overline{\mathrm{u}^{*\perp}}=\mathrm{b}^{*}$, si $n=2m$,

${\mathrm{u}}^{*}=\varpi^{m+1}B^{*}/\varpi^{n}B$, 
${\mathrm{u}}^{*\perp}=\varpi^{m}B/\varpi^{n}B$ et 
$\overline{\mathrm{u}^{*\perp}}=\mathrm{b}$, 
si $n=2m+1$.
\end{lem}
\begin{dem}
Comme le morphisme naturel de $K_{B}$ dans les groupes sym\-plectiques
$Sp({\mathrm{b}}_{n})$ et $Sp({\mathrm{b}}_{n}^{*})$ est surjectif et comme
les seuls réseaux de $W$ qui sont $K_{B}$-invariants sont ceux
proportionnels à $B$ ou $B^{*}$, on voit que les sous-modules de
${\mathrm{b}}_{n}$ ou ${\mathrm{b}}^{*}_{n}$ invariants par le groupe
symplectique sont image de $\varpi^{s}B$ ou $\varpi^{s}B^{*}$,
$s\geq0$ par la projection naturelle de $B$ (resp. $B^{*}$) sur
${\mathrm{b}}_{n}$ ou ${\mathrm{b}}^{*}_{n}$. On en déduit facilement le lemme.
\end{dem}

\begin{theo}\label{theo3.5.1}
Soit $n\in\Bbb{N}$ et $\psi$ un caractère primitif de $O_{n}$. On
reprend les notations du lemme \ref{lem3.5.1}.

(i) Soit $S$ une représentation de Weil de type $\psi$ de
$Sp({\mathrm{b}}_{n})$. Alors, il existe une unique représentation de
Weil, $\sigma$, de type $\psi$ de $Sp({\mathrm{b}}_{n})$
relative au morphisme naturel de $Sp({\mathrm{b}}_{n})$ dans
$Sp({\mathrm{b}})$ (resp. $Sp({\mathrm{b}}^{*})$) si $n$ est pair
(resp. impair) telle que $S=S^{\mathrm{u},\sigma}$.

Dans tous les cas, les représentations apparaissant dans la formule
\ref{eq3.4.3} du lemme \ref{lem3.4.1} sont irréductibles et deux à
deux inéquivalentes. Lorsque $B\subsetneq
B^{*}$, il y en a $2n+2$~; lorsque  $B=B^{*}$, il y en a $n+2$. Elles
sont monomiales si $B=B^{*}$ et $n$ est impair.

(ii)  Soit $S$ une représentation de Weil de type $\psi$ de
$Sp({\mathrm{b}}^{*}_{n})$. Alors, il existe une unique représentation de
Weil, $\sigma$, de type $\psi$ de $Sp({\mathrm{b}}^{*}_{n})$
relative au morphisme naturel de $Sp({\mathrm{b}}^{*}_{n})$ dans
$Sp({\mathrm{b}}^{*})$ (resp. $Sp({\mathrm{b}})$) si $n$ est pair
(resp. impair) telle que $S=S^{\mathrm{u}^{*},\sigma}$.

Dans tous les cas, les représentations apparaissant dans la formule
\ref{eq3.4.3} du lemme \ref{lem3.4.1} sont irréductibles et deux à
deux inéquivalentes. Lorsque $\varpi B^{*}\subsetneq
B$, il y en a $2n+2$~; lorsque  $\varpi B^{*}=B$, il y en a $n+2$. Elles
sont monomiales si $\varpi B^{*}=B$ et $n$ est impair.
\end{theo}
\begin{dem}
Nous faisons la démonstration dans le cas (i), la démonstration dans
l'autre cas étant similaire. La première assertion est claire, compte
tenu des résultats du paragraphe précédent et de ceux du lemme
\ref{lem3.5.1}. Soit $c$ (resp. $d$) le nombre de
$Sp({\mathrm{b}}_{n})$-orbites dans ${\mathrm{b}}_{n}$
(resp. ${\mathrm{b}}_{n}/{\mathrm{u}}^{\perp}$). Compte tenu du
corollaire \ref{co3.4.1} et du lemme \ref{lem3.4.1}, il suffit de
montrer que $c=2d$, si $B\subsetneq B^{*}$ ou $n$ est pair, et
$c=2d-1$ dans le cas contraire. Cependant, comme le morphisme naturel
de $K_{B}$ dans $Sp({\mathrm{b}}_{n})$ est surjectif, on voit que
$c-1$ est le nombre des $K_{B}$-orbites dans
$B\backslash\varpi^{n+1}B^{*}$, tandis que $d-1$ est celui des
$K_{B}$-orbites dans $B\backslash\varpi^{m}B$
(resp. $B\backslash\varpi^{m+1}B^{*}$), si $n=2m$
(resp. $n=2m+1$). Or, d'après le lemme \ref{lem2.7.1}, l'ensemble des
$K_{B}$-orbites dans $B\backslash\varpi^{m}B$ est
$\{\varpi^{t}B\backslash\varpi^{t+1}B^{*}\vert 0\leq t\leq
m-1\}\sqcup\{\varpi^{t}B^{*}\backslash\varpi^{t}B\vert 1\leq t\leq
m\}$, de cardinal $2m$
(resp. $\{\varpi^{t}B\backslash\varpi^{t+1}B\vert 0\leq t\leq m-1\}$,
de cardinal $m$) si $B\subsetneq B^{*}$ (resp. $B=B^{*}$). Pour la
même raison, l'ensemble des $K_{B}$-orbites dans
$B\backslash\varpi^{m+1}B^{*}$ est
$\{\varpi^{t}B\backslash\varpi^{t+1}B^{*}\vert 0\leq t\leq
m\}\sqcup\{\varpi^{t}B^{*}\backslash\varpi^{t}B\vert 1\leq t\leq m\}$,
de cardinal $2m+1$ (resp. $\{\varpi^{t}B\backslash\varpi^{t+1}B\vert
0\leq t\leq m\}$, de cardinal $m+1$) si $B\subsetneq B^{*}$
(resp. $B=B^{*}$). Le théorème en découle.
\end{dem}
\begin{rem}
Dans \cite{cliff-mcneilly-szechtman-2003} G. Cliff, D. McNeilly et
F. Szechtman établissent que la formule \ref{eq3.4.3} du lemme
\ref{lem3.4.1} donne la décomposition en irréductibles de la
représentation de Weil de $Sp(W)$ considérée, lorsque $R$ est un anneau
local fini et principal et $W$ est un $R$-module libre. Notre théorème
en est une conséquence, uniquement lorsque $B=B^{*}$ ou $B=\varpi
B^{*}$. 

Par ailleurs, dans \cite{dutta-prasad-2010} Dutta et Prasad démontrent
que les représentations de Weil du groupe symplectique d'un module
symplectique de type fini sur un anneau local fini et principal se
décompose en irréductibles avec multiplicité $1$ et donnent une
description de ces derniers en terme de la combinatoire des orbites du
groupe symplectique dans le module symplectique. Cependant, leur
description ne fait pas apparaître explicitement ces irréductibles
comme modules induits.
\end{rem}

\subsection{}\label{3.6}
Dans ce paragraphe, on suppose que $F$ est un corps local non
archimédien. Soit $X\subset W$ un sous-espace lagrangien, $P_{X}$ le
sous-groupe parabolique de $Sp(W)$ stabilisateur de $X$, $N_{X}$ son
radical unipotent et $Y\subset W$ un lagrangien supplémentaire de $X$.

Un élément $g$ de $Sp(W)$ s'identifie à une matrice 
$\left( 
\begin{smallmatrix}
a & b\\
c & d
\end{smallmatrix}\right)
$ dans laquelle $a\in End(X)$, $b\in Hom (Y,X)$, $c\in Hom (X,Y)$ et $d\in
End(Y)$. Remarquons que la forme symplectique $\beta$ induit une
dualité entre $X$ et $Y$, permettant d'identifier de manière naturelle
$Y$ avec l'espace dual de $X$~: on a donc une notion d'opérateur
symétrique $b:Y\longrightarrow X$. Le sous-groupe $N_{X}$ est alors le
sous-groupe des matrices de la forme $ x(b)=\left(\begin{smallmatrix}
  1 & b\\ 0 & 1
\end{smallmatrix}\right)
$
avec $b\in Hom(Y,X)$ symétrique.

Il résulte de \cite[Chapitre II, 11, Lemme]{duflo-1982} que la
restriction du revêtement métaplectique à $N_{X}$ admet un scindage
unique. On peut donc considérer le groupe $N_{X}$ comme un sous-groupe du
groupe métaplectique $Mp(W)$.

Rappelons que nous avons réalisé la représentation $\rho_{\psi}^{X}$
de $H(W)$ dans l'espace $L^{2}(Y,\diff y)$, où $\diff y$ est une
mesure de Haar sur $Y$~: on en déduit une réalisation $S_{\psi}^{X}$
de la représentation métaplectique dans ce même espace. 

Le résultat suivant est alors conséquence de \cite[Chapitre 2,II.6 et II.9,
  Lemme]{moeglin-vigneras-waldspurger-1987}~:
\begin{lem}\label{lem3.6.1}
La restriction de la représentation métaplectique $S_{\psi}^{X}$ à $N_{X}$
est donnée par
\begin{equation*}
S_{\psi}^{X}(x(b))\var(y)=\psi(\frac{1}{2}\beta(by,y))\var(y)\mbox{,
}y\in Y\mbox{, }\var\in L^{2}(Y,\diff y),
\end{equation*}
où $b$ parcourt l'ensemble des homomorphismes symétriques de $Y$ dans $X$.
\end{lem}

\section{\'Etude du revêtement métaplectique au dessus d'un sous-groupe
  compact maximal}\label{4} 
Dans cette section, nous montrons que le revêtement métaplectique est
scindé au dessus de tout sous-groupe compact maximal et nous donnons
une description de la restriction de la représentation métaplectique à
un tel sous-groupe. Pour ce faire, nous présentons le modèle latticiel
généralisé de la représentation métaplectique introduit par Shu Yen
Pan dans \cite[numéro 2.3]{pan-2001}. Désormais, nous reprenons les
notations de la section \ref{2}.

\subsection{}\label{4.1}
Soit $B\subset W$ un bon réseau. On a vu que ${\mathrm{b}
}^{*}=B^{*}/B$ est
un $\Bbb{F}_{q}$-espace vectoriel symplectique de dimension $2l$, où
$l=l(B)$. On peut donc considérer la représentation
$(\rho_{\overline{\psi}},\mathcal{H}_{\overline{\psi}})$ du groupe de
Heisenberg $H({\mathrm{b}
}^{*})$ de caractère central $\overline{\psi}$. Il
suit de ce que $B$ est un bon réseau que
$H(B^{*})=B^{*}\times\varpi^{\lambda_{\psi}-1}\O$ est un sous-groupe
de $H(W)$. Alors, l'application
$p_{H({\mathrm{b}
}^{*})}:H(B^{*})\longrightarrow H({\mathrm{b}
}^{*})$ définie par
\begin{equation*}
  p_{H({\mathrm{b}
}^{*})}(w,t)=(p_{{\mathrm{b}
}^{*}}(w),p_{\Bbb{F}_{q}}(\varpi^{1-\lambda_{\psi}}t))
\end{equation*}
est un homomorphisme surjectif de groupes localement compacts. On
désigne par $\widetilde{\rho}_{\overline{\psi}}$ la représentation de
$H(B^{*})$ relevant la représentation $\rho_{\overline{\psi}}$ via le
morphisme $p_{H({\mathrm{b} }^{*})}$. On note de même le prolongement de
$\widetilde{\rho}_{\overline{\psi}}$ à
$\overline{H}(B^{*})=B^{*}\times k$ dont la restriction à $k$ est un
multiple du caractère $\psi$.

On considère alors la représentation $\rho_{\psi}^{B}=
\Ind_{\overline{H}(B^{*})}^{H(W)}\widetilde{\rho}_{\overline{\psi}}$
que l'on réalise dans l'espace de Hilbert $\mathcal{H}_{\psi}^{B}$ des
fonctions $\var:W\longrightarrow\mathcal{H}_{\overline{\psi}}$ qui
vérifient
\begin{gather}
\var(b + w)= \psi(\frac{1}{2}\beta( w,
b))\widetilde{\rho}_{\overline{\psi}}(b)( \var(w))\mbox{, }b\in
B^{*}\mbox{, }w\in W,\label{eq4.1.1}\\ 
\int_{W/B^*}\Vert\var(w)\Vert^2
\diff\dot w < + \infty,\label{eq4.1.2}
\end{gather}
dans lequel $H(W)$ agit comme suit~:
\begin{equation*}
(\rho_\psi^B (w,t) \var)(w')= \psi( t + \frac{1}{2}\beta(w', w))
  \var(w' +w)\mbox{, } w, w' \in W\mbox{, }  t \in k .
\end{equation*}

\begin{pr}\label{pr4.1.1}
  La représentation $\rho_{\psi}^{B}$ est unitaire irréductible,
  équivalente à $\rho_{\psi}$.
\end{pr}
\begin{dem}
Soit $A$ un réseau autodual tel que $B\subset A\subset B^{*}$ et soit
$\mathrm{x}=A/B$ dont on a vu que c'est un lagrangien de
${\mathrm{b}
}^{*}$. Comme $p_{H({\mathrm{b}
}^{*})}^{-1}(\mathrm{x}\times\Bbb{F}_{q})=
A\times\varpi^{\lambda_{\psi}-1}\O$, il est clair que la
représentation $\Ind_{A\times k}^{\overline{H}(B^{*})}\psi_{A}$ est
équivalente à la représentation
$\widetilde{\rho}_{\overline{\psi}}$. Notre résultat est alors
conséquence du théorème d'induction par étage.
\end{dem}
On dit que $\rho_{\psi}^{B}$ est le modèle latticiel généralisé de la
représentation $\rho_{\psi}$ associé au réseau $B$. On note
$\widehat{Sp(W)}_{\psi}^{B}$, $Mp(W)_{\psi}^{B}$ et
$S_{\psi}^{B}$ les objets $\widehat{Sp(W)}_{\psi}$,
$Mp(W)_{\psi}$ et $S_{\psi}$ obtenus en utilisant la
réalisation $\rho_{\psi}^{B}$ de la représentation $\rho_{\psi}$~; on
dit que $S_{\psi}^{B}$ est le modèle latticiel généralisé de la
représentation métaplectique associé au réseau $B$.

Soit $A$ un réseau autodual tel que $B\subset A\subset B^{*}$,
$\mathrm{x}=A/B$ le lagrangien correspondant de ${\mathrm{b}
}^{*}$ et
$\rho_{\overline{\psi}}^{\mathrm{x}}$ la représentation d'espace
$\mathcal{H}_{\overline{\psi}}^{\mathrm{x}}$ équivalente à la
représentation $\rho_{\overline{\psi}}$ de $H({\mathrm{b}
}^{*})$. On réalise
alors l'espace $\mathcal{H}_{\psi}^{B}$ comme l'espace des fonctions à
valeurs dans $\mathcal{H}_{\overline{\psi}}^{\mathrm{x}}$ vérifiant les
conditions \ref{eq4.1.1} et \ref{eq4.1.2}. Il suit facilement de la
démonstration de la proposition \ref{pr4.1.1} que l'on a alors le
résultat suivant (voir aussi \cite[4.3, Lemma]{pan-2001})
\begin{co}\label{co4.1.1}
Soit $B$ un bon réseau de $W$ et $A$ un réseau autodual tel que
$B\subset A\subset B^{*}$. Alors, pour tout $\var\in\mathcal{H}_{\psi}^{B}$ la
fonction $\mathcal{I}_{A,B}(\var)$ définie sur $W$ par
\begin{equation*}
\mathcal{I}_{A,B}(\var)(w)=\var(w)(0)\mbox{, }w\in W,
\end{equation*}
est un élément de $\mathcal{H}_{\psi}^{A}$ et l'application
$\mathcal{I}_{A,B}:\mathcal{H}_{\psi}^{B}\longrightarrow\mathcal{H}_{\psi}^{A}$
 est un opérateur entrelaçant les représentations
$\rho_{\psi}^{B}$ et $\rho_{\psi}^{A}$. 
\end{co}

\subsection{}\label{4.2}
On garde les notations du paragraphe précédent et on reprend celles du
paragraphe \ref{2.4}. On note $\widetilde{S}_{\overline{\psi}}$ la
représentation de $K_{B}$ relevant la représentation
$S_{\overline{\psi}}$ de $Sp({\mathrm{b}
}^{*})$ via le morphisme
$p_{Sp({\mathrm{b}
}^{*})}$.

Considérons le produit semi-direct $K_{B}\ltimes H(W)$. Il contient
comme sous-groupe $K_{B}\ltimes\overline{H}(B^{*})$. On définit une
représentation
$\widetilde{S}_{\overline{\psi}}\widetilde{\rho}_{\overline{\psi}}$ de
$K_{B}\ltimes\overline{H}(B^{*})$ en décidant que
\begin{equation}\label{eq4.2.1}
\widetilde{S}_{\overline{\psi}}\widetilde{\rho}_{\overline{\psi}}(gh)=
\widetilde{S}_{\overline{\psi}}(g)\widetilde{\rho}_{\overline{\psi}}(h)
\mbox{, }g\in K_{B}\mbox{, }h\in\overline{H}(B^{*}).
\end{equation}
On considère la représentation $R^{B}_{\psi}$ de $K_{B}\ltimes H(W)$ définie par
\begin{equation}\label{eq4.2.0}
R^{B}_{\psi}=\Ind_{K_{B}\ltimes\overline{H}(B^{*})}^{K_{B}\ltimes H(W)}
\widetilde{S}_{\overline{\psi}}\widetilde{\rho}_{\overline{\psi}}
\end{equation}
dont l'espace est clairement $\mathcal{H}_{\psi}^{B}$ et on note
$S_{\psi}^{B}$ sa restriction à $K_{B}$.
\begin{lem}\label{lem4.2.1}
(i) Si $g\in K_{B}$,  l'opérateur $S_{\psi}^{B}(g)$ de
  $\mathcal{H}_{\psi}^{B}$ est donné par 
\begin{equation}\label{eq4.2.2}
  S_{\psi}^{B}(g)\var(w)=
  \widetilde{S}_{\overline{\psi}}(g)[\var(g^{-1}w)]\mbox{,
  }\var\in\mathcal{H}_{\psi}^{B}\mbox{, }w\in W.
\end{equation}

(ii) L'application
\begin{equation}\label{eq4.2.3}
g\longmapsto(g,S_{\psi}^{B}(g))
\end{equation}
est un homomorphisme de groupes de $K_{B}$ dans $\widehat{Sp(W)}_{\psi}^{B}$
scindant la suite exacte \ref{eq3.2.2}.
\end{lem}
\begin{dem}
Le point (i) est conséquence immédiate de la définition de la
représentation $R^{B}_{\psi}$ comme représentation induite.

Le point (ii) résulte de ce que la restriction de $R^{B}_{\psi}$ à
$H(W)$ est la représentation de Schrödinger $\rho_{\psi}^{B}$.
\end{dem}
En fait, on a le résultat plus précis suivant~:
\begin{theo}\label{theo4.2.1}
Le scindage au dessus de $K_{B}$ défini par \ref{eq4.2.3} est à
valeurs dans le groupe métaplectique de $W$.
\end{theo}

\subsection{}\label{4.3}
On garde les notations du paragraphe précédent. Soit $A$ un réseau
autodual tel que $B\subset A\subset B^{*}$. Pour démontrer le théorème
\ref{theo4.2.1}, nous devons relier la réalisation $S_{\psi}^{B}$ de
la représentation métaplectique à la réalisation
$S_{\psi}^{A}$. Commençons par donner une description de cette
dernière.

Soit $\diff w$ une mesure de Haar sur $W$ et soit $\mathcal{P}$ le
projecteur orthogonal de $L^{2}(W,\diff w)$ sur
$\mathcal{H}_{\psi}^{A}$~; en fait, on a
\begin{equation*}
\mathcal{P}\var(w)=
\int_{A}\psi(\frac{1}{2}\beta(a,w))\var(w+a)\diff a\mbox{, }w\in W,
\end{equation*}
où $\diff a$ est la mesure de Haar normalisée sur le groupe compact
$A$.  On note $L$ la représentation unitaire naturelle de $Sp(W)$
dans $L^{2}(W,\diff w)$ définie par $L(g)\var(w)=\var(g^{-1}w)$. Si
$g\in Sp(W)$, on pose $M_{A}(g)=P\circ L(g)_{|\mathcal{H}_{\psi}^{A}}$. On
vérifie immédiatement que, pour tout $\var\in\mathcal{H}_{\psi}^{A}$, on a
\begin{equation}\label{eq4.3.1}
  (M_{A}(g)\var)(w)= [A/A\cap gA]^{-1}\sum_{a \in A /A \cap gA}
  \psi(\frac{1}{2}\beta( a, w)) \var(g^{-1}(a + w)).
\end{equation}

On pose $K=K_{A}$. La proposition suivante rassemble les résultats de
\cite[Chapitre 2, II.8 et II.10, Lemme]{moeglin-vigneras-waldspurger-1987}~:
\begin{pr}\label{pr4.3.1}
(i) Pour tout $g\in Sp(W)$, il existe un nombre strictement positif
  $\nu_{A}(g)$ tel que $(g,\nu_{A}(g)M_{A}(g))$ soit un élément de
  $\widehat{Sp(W)}_{\psi}^{A}$.

(ii) Si $g\in K$, on a $\nu_{A}(g)=1$ et
\begin{equation*}
M_{A}(g)\var(w)=\var(g^{-1}w)\mbox{, }\var\in\mathcal{H}_{\psi}^{A}\mbox{,
}w\in W.
\end{equation*}
(iii) L'application $g\longmapsto(g,M_{A}(g))$ est un homomorphisme de
groupes de $K$ dans le groupe métaplectique de $W$. En particulier, il
existe un relèvement du groupe $K$ dans le groupe métaplectique pour
lequel la restriction de la représentation métaplectique
$S_{\psi}^{A}$ à $K$ soit donnée par
\begin{equation*}
S_{\psi}^{A}(g)=M_{A}(g)\mbox{, }g\in K.
\end{equation*}
\end{pr}
\begin{rems}
(i) Lorsque le réseau $B$ est autodual, l'affirmation du théorème
\ref{theo4.2.1} n'est autre que l'assertion (iii) de la proposition
\ref{pr4.3.1}.

(ii) Lorsque le corps résiduel est de cardinal au moins 4, le groupe
$K$ admet un unique relèvement dans le groupe métaplectique (voir
\cite[Chapitre 2, II 10, Remarque]{moeglin-vigneras-waldspurger-1987}).
\end{rems}

\subsection{}\label{4.4}
On garde les notations du paragraphe précédent. Dans celui-ci, nous
allons donner une description plus concrète des opérateurs
$M_{A}(g)$. Tout d'abord, si $v\in W$, on note $\delta_{v}$ la
fonction sur $W$ telle que
\begin{equation}\label{eq4.4.1}
\delta_{v}(w)=
\begin{cases}
\psi(\frac{1}{2}\beta(v, w-v))& \mbox{ si } w \in v+ A \\
0 & \mbox{ si } w \notin v+ A.
\end{cases}
\end{equation}
Remarquons que, si $v\in W$ et $a\in A$, on a 
\begin{equation}\label{eq4.4.2}
\delta_{v+a}=\psi(\frac{1}{2}\beta(a, v))\delta_{v}.
\end{equation}
De plus, si $v\in W$ et $g\in K$, on a 
\begin{equation}\label{eq4.4.3}
M_{A}(g)\delta_{v}=\delta_{g.v}.
\end{equation}
D'autre part, il est clair que si $\Sigma\subset W$ est un système de
représentants des éléments de $W/A$, la famille
$(\delta_{v})_{v\in\Sigma}$ est une base hilbertienne de l'espace
$\mathcal{H}_{\psi}^{A}$.

On a le résultat suivant~:
\begin{lem}\label{lem4.4.1}
Soit $g\in Sp(W)$ et $v\in W$. Alors, on a
\begin{equation}\label{eq4.4.4}
M_{A}(g)\delta_{v}=[A/A\cap gA]^{-1}
\sum_{a\in gA/A\cap gA}\psi(\frac{1}{2}\beta(v,g^{-1}a))\delta_{gv+a}.
\end{equation}
\end{lem}
\begin{dem}
Utilisant la formule \ref{eq4.3.1}, on voit que
$M_{A}(g)\delta_{v}(w)$ est non nul seulement si $w\in gv+A+gA$. Par
suite, il existe des nombres complexes $c(a)$, pour $a$ parcourant un
système de représentants dans $gA$ des éléments de
$(A+gA)/A=gA/A\cap gA$, tels que 
\begin{equation}\label{eq4.4.5}
M_{A}(g)\delta_{v}=\sum_{a\in gA/A\cap gA} c(a)\delta_{gv+a}.
\end{equation}
Or, il suit des relations \ref{eq4.4.5} et \ref{eq4.3.1} que pour
$a\in gA/A\cap gA$, on a
\begin{eqnarray*}
c(a)&=&M_{A}(g)\delta_{v}(gv+a)\\
&=& [A/A\cap gA]^{-1} \sum_{b \in A /A \cap gA}
  \psi(\frac{1}{2}\beta( b, gv+a)) \delta_{v}(v+g^{-1}(a + b)).
\end{eqnarray*}
Or $\delta_{v}(v+g^{-1}(a+b))$ est non nul si et seulement si
$g^{-1}(a+b)\in A$, soit encore $a+b\in gA$, soit finalement $b\in
A\cap gA$. On en déduit que 
\begin{equation*}
c(a)=[A/A\cap gA]^{-1}\psi(\frac{1}{2}\beta(v,g^{-1}a))
\end{equation*}
comme voulu.
\end{dem}
On déduit facilement du lemme \ref{lem4.4.1} le résultat suivant~:
\begin{co}\label{co4.4.1}
Pour tout $g\in Sp(W)$, on a 
\begin{equation*}
\nu_{A}(g)=[gA/A\cap gA]^{-\frac{1}{2}}[A/A\cap gA].
\end{equation*}
\end{co}

\subsection{}\label{4.5}
On garde les notations des paragraphes précédents. On rappelle
l'opérateur d'entrelacement $\mathcal{I}_{A,B}$ défini au paragraphe
\ref{4.1}, corollaire \ref{co4.1.1}, ainsi que les fonctions $j$ et
$m$ sur le groupe $Sp({\mathrm{b} }^{*})$ relative au lagrangien
$\mathrm{x}$ et au caractère $\overline{\psi}$ définies au paragraphe
\ref{3.3}, équations \ref{eq3.3.0} et \ref{eq3.3.1}. En particulier,
l'espace $\mathcal{H}_{\psi}^{B}$ est l'espace des fonctions sur $W$ à
valeurs dans l'espace $\mathcal{H}_{\overline{\psi}}^{\mathrm{x}}$ de
la représentation $\rho_{\overline{\psi}}^{\mathrm{x}}$ vérifiant les
conditions \ref{eq4.1.1} et \ref{eq4.1.2}.

En fait, on peut réaliser l'espace
$\mathcal{H}_{\overline{\psi}}^{\mathrm{x}}$ comme l'espace des fonctions
$\var$ sur $B^{*}$ telles que 
\begin{equation*}
\var(b+a)=\psi(\frac{1}{2}\beta(b,a))\var(b)\mbox{, }b\in B^{*}\mbox{,
}a\in A.
\end{equation*}
Dans ces conditions, on déduit facilement du lemme \ref{lem3.3.1} et
des formules \ref{eq3.3.2} et \ref{eq3.3.3} que, pour $g\in K_{B}$ et
$\var\in\mathcal{H}_{\overline{\psi}}^{\mathrm{x}}$, on a
\begin{equation}\label{eq4.5.1}
\begin{split}
  \widetilde{S}_{\overline{\psi}}(g)\var(b)&=
  m(p_{Sp({\mathrm{b}
}^{*})}(g))q^{\frac{j(p_{Sp({\mathrm{b}
}^{*})}(g))}{2}}
  [A/gA\cap A]^{-1}\\ 
&\quad\times\sum_{a\in A/gA\cap
    A}\psi(\frac{1}{2}\beta(a,b)) \var(g^{-1}(b+a))\mbox{, }b\in
  B^{*}.
\end{split}
\end{equation}
On a alors le résultat suivant (comparer avec \cite[4.3.e]{pan-2001})~:
\begin{theo}\label{theo4.5.1}
Pour tout $g\in K_{B}$, on a
\begin{equation*}
  \mathcal{I}_{A,B}\circ S_{\psi}^{B}(g)=
  m(p_{Sp({\mathrm{b}
}^{*})}(g))q^{\frac{j(p_{Sp({\mathrm{b}
}^{*})}(g))}{2}}
M_{A}(g)\circ \mathcal{I}_{A,B}.
\end{equation*}
\end{theo}
\begin{dem}
Commençons par remarquer que $\mathcal{H}_{\psi}^{B}$ est l'espace
des fonctions $\var$ définies sur $W$, à valeurs dans l'espace des
fonctions numériques sur $B^{*}$, qui vérifient les relations
\begin{gather}
\var(w)(a+b)  = \psi(\frac{1}{2}\beta(b,a))\var(w)(b)\mbox{, }w\in
W\mbox{, }b\in B^{*}\mbox{, }a\in A\label{eq4.5.2}\\
\var(w)(b)  =  \psi(\frac{1}{2}\beta(b,w))\var(w+b)(0)\mbox{, }w\in
W\mbox{, }b\in B^{*}\label{eq4.5.3}\\
\int_{W}\sum_{B^{*}/A}\vert\var(w)(b)\vert^{2}\diff
w<+\infty.\label{eq4.5.4}
\end{gather}
En effet, la relation \ref{eq4.5.2} dit simplement que $\var$ est à
valeurs dans l'espace $\mathcal{H}_{\overline{\psi}}^{\mathrm{x}}$. Notons
$\widetilde{\rho}_{\overline{\psi}}^{\mathrm{x}}$ la représentation de
$H(B^{*})$ relevant la représentation
$\rho_{\overline{\psi}}^{\mathrm{x}}$ de $H({\mathrm{b}
}^{*})$ via le morphisme
$p_{H({\mathrm{b}
}^{*})}$ et soit $\var\in \mathcal{H}_{\psi}^{B}$. Alors,
il suit de la relation \ref{eq4.1.1} satisfaite par $\var$ que, l'on a
pour $b\in B^{*}$ et $w\in W$
\begin{eqnarray*}
  \var(w+b)(0)&=&\psi(\frac{1}{2}\beta(w,b))
  (\widetilde{\rho}_{\overline{\psi}}^{\mathrm{x}}(b)\var(w))(0)\\
&=&\psi(\frac{1}{2}\beta(w,b))\var(w)(b),
\end{eqnarray*}
qui est la relation \ref{eq4.5.3}. Enfin, la relation \ref{eq4.5.4}
dit simplement que la fonction $w\longmapsto\Vert\var(w)\Vert^{2}$
est de carré sommable sur $W$.

Réciproquement, supposons les relations \ref{eq4.5.3} et \ref{eq4.5.4}
vérifiées par la fonction $\var$ définie sur $W$ et à valeurs dans
l'espace $\mathcal{H}_{\overline{\psi}}^{\mathrm{x}}$. Alors, pour $w\in
W$ et $b,b'\in B^{*}$, on a
\begin{eqnarray*}
(\widetilde{\rho}_{\overline{\psi}}^{\mathrm{x}}(b)\var(w))(b')   
& =&\psi(\frac{1}{2}\beta(b',b))\var(w)(b+b')\\
& =&  \psi(\frac{1}{2}\beta(b',b))\psi(\frac{1}{2}\beta(b+b',w))
\var(w+b+b')(0)\\
& =&  \psi(\frac{1}{2}(\beta(b',b)+\beta(b+b',w)+\beta(w+b,b')))
\var(w+b)(b')\\
& =& \psi(\frac{1}{2}\beta(b,w))\var(w+b)(b')
\end{eqnarray*}
qui est exactement la relation \ref{eq4.1.1}. On a alors
$\sum_{B^{*}/A}\vert\var(w)(b)\vert^{2}=\Vert\var(w)\Vert^{2}$, $w\in W$, si bien
qu'il suit de la relation \ref{eq4.5.4} que la fonction
$w\longmapsto\Vert\var(w)\Vert^{2}$ est de carré sommable sur $W$ et
que $\var$ est bien un élément de $\mathcal{H}_{\psi}^{B}$.

Maintenant, montrons qu'il existe une fonction
$\var\in\mathcal{H}_{\psi}^{B}$ 
telle que $\var(0)\neq 0$. En effet, soit
$\var\in\mathcal{H}_{\psi}^{B}$ un élément non nul : il existe $w\in
W$ tel que $\var(w)\neq 0$. Appliquant à $\var$ l'opérateur
$\rho_{\psi}^{B}(-w)$ on se ramène au cas où $w=0$, i.e. $\var(0)\neq
0$. 

Soit $g\in K_{B}$. Comme $(g,S_{\psi}^{B}(g))$ (resp. $(g,M_{A}(g))$)
est un élément de $\widehat{Sp(W)}_{\psi}^{B}$
(resp. $\widehat{Sp(W)}_{\psi}^{A}$), il existe un nombre complexe
$\mu(g)$ tel que 
\begin{equation*}
\mathcal{I}_{A,B}\circ S_{\psi}^{B}(g)=\mu(g)
M_{A}(g)\circ \mathcal{I}_{A,B}.
\end{equation*}
D'après ce qu'on a vu plus haut et compte tenu de la relation
\ref{eq4.2.2}, il existe une fonction $\var\in\mathcal{H}_{\psi}^{B}$
telle que $S_{\psi}^{B}(g)\var(0)\neq 0$. Alors, compte tenu de la
définition de l'opérateur $\mathcal{I}_{A,B}$ (voir le corollaire
\ref{co4.1.1}) et des relations \ref{eq4.5.3}, \ref{eq4.2.2} et
\ref{eq4.5.1}, on a, pour tout $b\in B^{*}$,
\begin{eqnarray*}
\mathcal{I}_{A,B}\circ S_{\psi}^{B}(g)\var(b) & = &
S_{\psi}^{B}(g)\var(b)(0)\\ 
& = & S_{\psi}^{B}(g)\var(0)(b)\\ 
& = &
\widetilde{S}_{\overline{\psi}}(g)(\var(0))(b)\\ 
& = & m(p_{Sp({\mathrm{b}
}^{*})}(g))q^{\frac{j(p_{Sp({\mathrm{b}
}^{*})}(g))}{2}}
[A/gA\cap A]^{-1}\\ 
& &\quad\times\sum_{a\in A/gA\cap
  A}\psi(\frac{1}{2}\beta(a,b)) \var(0)(g^{-1}(b+a))\\
& = & m(p_{Sp({\mathrm{b}
}^{*})}(g))q^{\frac{j(p_{Sp({\mathrm{b}
}^{*})}(g))}{2}}
[A/gA\cap A]^{-1}\\ 
& &\quad\times\sum_{a\in A/gA\cap
  A}\psi(\frac{1}{2}\beta(a,b)) \var(g^{-1}(b+a))(0).
\end{eqnarray*}
D'autre part, il suit de la formule \ref{eq4.3.1} donnant l'expression
de l'opérateur $M_{A}(g)$, $g\in Sp(W)$, que, pour tout $b\in B^{*}$,
\begin{equation*}
\begin{split}
M_{A}(g)\!\circ\mathcal\!\mathcal{I}_{A,B}\var(b)
& = [A/gA\cap A]^{-1}
\sum_{a\in A/gA\cap A}\psi(\frac{1}{2}\beta(a,b))
\mathcal{I}_{A,B}\var(g^{-1}(b+a))\\
&= [A/gA\cap A]^{-1}
\sum_{a\in A/gA\cap A}\psi(\frac{1}{2}\beta(a,b))
\var(g^{-1}(b+a))(0),
\end{split}
 \end{equation*}
d'où le théorème.
\end{dem}

Pour $g\in K_{B}$, on définit l'opérateur $S_{\psi}^{A}(g)$ dans
l'espace $\mathcal{H}_{\psi}^{A}$ en posant
\begin{equation}\label{eq4.5.0}
S_{\psi}^{A}(g)=m(p_{Sp({\mathrm{b}
}^{*})}(g))q^{\frac{j(p_{Sp({\mathrm{b}
}^{*})}(g))}{2}}
M_{A}(g)
\end{equation}
Il est clair d'après le théorème \ref{theo4.5.1} et le lemme
\ref{lem4.2.1} que l'application $g\longmapsto(g,S_{\psi}^{A}(g))$ est
un homomorphisme de groupes de $K_{B}$ dans
$\widehat{Sp(W)}_{\psi}^{A}$, scindant la suite exacte
\ref{eq3.2.2}.

Rappelons que $P_{B}$ est le stabilisateur du lagrangien $\mathrm{x}$ pour
l'action de $K_{B}$ dans le $\Bbb{F}_{q}$-espace symplectique ${\mathrm{b}
}^{*}$. On désigne par $\zeta$ le caractère de $P_{B}$ tel que
\begin{equation}\label{eq4.5.5}
\zeta(p)=\left(\frac{\det_{\mathrm{x}}p_{Sp({\mathrm{b}
}^{*})}(p)}{q}\right),
\end{equation}
où $\left(\frac{\,}{q}\right)$ désigne le symbole de Legendre relatif au
corps $\Bbb{F}_{q}$~: pour tout $a\in \Bbb{F}_{q}^{\times}$
\begin{equation}\label{eq4.5.6}
\left(\frac{a}{q}\right)=
\begin{cases}
  1&\mbox{ si
  }a\in(\Bbb{F}_{q}^{\times})^{2}\\
 -1&\mbox{ si }a\notin(\Bbb{F}_{q}^{\times})^{2}.
\end{cases}
\end{equation}

On se donne une base autoduale
$(e_{1},\ldots,e_{r},f_{1},\ldots,f_{r})$ de $W$ telle que $B =
\varpi\O e_1 \oplus \cdots \oplus \varpi\O e_l\oplus\O e_{l+1} \oplus
\cdots \oplus \O e_r \oplus \O f_1 \oplus \cdots \oplus \O f_r$ et $ A
= \O e_{1}\oplus\cdots\oplus \O e_{r}\oplus\O
f_{1}\oplus\cdots\oplus\O f_{r}$ et on reprend les notations des
paragraphes \ref{2.6} et \ref{3.3}. En particulier, à partir d'ici $\omega$
désigne l'indice de Weil relatif au caractère
$\frac{1}{2}\overline{\psi}$ du corps $\mathbb{F}_{q}$. On a alors le
résultat suivant~:
\begin{co}\label{co4.5.1}
La représentation  $S_{\psi}^{A}$ de $K_{B}$ vérifie les relations
suivantes qui la déterminent entièrement~:
\begin{eqnarray}
S_{\psi}^{A}(p) &=& \zeta(p)M_{A}(p)\mbox{, }p\in
P_{B}\label{eq4.5.7}\\ 
S_{\psi}^{A}(\varsigma_{B}) &=&\left(\frac{-1}{q}\right)^{l}
\omega(1)^{-l}q^{\frac{l}{2}}M_{A}(\varsigma_{B}).\label{eq4.5.8}
\end{eqnarray}
\end{co}
\begin{dem}
Remarquons tout d'abord que $j(p_{Sp({\mathrm{b} }^{*})}(p))=0$, $p\in
P_{B}$, tandis que, avec les notations du paragraphe \ref{3.3},
$p_{Sp({\mathrm{b} }^{*})}(\varsigma_{B})=J_{{\mathrm{b}
  }^{*}}=-\tau_{\{1,\ldots,l\}}$ et $j(p_{Sp({\mathrm{b}
  }^{*})}(\varsigma_{B}))=j(J_{{\mathrm{b} }^{*}})=l$.  Il suit alors
de la définition de la représentation $S_{\psi}^{A}$ de $K_{B}$ et des
propriétés de la fonction $\theta$ énoncées au paragraphe \ref{3.3}
\begin{eqnarray*}
  S_{\psi}^{A}(p) &=&
  \omega(\sideset{}{_{\mathrm{x}}}\det(p_{Sp({\mathrm{b}
    }^{*})}(p))^{-1} \omega(1)M_{A}(p)\mbox{, }p\in
  P_{B}\\ S_{\psi}^{A}(\varsigma_{B}) &=&\left(\frac{-1}{q}\right)^{l}
  \omega(1)^{-l}q^{\frac{l}{2}}M_{A}(\varsigma_{B})
\end{eqnarray*}
Or, d'après \cite[Proposition A.9 (1)]{rao-1993}, on a
$\omega(a)^{-1}\omega(1)=\left(\frac{a}{q}\right)$, pour
$a\in\Bbb{F}_{q}$. Les formules \ref{eq4.5.7} et \ref{eq4.5.8}
résultent de ces considérations.  Le lemme \ref{lem2.6.1} entraîne que
ces formules déterminent entièrement la représentation $S_{\psi}^{A}$
de $K_{B}$.
\end{dem}

\subsection{}\label{4.6}
On garde les notations du paragraphe précédent. Il suit du théorème
\ref{theo4.5.1} que démontrer le théorème \ref{theo4.2.1} revient à
démontrer le
\begin{theo}\label{theo4.6.1}
  L'application $s_{B}:g\longmapsto(g,S_{\psi}^{A}(g))$ est un homomorphisme
  de grou\-pes de $K_{B}$ dans le groupe métaplectique de $W$.
\end{theo}
\begin{dem}
Selon le corollaire \ref{co2.4.1}, le groupe $K_{B}$ est engendré par
$P_{B}$ et $\overline{N}_{B}$. Or, il suit de la proposition
\ref{pr4.3.1} et du corollaire \ref{co4.5.1} que, pour tout $p\in
P_{B}$, $s_{B}(p)$ est un élément de $Mp(W)$. Il nous reste à démontrer
que, pour tout $n\in\overline{N}_{B}$, $s_{B}(n)\in Mp(W)$.

On a $W=X\oplus Y$ où $X$ (resp. $Y$) désigne le sous-espace
lagrangien de $W$ engendré par les vecteurs $e_{1},\ldots,e_{r}$
(resp. $f_{1},\ldots,f_{r}$). Par suite, l'application
$\var\longmapsto\var_{\vert X}$ permet d'identifier l'espace
$\mathcal{H}_{\psi}^{Y}$ de la représentation $\rho_{\psi}^{Y}$ avec
$L^{2}(X,\diff x)$ et l'opérateur
$\mathcal{I}_{Y,A}:\mathcal{H}_{\psi}^{A}\longrightarrow\mathcal{H}_{\psi}^{Y}
=L^{2}(X,\diff x)$ entrelaçant les représentations $\rho_{\psi}^{A}$
et $\rho_{\psi}^{Y}$ est donné par
\begin{equation}\label{eq4.6.1}
\mathcal{I}_{Y,A}\var(x)=\int_{Y/A\cap Y}\psi(\frac{1}{2}\beta(y,x))
\var(x+y)\diff\dot{y}\mbox{, }\var\in\mathcal{H}_{\psi}^{A}\mbox{,
}x\in X.
\end{equation}
Il est immédiat que $1_{A}$, la fonction caractéristique du réseau
$A$, est un élément de $\mathcal{H}_{\psi}^{A}$. De plus, comme
$A=A\cap X\oplus A\cap Y$, on a $\mathcal{I}_{Y,A}1_{A}=1_{A\cap X}$,
où $1_{A\cap X}$ est la fonction caractéristique du réseau $A\cap X$
de $X$.

Comme $\overline{N}_{B}$ est un sous-groupe du radical unipotent
$N_{Y}$ du stabilisateur de $Y$ dans $Sp(W)$, lequel se relève de
manière unique en un sous-groupe du groupe métaplectique (voir le
paragraphe \ref{3.6}),
pour établir notre résultat, il suffit de montrer que pour tout
$n\in\overline{N}_{B}$, on a
\begin{equation}\label{eq4.6.2}
\mathcal{I}_{Y,A}\circ S_{\psi}^{A}(n)=
S_{\psi}^{Y}(n)\circ\mathcal{I}_{Y,A}.
\end{equation}
Or, il est clair que les deux membres de
l'équation \ref{eq4.6.2} diffèrent d'un facteur scalaire. Il suffit
donc de montrer que pour tout $n\in\overline{N}_{B}$, on a
\begin{equation}\label{eq4.6.3}
\mathcal{I}_{Y,A}\circ S_{\psi}^{A}(n)1_{A}=S_{\psi}^{Y}(n)1_{A\cap X}.
\end{equation}
Mais, tout élément de $\overline{N}_{B}$ est de la forme $y_{B}(a)$
avec $a\in \mathcal{M}_{l}(\O)$ une matrice symétrique. Soit donc $a$
une telle matrice. Pour ce qui est du membre de droite de
\ref{eq4.6.3}, il résulte du lemme \ref{lem3.6.1} que l'on a
\begin{equation}\label{eq4.6.4}
  S_{\psi}^{Y}(y_{B}(a))1_{A\cap X}(x)=
\psi(-\frac{\varpi^{\lambda_{\psi}-1}}{2}\transposee{x}ax)1_{A\cap
  X}(x)\mbox{, }x\in X.
\end{equation}

Explicitons le membre de gauche de \ref{eq4.6.3}. Remarquons tout
d'abord que $y_{B}(a)=\varsigma_{B}x_{B}(-a)\varsigma_{B}^{-1}$ de sorte que
\begin{equation*}
S_{\psi}^{A}(y_{B}(a))=S_{\psi}^{A}(\varsigma_{B})S_{\psi}^{A}(x_{B}(-a))
S_{\psi}^{A}(\varsigma_{B}^{-1}).
\end{equation*}
Or $\varsigma_{B}^{4}=Id$, $\varsigma_{B}^{2}\in P_{B}$ et
$\zeta(\varsigma_{B}^{2})=\left(\frac{-1}{q}\right)^{l}$. On en déduit que 
\begin{equation}\label{eq4.6.5}
S_{\psi}^{A}(y_{B}(a))1_{A}=\left(\frac{-1}{q}\right)^{l}
S_{\psi}^{A}(\varsigma_{B})S_{\psi}^{A}(x_{B}(-a))
S_{\psi}^{A}(\varsigma_{B})1_{A}
\end{equation}
D'autre part 
\begin{multline*}
\varsigma_{B}A=\varpi\O e_{1}\oplus\cdots\oplus\varpi\O e_{l}
\oplus\O e_{l+1}\oplus\cdots\oplus\O e_{r}\\
\oplus\varpi^{-1}\O f_{1}\oplus\cdots\oplus\varpi^{-1}\O f_{l}
\oplus\O f_{l+1}\oplus\cdots\oplus\O f_{r}
\end{multline*}
de sorte que,  
\begin{equation*}
\varsigma_{B}A/A\cap\varsigma_{B}A= \Bbb{F}_{q}^{l}.
\end{equation*}

Si $u\in\O^{l}$, on
désigne par $\dot{u}$ son image dans $\Bbb{F}_{q}^{l}$ et on pose
$u.e=u_{1}e_{1}+\cdots+u_{l}e_{l}$ et
$u.f=u_{1}f_{1}+\cdots+u_{l}f_{l}$.

Comme $1_{A}=\delta_{0}$, il suit du corollaire
\ref{co4.5.1} et du lemme \ref{lem4.4.1} que l'on a
\begin{equation}\label{eq4.6.6}
S_{\psi}^{A}(\varsigma_{B})1_{A}=q^{-\frac{l}{2}}
\left(\frac{-1}{q}\right)^{l}\omega(1)^{-l}
\sum_{\dot{u}\in\Bbb{F}_{q}^{l}}\delta_{\varpi^{-1}u.f}.
\end{equation}
Or, si $u\in\O^{l}$, il résulte du corollaire \ref{co4.5.1}, de la
proposition \ref{pr4.3.1}, du lemme \ref{lem4.4.1} et de la relation
\ref{eq4.4.2} que l'on a
\begin{eqnarray*}
S_{\psi}^{A}(x_{B}(-a))\delta_{\varpi^{-1}u.f} & = &
\delta_{x_{B}(-a)\varpi^{-1}u.f}\\
 & = & \delta_{-au.e+\varpi^{-1}u.f}\\
 & =
&\psi(-\frac{1}{2}\beta(au.e,\varpi^{-1}u.f))\delta_{\varpi^{-1}u.f}\\
& = &\psi(-\frac{\varpi^{\lambda_{\psi}-1}}{2}\transposee{u}au)\delta_{\varpi^{-1}u.f}.
\end{eqnarray*}
Compte tenu de \ref{eq4.6.6}, il vient
\begin{equation*}
S_{\psi}^{A}(x_{B}(-a))S_{\psi}^{A}(\varsigma_{B})1_{A}=
q^{-\frac{l}{2}}\left(\frac{-1}{q}\right)^{l}\omega(1)^{-l}
\sum_{\dot{u}\in\Bbb{F}_{q}^{l}}\psi(-\frac{\varpi^{\lambda_{\psi}-1}}{2}\transposee{u}au)
\delta_{\varpi^{-1}u.f}.
\end{equation*}
Utilisant à nouveau le corollaire \ref{co4.5.1}, le lemme
\ref{lem4.4.1}, la formule \ref{eq4.4.2}, et compte tenu de
\ref{eq4.6.5} et de la relation
$\left(\frac{-1}{q}\right)=\omega(1)^{2}$ (voir \cite[Proposition A.9
  (1) et Theorem A.2 (2)]{rao-1993}), on obtient
\begin{equation}\label{eq4.6.7}
\begin{split}
S_{\psi}^{A}(y_{B}(a))1_{A} &=
\left(\frac{-1}{q}\right)^{l}\omega(1)^{-2l}q^{-l}
\sum_{\dot{u},\dot{v}\in\Bbb{F}_{q}^{l}}
\psi(-\frac{\varpi^{\lambda_{\psi}-1}}{2}\transposee{u}au)\\
& \qquad\times
\psi(\frac{1}{2}\beta(\varpi^{-1}u.f,-v.e))\delta_{u.e+\varpi^{-1}v.f}\\
&= q^{-l}\sum_{\dot{v}\in\Bbb{F}_{q}^{l}}\left(\sum_{\dot{u}\in\Bbb{F}_{q}^{l}}(
\psi(-\frac{\varpi^{\lambda_{\psi}-1}}{2}\transposee{u}au)
\psi(\varpi^{\lambda_{\psi}-1}\transposee{u}v)\right)\delta_{\varpi^{-1}v.f}
\end{split}
\end{equation}

Maintenant, nous allons calculer l'action de l'opérateur
d'entrelacement $\mathcal{I}_{Y,A}$ sur les fonctions $\delta_{v}$,
$v\in W$. Soit $v\in W$. \'Ecrivons $v=v_{X}+v_{Y}$ avec $v_{X}\in X$
et $v_{Y}\in Y$. Utilisant la formule \ref{eq4.6.1} on voit que, si
$x\notin v_{X}+A\cap X$, on a
\begin{equation*}
\mathcal{I}_{Y,A}\delta_{v}(x)=0
\end{equation*}
et, si $x\in v_{X}+A\cap X$,
\begin{equation*}
\begin{split}
\mathcal{I}_{Y,A}\delta_{v}(x) &=
\delta_{v}(x+v_{Y})\psi(\frac{1}{2}\beta(v_{Y},x))\\
& = \psi(\frac{1}{2}\beta(v,x-v_{X}))\psi(\frac{1}{2}\beta(v_{Y},x))\\
& = \psi(\frac{1}{2}\beta(v_{X},v_{Y}))\psi(\beta(v_{Y},x)).
\end{split}
\end{equation*}
On voit donc que, si $v\in\O^{l}$, on a 
\begin{equation}
\mathcal{I}_{Y,A}\delta_{\varpi^{-1}v.f}(x) =
\begin{cases}
0 &\mbox{ si }x\notin A\cap X\\
\psi(-\varpi^{\lambda_{\psi}-1}\transposee{v}x) &\mbox{ si }x\in A\cap X.
\end{cases}
\end{equation}
Compte tenu de la relation \ref{eq4.6.7}, il vient, pour $x\notin
A\cap X$,
\begin{equation*}
\mathcal{I}_{Y,A}S_{\psi}^{A}(y_{B}(a))1_{A}(x)=0
\end{equation*}
et, pour $x\in A\cap X$,
\begin{equation*}
\begin{split}
\mathcal{I}_{Y,A}S_{\psi}^{A}(y_{B}(a))1_{A}(x) &=
q^{-l}\sum_{\dot{u}\in\Bbb{F}_{q}^{l}}
\psi(-\frac{\varpi^{\lambda_{\psi}-1}}{2}\transposee{u}au)
\sum_{\dot{v}\in\Bbb{F}_{q}^{l}}
\psi(\varpi^{\lambda_{\psi}-1}\transposee{v}(u-x))
\\
&=\psi(-\frac{\varpi^{\lambda_{\psi}-1}}{2}\transposee{x}ax)
\end{split}
\end{equation*}
Comparant avec \ref{eq4.6.4}, on voit que la relation \ref{eq4.6.3}
est vraie. Ceci achève la démonstration du théorème.
\end{dem}
\begin{rem}
(i) En utilisant les résultats de \cite{bruhat-tits-1972} et
  \cite{bruhat-tits-1987} et en raisonnant comme dans la démonstration
  de \cite[Lemma 11.1]{moore-1968}, on peut montrer que, si $q\geq4$,
  le groupe $K_{B}$ est parfait, i.e. égal à son sous-groupe des
  commutateurs. Il s'ensuit que, dans ce cas, la suite exacte
  \ref{eq3.2.2} admet un unique scindage au dessus de $K_{B}$, lequel
  est à valeurs dans $Mp(W)$.

(ii) Désormais, nous utilisons la section $s_B$ pour identifier $K_{B}$
avec un sous-groupe de $Mp(W)$. Elle est caractérisée par le fait que
la représentation $S_{\psi}^{A}\circ s_{B}$ de $K_{B}$ est donnée par
la formule \ref{eq4.5.0}. Avec cette convention, la restriction
de la représentation métaplectique $S_{\psi}^{A}$ à $K_{B}$ est
entièrement déterminée par les formules \ref{eq4.5.7} et
\ref{eq4.5.8}.
\end{rem}

\subsection{}\label{4.7}

On se donne une décomposition $W=W_{1}\oplus\cdots\oplus W_{n}$ de $W$
en somme directe orthogonale de sous-espaces symplectiques tous non nuls.

Les groupes de Heisenberg $H(W_{i})$, $1\leq i\leq n$ sont des
sous-groupes qui commutent deux à deux du groupe $H(W)$ et
l'application $(h_{1},\ldots,h_{n})\mapsto h_{1}\cdots h_{n}$ est un
morphisme surjectif du produit direct $H(W_{1})\times\cdots\times
H(W_{n})$ sur $H(W)$.

Pour $1\leq i\leq n$, soit $A_{i}$ un réseau autodual de $W_{i}$. Alors
$A=A_{1}\oplus\cdots\oplus A_{n}$ est un réseau autodual de $W$. Pour
$1\leq i\leq n$, soit $\varphi_{i}\in\mathcal{H}_{\psi}^{A_{i}}$. Alors
la fonction $\varphi_{1}\otimes\cdots\otimes\varphi_{n}$ définie sur $W$ par
la relation
\begin{equation}\label{eq4.7.1}
\varphi_{1}\otimes\cdots\otimes\varphi_{n}(w_{1}+\cdots+w_{n})=
\varphi_{1}(w_{1})\cdots\varphi_{n}(w_{n})\mbox{, }w_{i}\in
W_{i}\mbox{, }1\leq i\leq n
\end{equation}
est un élément de $\mathcal{H}_{\psi}^{A}$. De plus, l'application
$(\varphi_{1},\ldots,\varphi_{n})\mapsto
\varphi_{1}\otimes\cdots\otimes\varphi_{n}$
induit un isomorphisme d'espaces
de Hilbert $\mathcal{I}_{A}$ du produit tensoriel hilbertien
$\mathcal{H}_{\psi}^{A_{1}}\widehat{\otimes}\cdots\widehat{\otimes}
\mathcal{H}_{\psi}^{A_{n}}$ sur $\mathcal{H}_{\psi}^{A}$.

Le lemme suivant est bien connu et en tout cas immédiat à démontrer.
\begin{lem}\label{lem4.7.1}
  La représentation
  $\rho_{\psi}^{A_{1}}\otimes\cdots\otimes\rho_{\psi}^{A_{n}}$ du
  produit direct $H(W_{1})\times\cdots\times H(W_{n})$ passe au
  quotient en une représentation de $H(W)$. De plus $\mathcal{I}_{A}$
  entrelace les représentations
  $\rho_{\psi}^{A_{1}}\otimes\cdots\otimes\rho_{\psi}^{A_{n}}$ et
  $\rho_{\psi}^{A}$ de $H(W)$.
\end{lem}

Soit $1\leq i\leq n$. On désigne par $\gamma_{i}$ l'injection
naturelle de $Sp(W_{i})$ dans $Sp(W)$ : si $x\in Sp(W_{i})$,
$\gamma_{i}(x)$ laisse stable $W_{i}$ en y agissant comme $x$, et
opère trivialement sur $W_{i}^{\perp}$. Alors $\gamma_{i}$ se relève
de manière unique en un homomorphisme injectif noté
$\widetilde{\gamma}_{i}$ de $Mp(W_{i})$ dans $Mp(W)$. On identifie
$Sp(W_{i})$ (resp.  $Mp(W_{i})$) avec son image dans $Sp(W)$
(resp. $Mp(W)$) par $\gamma_{i}$ (resp. $\widetilde{\gamma}_{i}$). Les
sous-groupes $Sp(W_{i})$ (resp. $Mp(W_{i})$), $1\leq i\leq n$ de
$Sp(W)$ (resp. $Mp(W)$) commutent deux à deux. On désigne par
$Sp(W_{1})\cdots Sp(W_{n})$ (resp. $Mp(W_{1})\cdots Mp(W_{n})$) leur
produit dans $Sp(W)$ (resp. $Mp(W)$), qui n'est autre que le
sous-groupe qu'ils engendrent. Evidemment, la multiplication induit un
isomorphisme (resp.  morphisme surjectif) du produit direct
$Sp(W_{1})\times\cdots\times Sp(W_{n})$ (resp.
$Mp(W_{1})\times\cdots\times Mp(W_{n})$) sur $Sp(W_{1})\cdots Sp(W_{n})$
(resp. $Mp(W_{1})\cdots Mp(W_{n})$).

\begin{lem}\label{lem4.7.2}
La représentation $S_{\psi}^{A_{1}}\otimes\cdots\otimes
S_{\psi}^{A_{n}}$ du produit direct $Mp(W_{1})\times\cdots\times
Mp(W_{n})$ passe au quotient en une représentation du sous-groupe
$Mp(W_{1})\cdots Mp(W_{n})$ de $Mp(W)$ notée de même. De plus
$\mathcal{I}_{A}$ entrelace les représentations
$S_{\psi}^{A_{1}}\otimes\cdots\otimes S_{\psi}^{A_{n}}$ et
$(S_{\psi}^{A})_{\vert Mp(W_{1})\cdots Mp(W_{n})}$.
\end{lem}
\begin{dem}
Le premier point est immédiat. Le deuxième point est une conséquence
facile du lemme \ref{lem4.7.1} et du fait que les groupes $Mp(W_{i})$
sont parfaits.
\end{dem}

Maintenant, pour $1\leq i\leq n$,  on se donne un bon réseau $B_{i}$
de $W_{i}$. Alors $B=B_{1}\oplus\cdots\oplus B_{n}$ est un bon réseau
de $W$. On peut considérer les groupes $K_{B_{i}}$, $1\leq i\leq n$
comme des sous-groupes de $Sp(W)$ et comme ils commutent deux à deux
former leur produit $K_{B_{1}}\cdots K_{B_{n}}$, qui est clairement un
sous-groupe de $K_{B}$. Pour $1\leq i\leq n$, on considère
$s_{B_{i}}$ comme une section de $K_{B_{i}}$ à valeurs dans le
sous-groupe $Mp(W_{i})$ de $Mp(W)$. 
\begin{pr}\label{pr4.7.1}
On a l'égalité des sections $s_{B}$ et $s_{B_{1}}\cdots s_{B_{n}}$ au
dessus du sous-groupe $K_{B_{1}}\cdots K_{B_{n}}$.
\end{pr}
\begin{dem}
Pour $1\leq i\leq n$, soit $A_{i}$ un réseau autodual de $W_{i}$
vérifiant $B_{i}\subset A_{i}\subset B_{i}^{*}$. Alors
$A=A_{1}\oplus\cdots\oplus A_{n}$ est un réseau autodual de $W$
vérifiant $B\subset A\subset B^{*}$.  Il suffit alors de démontrer que
l'opérateur $\mathcal{I}_{A}$ du lemme \ref{lem4.7.1} entrelace les deux
représentations $(S_{\psi}^{A_{1}}\circ s_{B_{1}})\otimes\cdots\otimes
(S_{\psi}^{A_{n}}\circ s_{B_{n}})$ et $(S_{\psi}^{A}\circ s_{B})_{\vert
  K_{B_{1}}\cdots K_{B_{n}}}$ de $K_{B_{1}}\cdots K_{B_{n}}$. Ceci se
vérifie par un calcul direct à partir de la formule \ref{eq4.5.0}
qui décrit explicitement la représentation $S_{\psi}^{A}\circ s_{B}$.
\end{dem}

\section{Restriction de la représentation de Weil à un
  sous-groupe compact maximal.}\label{5}
Soit $B$ un bon réseau et $l=l(B)$. Dans cette section, nous
déterminons la décomposition en composantes irréductibles de la
restriction de la représentation métaplectique à $K_{B}$. On reprend
les notations de la section précédente.
\subsection{}\label{5.1}
Soit $A$ un réseau autodual tel que $\varpi B^{*}\subset B\subset
A\subset B^{*}$. On réalise la restriction de la représentation
métaplectique, notée $S_{\psi}^{A}$, à $K_{B}$ dans l'espace
$\mathcal{H}_{\psi}^{A}$ par les formules \ref{eq4.5.7} et
\ref{eq4.5.8} du corollaire \ref{co4.5.1}.

Si $\mathcal{E}$ est un sous-espace de $\mathcal{H}_{\psi}^{A}$, on
désigne par $\mathcal{E}^{+}$ (resp. $\mathcal{E}^{-}$) le sous-espace
de $\mathcal{E}$ constitué des fonctions paires (resp. impaires). On
sait que les sous-espaces $\mathcal{H}_{\psi}^{A,\pm}$ sont invariants
et irréductibles sous l'action de la représentation $S_{\psi}^{A}$ de
$Mp(W)$.

On considère les sous-espaces suivants de
$\mathcal{H}_{\psi}^{A}$~: 
\begin{eqnarray*}
\mathcal{E}^{B}_{0,0}&=&\Bbb{C}\delta_{0}\\
\mathcal{E}^{B}_{n,0}&=&\{f\in\mathcal{H}_{\psi}^{A}\vert
\Supp f\subset\varpi^{-n}A\backslash\varpi^{-n}B\}\mbox{,
}n\in\Bbb{N}\mbox{, }n\geq 1\\
\mathcal{E}^{B}_{n,2}&=&\{f\in\mathcal{H}_{\psi}^{A}\vert
\Supp f\subset\varpi^{-n}B^{*}\backslash\varpi^{-n}A\}\mbox{,
}n\in\Bbb{N}\\
\mathcal{E}^{B}_{n,1}&=&\{f\in\mathcal{H}_{\psi}^{A}\vert
\Supp f\subset\varpi^{-(n+1)}B\backslash\varpi^{-n}B^{*}\}\mbox{,
}n\in\Bbb{N}
\end{eqnarray*}

Pour $n\in\Bbb{N}$, soit $K_{n}$ le sous-groupe des éléments $g\in
K=K_{A}$ tels que $(g-1)A\subset\varpi^{n} A$. Alors, on a $K_{0}=K$
et, pour $n\geq1$, $K_{n}\subset P_{B}$. 

Soit $n\in\Bbb{N}$, $n\geq1$. Alors le groupe quotient
$\varpi^{-n}A/A$ est naturellement muni d'une structure de
$O_{n-1}$-module symplectique pour la forme bilinéaire alternée
$\beta_{n}$ définie par
\begin{equation*}
  \beta_{n}(v+A,w+A)=
\varpi^{-\lambda_{\psi}+2n}\beta(v,w)+\varpi^{n}\O\mbox{, }v,w\in\varpi^{-n}A.
\end{equation*}
Il est isomorphe, via l'application $v+A\longmapsto
p_{\mathrm{a}_{n}}(\varpi^{n}v)$, au $O_{n-1}$-module symplectique
${\mathrm{a}}_{n}=A/\varpi^{n}A$, défini au paragraphe \ref{2.5}. Le groupe $K$
agit naturellement dans $\varpi^{-n}A/A=\mathrm{a}_{n}$ et cette action
induit un morphisme surjectif de groupes de $K$ sur le groupe
symplectique $Sp({\mathrm{a}_{n}})$ dont le noyau est $K_{n}$.

Soit $v\in\varpi^{-n}A\backslash\varpi^{-n+1}A$. On désigne par
$P_{B,v}$ le stabilisateur dans le groupe $P_{B}$ de
$v+A\in\varpi^{-n}A/A$ et par $\tilde{P}_{B,v}$ le sous-groupe $\{\pm
Id\}P_{B,v}$ de $P_{B}$. Il est clair que $K_{n}$ est un sous-groupe
de $P_{B,v}$.  Il résulte de \ref{eq4.4.2} et \ref{eq4.4.3} que la
fonction $\chi_{B,v}$ définie sur $P_{B,v}$ en posant
\begin{equation*}
\chi_{B,v}(g)=\psi(\frac{1}{2}\beta(gv,v))\mbox{, }g\in P_{B,v},
\end{equation*}
est un caractère et que l'on a 
\begin{equation*}
  M_{A}(g)\delta_{v}=\chi_{B,v}(g)\delta_{v}\mbox{, }g\in P_{B,v}.
\end{equation*}
Le caractère $\chi_{B,v}$ de $P_{B,v}$ s'étend en les deux seuls
caractères $\chi^{\pm}_{B,v}$ de $\tilde{P}_{B,v}$ tels que
$\chi^{\pm}_{B,v}(-Id)=\pm1$.
Enfin, on désigne par $\chi_{v}$ la restriction du caractère
$\chi_{B,v}$ à $K_{n}$.
\begin{lem}\label{lem5.1.1}
Soit $n\in\Bbb{N}$, $n\geq1$.

(i) Soit $v,v'\in\varpi^{-n}A\backslash\varpi^{-n+1}A$. Alors, on a
$\chi_{v}=\chi_{v'}$ si et seulement si $v'\in\pm v+A$.

(ii) On a
\begin{equation*}
K_{2n}\subset\bigcap_{v\in C}\ker\chi_{v}
\subsetneq K_{2n-2}
\end{equation*}
où $C$ désigne $\varpi^{-n+1}B^{*}\backslash\varpi^{-n+1}A$,
$\varpi^{-n}A\backslash\varpi^{-n}B$ ou
$\varpi^{-n}B\backslash\varpi^{-n+1}B^{*}$.
\end{lem}
\begin{dem}
On se donne une base autoduale
$(e_{1},\ldots,e_{r},f_{1},\ldots,f_{r})$ de $W$ telle que $A=\O e_1
\oplus \cdots \oplus \O e_r \oplus \O f_{1}\cdots\oplus\O f_{r}$ et
$B=\varpi\O e_1 \oplus \cdots \oplus\varpi\O e_l \oplus\O e_{l+1} \oplus
\cdots \oplus \O e_r\oplus \O f_{1}\cdots\oplus\O f_{r}$. On identifie les
éléments (resp. les endomorphismes) de $W$ avec le vecteur colonne de
leurs coordonnées (resp. leur matrice) dans cette base. Si $\lambda$
est un élément de $k^{r}$, on pose
$\Vert\lambda\Vert_{p}=\max\{\vert\lambda_{i}\vert_{p}\,\vert\, 1\leq
i\leq r\}$, où $\vert\,\vert_{p}$ désigne la valeur absolue normalisée
sur le corps $k$.

Alors, les éléments de $K_{n}$ sont ceux de $Sp(W)$
s'écrivant $\left(
\begin{smallmatrix}
Id+\varpi^{n}u & \varpi^{n}x\\
\varpi^{n}y & Id+\varpi^{n}z
\end{smallmatrix}\right)
$ avec $u,x,y,z$ des matrices d'ordre $r$ à coefficients dans $\O$ et
vérifiant
\begin{equation}\label{eq5.1.1}
\begin{split}
\transposee{x}-x+\varpi^{n}(\transposee{x}z-\transposee{z}x)&=0\\
y-\transposee{y}+\varpi^{n}(\transposee{u}y-\transposee{y}u)&=0\\
z+\transposee{u}+\varpi^{n}(\transposee{u}z-\transposee{y}x)&=0.
\end{split}
\end{equation}
De même, les éléments de $\varpi^{-n}A\backslash\varpi^{-n+1}A$ sont
ceux de $W$  s'écrivant $\varpi^{-n}\left(
\begin{smallmatrix}
\lambda\\\mu
\end{smallmatrix}\right)
$ avec $\lambda,\mu\in k^{r}$ tels
$\max\{\Vert\lambda\Vert_{p},\Vert\mu\Vert_{p}\}=1$. 

(i) Un calcul immédiat montre que si $v=\varpi^{-n}\left(
\begin{smallmatrix}
\lambda\\ \mu
\end{smallmatrix}\right)$ est un élément de $
\varpi^{-n}A\backslash\varpi^{-n+1}A$ et $g=\left(
\begin{smallmatrix}
Id+\varpi^{n}u & \varpi^{n}x\\
\varpi^{n}y & Id+\varpi^{n}z
\end{smallmatrix}\right)$ est un élément de $ K_{n}$, on a 
\begin{equation}\label{eq5.1.2}
  \chi_{v}(g)=
\psi(\frac{1}{2}\varpi^{\lambda_{\psi}-n}
(\transposee{\lambda}\transposee{y}\lambda+
\transposee{\lambda}(\transposee{u}-z)\mu+
\transposee{\mu}x\mu)).
\end{equation}
On en déduit facilement que si $v=\varpi^{-n}\left(
\begin{smallmatrix}
\lambda\\\mu
\end{smallmatrix}\right)$ et $v'=\varpi^{-n}\left(
\begin{smallmatrix}
\lambda'\\ \mu'
\end{smallmatrix}\right)$ sont deux éléments de
$\varpi^{-n}A\backslash\varpi^{-n+1}A$, on a $\chi_{v}=\chi_{v'}$, si et
seulement si pour toute paire $(x,z)$ de matrices à coefficients dans
$\O$ avec $x$ symétrique, on a 
\begin{eqnarray*}
\transposee{\lambda}\transposee{x}\lambda-
\transposee{\lambda}'\transposee{x}\lambda'&\in&\varpi^{n}\O\\
\transposee{\mu}x\mu-\transposee{\mu}'x\mu'&\in&\varpi^{n}\O\\
\transposee{\lambda}z\mu-\transposee{\lambda}'z\mu'&\in&\varpi^{n}\O.
\end{eqnarray*}
L'assertion (i) du lemme en découle. 

(ii) Il est immédiat que si
$v\in\varpi^{-n}A\backslash\varpi^{-n+1}A$,
$K_{2n}\subset\ker\chi_{v}$. Ceci montre la première inclusion de
l'assertion (ii).

Soit $m\leq s$ deux entiers naturels. Dans la suite, on écrit les
vecteurs colonnes $\lambda\in k^{s}$ (resp. les matrices
$a\in\mathcal{M}_{s}(k)$) sous la forme $\lambda=\left(
\begin{smallmatrix}
\lambda_{1}\\\lambda_{2}
\end{smallmatrix}\right)
$ (resp. $a=\left(
\begin{smallmatrix}
a_{11}&a_{12}\\a_{21}&a_{22}
\end{smallmatrix}\right)
$) avec $\lambda_{1}\in k^{m}$ et $\lambda_{2}\in k^{s-m}$
(resp. $a_{11}\in\mathcal{M}_{m}(k)$, $a_{12}\in\mathcal{M}_{m,
  s-m}(k)$, $a_{21}\in\mathcal{M}_{s-m,m}(k)$ et
$a_{22}\in\mathcal{M}_{s-m}(k)$).  L'assertion suivante est
immédiate~: on suppose que la matrice $a\in\mathcal{M}_{s}(\O)$ est
telle que
\begin{equation}\label{eq5.1.3}
\transposee{a}-a\in\varpi^{n}\mathcal{M}_{s}(\O)
\end{equation} 
et que la forme quadratique
$Q_{a}(\lambda)=\transposee{\lambda}a\lambda$ vérifie
\begin{equation}\label{eq5.1.4}
Q_{a}(\lambda)\in\varpi^{n}\O\mbox{, pour tout }
\lambda\in\O^{s}\mbox{ tel que }\Vert\lambda_{1}\Vert_{p}=1\mbox{ et
}\Vert\lambda_{2}\Vert_{p}<1. 
\end{equation}
Alors, on a $a\in\varpi^{n-2}\mathcal{M}_{s}(\O)$ et
$a_{11}\in\varpi^{n}\mathcal{M}_{m}(\O)$.

Soit  $g=\left(
\begin{smallmatrix}
Id+\varpi^{n}u & \varpi^{n}x\\
\varpi^{n}y & Id+\varpi^{n}z
\end{smallmatrix}\right)\in K_{n}$. D'après la formule \ref{eq5.1.2},
si $v=\varpi^{-n} 
\lambda\in\varpi^{-n}A\backslash\varpi^{-n+1}A $, on a
\begin{equation*}
  \chi_{v}(g)=\psi(\frac{1}{2}\varpi^{\lambda_{\psi}-n}Q_{a}(\lambda))
\end{equation*}
avec $a=\left(
\begin{smallmatrix}
\transposee{y}&-z\\
u&x
\end{smallmatrix}\right)
$. Il suit des relations \ref{eq5.1.1} que $a$ est telle que
$\transposee{a}-a\in\varpi^{n}\mathcal{M}_{2r}(\O)$. 

D'autre part, si
$C$ est comme dans l'assertion (ii) du lemme, on a $\varpi^{n}C=\varpi
B^{*}\backslash\varpi A$, $A\backslash B$ ou $B\backslash\varpi
B^{*}$. Si $\lambda\in k^{r}$ est un vecteur
colonne, on écrit
$
\lambda=\left(
\begin{smallmatrix}
\lambda_{1}\\
\lambda_{2}
\end{smallmatrix}\right)
$
où $\lambda_{1}\in k^{l}$ et $\lambda_{2}\in k^{r-l}$. Avec ces
notations, on a
$
A\backslash B =\{v=\left(
\begin{smallmatrix}
\lambda\\
\mu
\end{smallmatrix}\right)\in A\,\vert\,\lambda,\mu\in k^{r},
\Vert\lambda_{1}\Vert_{p}=1\}$, $ 
B\backslash\varpi
B^{*} =\{v=\left(
\begin{smallmatrix}
\lambda\\
\mu
\end{smallmatrix}\right)\in A\,\vert\, \lambda,\mu\in k^{r}, 
\Vert\left(
\begin{smallmatrix}
\lambda_{2}\\
\mu_{2}
\end{smallmatrix}\right)
\Vert_{p}=1, \Vert\lambda_{1}\Vert_{p}<1\}$ et $ 
B^{*}\backslash A = \varpi^{-1}\{v=\left(
\begin{smallmatrix}
\lambda\\
\mu
\end{smallmatrix}\right)\in
W\,\vert\, \lambda,\mu\in k^{r}, \Vert v\Vert_{p}=1,\Vert\left(
\begin{smallmatrix}
\lambda_{2}\\
\mu_{2}
\end{smallmatrix}\right)
\Vert_{p}<1, \Vert\lambda_{1}\Vert_{p}<1\} 
$

Alors, si $C$ est comme dans l'assertion (ii) du lemme, quitte à
effectuer une permutation sur les coordonnées dans $k^{2r}$ et pour un
bon choix de l'entier $m\leq 2r$, la condition $g\in\cap_{v\in
  C}\ker\chi_{v}$ implique la condition \ref{eq5.1.4}. Il résulte
alors de notre assertion, que $a\in\varpi^{n-2}\mathcal{M}_{2r}(\O)$
et $a_{11}$ est à coefficients dans $\varpi^{n}\O$. On en déduit que
$\cap_{v\in C}\ker\chi_{v}\subsetneq K_{2n-2}$.
\end{dem}

Soit donc $n\in\Bbb{N}$, $n\geq1$ et
$v\in\varpi^{-n}A\backslash\varpi^{-n+1}A$. Il suit du lemme
\ref{lem5.1.1} (i) que $\tilde{P}_{B,v}$ est le stabilisateur dans $P_{B}$
du caractère $\chi_{v}$ du sous-groupe distingué $K_{n}$. La méthode
des petits groupes de Mackey montre alors que la représentation
$\Ind_{\tilde{P}_{B,v}}^{P_{B}}\chi_{B,v}^{\pm}$ est irréductible. On en
désigne par $\mathcal{E}_{v,\pm}$ l'espace. 

On rappelle le caractère $\zeta$ de $P_{B}$ défini par les formules
\ref{eq4.5.5} et \ref{eq4.5.6}. Pour $v\in W$, on désigne par
$\zeta\chi^{\pm}_{B,v}$ le caractère de $\tilde{P}_{B,v}$ produit de
$\chi_{B,v}^{\pm}$ par $\zeta_{\vert \tilde{P}_{B,v}}$.
\begin{lem}\label{lem5.1.2}
(i) L'espace de Hilbert $\mathcal{H}_{\psi}^{A}$ est somme
  hilbertienne des sous-espaces $\mathcal{E}^{B}_{n,j}$,
  $n\in\Bbb{N}$, $j=0,1,2$. Le sous-espace $\mathcal{E}^{B}_{0,0}$ est
  non nul. Les sous-espaces
  $\mathcal{E}^{B,\pm}_{n,0}$, $n\in\Bbb{N}\backslash\{0\}$, et
  $\mathcal{E}^{B,\pm}_{n,2}$, $n\in\Bbb{N}$, sont non nuls si et
  seulement si $l>0$. Les espaces $\mathcal{E}^{B,\pm}_{n,1}$,
  $n\in\Bbb{N}$ sont non nuls si et seulement si $l<r$.

(ii) Lorsqu'ils sont non nuls, les sous-espaces
  $\mathcal{E}^{B}_{0,0}$, $\mathcal{E}^{B,\pm}_{n,0}$,
  $n\in\Bbb{N}\backslash\{0\}$, et $\mathcal{E}^{B,\pm}_{n,j}$,
  $n\in\Bbb{N}$, $j=1,2$ sont invariants et irréductibles sous la
  restriction de la représentation $M_{A}$ de $K$ à $P_{B}$.  On note
  $M^{B}_{0,0}$, $M^{B,\pm}_{n,0}$, $n\in\Bbb{N}$, $n\geq1$ et
  $M^{B,\pm}_{n,j}$, $n\in\Bbb{N}$, $j=1,2$ les représentations de
  $P_{B}$ ainsi obtenues. Elles sont deux à deux non équivalentes et
  monomiales. De plus, $M^{B}_{0,0}$ est la représentation triviale et
  l'on a~:
\begin{eqnarray*}
  M^{B,\pm}_{n,0}&=&
  \Ind_{\tilde{P}_{B,v}}^{P_{B}}\chi_{B,v}^{\pm}\mbox{,
  }v\in\varpi^{-n}A\backslash\varpi^{-n}B\mbox{, }n\in\Bbb{N}\mbox{,
  }n\geq1\\
M^{B,\pm}_{n,2}&=&
  \Ind_{\tilde{P}_{B,v}}^{P_{B}}\chi_{B,v}^{\pm}\mbox{,
  }v\in\varpi^{-n}B^{*}\backslash\varpi^{-n}A\mbox{, }n\in\Bbb{N}\\
M^{B,\pm}_{n,1}&=&\Ind_{\tilde{P}_{B,v}}^{P_{B}}\chi_{B,v}^{\pm}
\mbox{, }v\in\varpi^{-(n+1)}B\backslash\varpi^{-n}B^{*}\mbox{,
  }n\in\Bbb{N}.
\end{eqnarray*}

(iii) Lorsqu'ils sont non nuls, les sous-espaces
$\mathcal{E}^{B}_{0,0}$, $\mathcal{E}^{B,\pm}_{n,0}$,
$n\in\Bbb{N}\backslash\{0\}$, et $\mathcal{E}^{B,\pm}_{n,j}$,
$n\in\Bbb{N}$, $j=1,2$ sont invariants et irréductibles sous la
restriction de la représentation métaplectique $S_{\psi}^{A}$ à
$P_{B}$. On note $S^{B}_{0,0}$, $S^{B,\pm}_{n,0}$, $n\in\Bbb{N}$,
$n\geq1$ et $S^{B,\pm}_{n,j}$, $n\in\Bbb{N}$, $j=1,2$ les
représentations de $P_{B}$ ainsi obtenues. Elles sont deux à deux non
équivalentes et l'on a~:
\begin{eqnarray*}
S^{B}_{0,0}&=&\zeta\\
  S^{B,\pm}_{n,0}
&=& \Ind_{\tilde{P}_{B,v}}^{P_{B}}\zeta\chi_{B,v}^{\pm}\mbox{,
  }v\in\varpi^{-n}A\backslash\varpi^{-n}B\mbox{, }n\in\Bbb{N}\mbox{,
  }n\geq1\\
S^{B,\pm}_{n,2}
&=& \Ind_{\tilde{P}_{B,v}}^{P_{B}}\zeta\chi_{B,v}^{\pm}\mbox{,
  }v\in\varpi^{-n}B^{*}\backslash\varpi^{-n}A\mbox{, }n\in\Bbb{N}\\
S^{B,\pm}_{n,1}
&=& \Ind_{\tilde{P}_{B,v}}^{P_{B}}\zeta\chi_{B,v}^{\pm}
\mbox{, }v\in\varpi^{-(n+1)}B\backslash\varpi^{-n}B^{*}\mbox{,
  }n\in\Bbb{N}.
\end{eqnarray*}
\end{lem}
\begin{dem}
(i) Compte tenu des inclusions $\varpi B^{*}\subset B\subset A\subset
  B^{*}$, il est clair que $\mathcal{H}_{\psi}^{A}$ est somme
  hilbertienne des sous-espaces $\mathcal{E}^{B}_{n,j}$,
  $n\in\Bbb{N}$, $j=0,1,2$. De plus, ces inclusions sont strictes si
  $0<l<r$ tandis que $\varpi B^{*}\subsetneq B= A= B^{*}$, si $l=0$, et
  $\varpi B^{*}= B\subsetneq A\subsetneq B^{*}$, si $l=r$.

(ii) Le fait que les sous-espaces $\mathcal{E}^{B}_{0,0}$,
  $\mathcal{E}^{B,\pm}_{n,0}$, $n\in\Bbb{N}\backslash\{0\}$, et
  $\mathcal{E}^{B,\pm}_{n,j}$, $n\in\Bbb{N}$, $j=1,2$, sont invariants
  résulte simplement de l'invariance des réseaux $A$, $B$ et $B^{*}$
  sous l'action de $P_{B}$.

Il est clair que l'action de $P_{B}$ dans $\mathcal{E}^{B}_{0,0}$ est
triviale. D'autre part, d'après le lemme \ref{lem2.7.1}, les orbites
de $P_{B}$ dans $B^{*}\backslash\varpi B^{*}$ sont $B^{*}\backslash
A$, $A\backslash B$ et $B\backslash\varpi B^{*}$. On déduit de ceci
que, pour $v\in\varpi^{-n}A\backslash\varpi^{-n}B$ et
$n\in\Bbb{N}\backslash\{0\}$,
(resp. $v\in\varpi^{-n}B^{*}\backslash\varpi^{-n}A$ et $n\in\Bbb{N}$,
$v\in\varpi^{-(n+1)}B\backslash\varpi^{-n}B^{*}$ et $n\in\Bbb{N}$),
l'application $\varphi\mapsto\tilde{\varphi}$ définie par
$\tilde{\varphi}(p)=\varphi(p^{-1}v)$ induit une bijection de l'espace
$\mathcal{E}^{B,\pm}_{n,0}$ (resp. $\mathcal{E}^{B,\pm}_{n,2}$,
$\mathcal{E}^{B,\pm}_{n,1}$) sur l'espace $\mathcal{E}_{v,\pm}$ qui
entrelace la restriction de $M_{A}$ à $P_{B}$ et la représentation
$\Ind_{\tilde{P}_{B,v}}^{P_{B}}\chi_{B,v}^{\pm}$.

Soit $n\geq1$ et $v\in\varpi^{-n}A\backslash\varpi^{-n+1}A$. Il suit
de l'assertion (ii) du lemme \ref{lem5.1.1} que le noyau de la
représentation $\Ind_{\tilde{P}_{B,v}}^{P_{B}}\chi_{B,v}^{\pm}$
contient $K_{2n}$ et est strictement contenu dans $K_{2n-2}$. Par
suite, si $n\neq m$ sont deux entiers non nuls et si
$v\in\varpi^{-n}A\backslash\varpi^{-n+1}A$,
$w\in\varpi^{-m}A\backslash\varpi^{-m+1}A$, les représentations
$\Ind_{\tilde{P}_{B,v}}^{P_{B}}\chi_{B,v}^{\pm}$ et
$\Ind_{\tilde{P}_{B,w}}^{P_{B}}\chi_{B,w}^{\pm}$ ne sont pas
équivalentes et distinctes de la représentation triviale. D'autre part,
si $n$ est un entier non nul,
$v,w\in\varpi^{-n}A\backslash\varpi^{-n+1}A$ et
$\epsilon,\epsilon'\in\{\pm\}$, il suit de la méthode des petits groupes
de Mackey appliquée au sous-groupe invariant $K_{n}$ que les
représentations $\Ind_{\tilde{P}_{B,v}}^{P_{B}}\chi_{B,v}^{\epsilon}$
et $\Ind_{\tilde{P}_{B,w}}^{P_{B}}\chi_{B,w}^{\epsilon'}$ ne peuvent
être équivalentes que si $v$ et $w$ sont dans la même $P_{B}$-orbite
et $\epsilon'=\epsilon$. Il en résulte que les représentations
$M^{B}_{0,0}$, $M^{B,\pm}_{n,0}$, $n\in\Bbb{N}$, $n\geq1$ et
$M^{B,\pm}_{n,j}$, $n\in\Bbb{N}$, $j=1,2$, sont deux à deux non
équivalentes.

(iii) Cette assertion est conséquence immédiate de l'assertion (ii) et
du fait que, d'après le corollaire \ref{co4.5.1}, la restriction de la
représentation métaplectique $S_{\psi}^{A}$ à $P_{B}$ est égale à
$\zeta\otimes M_{A\vert P_{B}}$.
\end{dem}

\subsection{}\label{5.2}
On garde les notations du paragraphe précédent. Pour $n\in\Bbb{N}$, on
pose
\begin{equation*}
\mathcal{E}^{B}_{n}
=\mathcal{E}^{B}_{n,0}\oplus\mathcal{E}^{B}_{n,2}.
\end{equation*}
On a
\begin{equation*}
\begin{split}
\mathcal{E}^{B}_{0}&=\{f\in\mathcal{H}_{\psi}^{A}\vert\Supp
f\subset B^{*}\},\\ \mathcal{E}^{B}_{n}&=
\{f\in\mathcal{H}_{\psi}^{A}\vert\Supp
f\subset\varpi^{-n}(B^{*}\backslash B)\}\mbox{,
}n\in\Bbb{N}\backslash\{0\}.
\end{split}
\end{equation*}

\begin{theo}\label{theo5.2.1}
(i) Le sous-espace $\mathcal{E}^{B}_{0}$ est non nul. Les sous-espaces
  $\mathcal{E}^{B,-}_{0}$ et $\mathcal{E}^{B,\pm}_{n}$,
  $n\in\Bbb{N}\backslash\{0\}$, sont non nuls si et seulement si
  $l>0$. Les sous-espaces $\mathcal{E}^{B,\pm}_{n,1}$, $n\in\Bbb{N}$,
  sont non nuls si et seulement si $l<r$.

(ii) Lorsqu'ils sont non nuls, les sous-espaces
  $\mathcal{E}^{B,\pm}_{n}$ et $\mathcal{E}^{B,\pm}_{n,1}$,
  $n\in\Bbb{N}$ sont invariants et irréductibles sous la restriction
  de la représentation métaplectique $S_{\psi}^{A}$ à $K_{B}$. Les
  représentations de $K_{B}$ ainsi obtenues sont deux à deux non
  équivalentes. La restriction de la représentation métaplectique à
  $K_{B}$ est sans multiplicité et somme directe de ces
  représentations.

(iii) On a~:
\begin{eqnarray*}
\dim\mathcal{E}^{B,+}_{0}&=&1+\frac{1}{2}(q^{l}-1)\\ 
\dim\mathcal{E}^{B,-}_{0}&=&\frac{1}{2}(q^{l}-1)\\ 
\dim\mathcal{E}^{B,\pm}_{n}&=&\frac{1}{2}q^{2rn-l}(q^{2l}-1)\mbox{,
}n\geq1\\
\dim\mathcal{E}^{B,\pm}_{n,1}&=&\frac{1}{2}q^{2rn+l}(q^{2(r-l)}-1)\mbox{,
}n\geq0.
\end{eqnarray*}
\end{theo}
\begin{dem}
(i) C'est une conséquence immédiate de l'assertion (i) du lemme
  \ref{lem5.1.2}.

(ii) Compte tenu du fait que $K_{B}$ est engendré par $P_{B}$ et
  $\varsigma_{B}$, de la formule \ref{eq4.5.8} du corollaire
  \ref{co4.5.1} et du lemme \ref{lem5.1.2}, il suffit de montrer que
  pour tout entier naturel $n$, les espaces $\mathcal{E}^{B}_{n}$ et
  $\mathcal{E}^{B}_{n,1}$ sont invariants sous l'opérateur
  $M_{A}(\varsigma_{B})$ et qu'aucun des espaces $\mathcal{E}^{B,\pm}_{n,0}$ ne
  l'est.

Commençons par remarquer qu'il suit du lemme \ref{lem4.4.1} que si
$v\in W$, $\varsigma_{B}\delta_{v}$ est combinaison linéaire de
$\delta_{w}$ avec $w\in\varsigma_{B}(v+A)$. Soit $n$ un entier
naturel. 
L'espace $\mathcal{E}^{B}_{n,1}$ est engendré par les
vecteurs $\delta_{v}$, $v\in\varpi^{-n}(\varpi^{-1}B\backslash
B^{*})$. L'ensemble $\varpi^{-1}B\backslash B^{*}$ est à la fois
$K_{B}$-invariant et réunion de classes modulo $A$. Notre remarque
montre alors que $\mathcal{E}^{B}_{n,1}$ est invariant sous l'action
de $M_{A}(\varsigma_{B})$, comme voulu.

L'espace $\mathcal{E}^{B}_{n,0}$ (resp. $\mathcal{E}^{B}_{n,2}$) est
engendré par les vecteurs $\delta_{v}$, $v\in\varpi^{-n}(A\backslash
B)$ (resp. $v\in\varpi^{-n}(B^{*}\backslash A)$). Or, il est clair que
$K_{B}$ laisse stable $B^{*}\backslash B$ et on vérifie
facilement que $\varsigma_{B}(A\backslash B)\subset B^{*}\backslash
A$. Ceci montre que l'espace $\mathcal{E}^{B}_{n}$ est invariant sous
l'action de $M_{A}(\varsigma_{B})$ et qu'aucun des espaces
$\mathcal{E}^{B,\pm}_{n,0}$ ne l'est. Ceci achève la
démonstration de l'assertion (ii).

(iii) La multiplication des scalaires induit une action du groupe
$\{\pm1\}$ dans $W/A$ et on considère l'espace quotient de cette
action $\{\pm1\}\backslash W/A$.

On a
$\mathcal{E}^{B,+}_{0}=\Bbb{C}\delta_{0}\oplus\mathcal{E}^{B,+}_{0,2}$
et $\mathcal{E}^{B,-}_{0}=\mathcal{E}^{B,-}_{0,2}$, tandis qu'une
base de $\mathcal{E}^{B,\pm}_{0,2}$ est constituée des
$\delta_{v}\pm\delta_{-v}$, pour $v$ parcourant un système de
représentants des classes de $\{\pm1\}\backslash W/A$ contenues dans
$B^{*}\backslash A$. On a donc
$\dim\mathcal{E}^{B,\pm}_{0,2}=\frac{1}{2}([B^{*}/A]-1)$. Mais la
suite exacte $0\longrightarrow A/B\longrightarrow
B^{*}/B\longrightarrow B^{*}/A\longrightarrow 0$ montre que
$q^{2l}=[B^{*}/B]=[A/B][B^{*}/A]$. Or, $A/B$ étant un lagrangien du
$\Bbb{F}_{q}$-espace symplectique ${\mathrm{b}
}^{*}$ de dimension $2l$, on a
$[A/B]=q^{l}$ et donc $[B^{*}/A]=q^{l}$. D'où la formule pour les
dimensions des espaces $\mathcal{E}^{B,\pm}_{0}$.

Soit $n\in\mathbb{N}\backslash\{0\}$. Comme 
$\varpi^{-n}(B^{*}\backslash B)$
(resp. $\varpi^{-n}(\varpi^{-1}B\backslash B^{*})$) ne rencontre pas
$A$, on voit que $\dim\mathcal{E}^{B,\pm}_{n}=
\frac{1}{2}([\varpi^{-n}B^{*}/A]-[\varpi^{-n}B/A])$
(resp. $\dim\mathcal{E}^{B,\pm}_{n,1}=
\frac{1}{2}([\varpi^{-(n+1)}B/A]-[\varpi^{-n}B^{*}/A])$). Or, on a
$[\varpi^{-n}B^{*}/A]=q^{2rn+l}$ et $[\varpi^{-n}B/A]=q^{2rn-l}$. D'où
le résultat cherché.
\end{dem}

\subsection{}\label{5.3}
On reprend les notations du paragraphe précédent et on désigne par
$S^{B,\pm}_{n}$ (resp. $S^{B,\pm}_{n,1}$) la représentation de $K_{B}$
induite par la représentation $S_{\psi}^{A}$ dans le sous-espace
$\mathcal{E}^{B,\pm}_{n}$ (resp. $\mathcal{E}^{B,\pm}_{n,1})$ lorsque
celui-ci est non nul. Nous allons utiliser les résultats de
\cite{cliff-mcneilly-szechtman-2003} rappelés dans le paragraphe
\ref{3.4} pour donner une description des représentations
$S^{B,\pm}_{n}$ et $S^{B,\pm}_{n,1}$ comme représentations induites.

\'Etant donné $x\in W$, on désigne par $\hat{x}$ (resp. $\dot{x}$) la
classe de $x$ dans le $\O$-module quotient $W/B^{*}$ (resp. $W/B$) et
on rappelle que $x$ s'identifie à l'élément $(x,0)$ de $H(W)$.
\begin{lem}\label{lem5.3.1}
Soit $x\in W$ et soit $n$ le plus petit entier naturel tel que
$x\in\varpi^{-(n+1)}B$.

(i) Pour tout $g\in K_{B}(\hat{x})$, le commutateur
$g^{-1}xgx^{-1}$ est contenu dans le sous-groupe
$B^{*}\times\varpi^{\lambda_{\psi}-n-1}\O$ de $H(W)$.

(ii) Si $x\notin B^{*}$, on a
$K_{B}(\hat{x})=K_{B}(x)K'_{B}(\hat{x})$.

(iii) Si $x\in\varpi^{-(n+1)}B\backslash\varpi^{-n}B^{*}$, on a
\begin{equation*}
p_{Sp({\mathrm{b} }^{*})}(K_{B}(x))=Sp({\mathrm{b} }^{*}).
\end{equation*}

(iv) Si  $x\in\varpi^{-n}B^{*}\backslash\varpi^{-n}B$, on a
\begin{equation*}
p_{Sp({\mathrm{b} }^{*})}(K_{B}(x))=Sp({\mathrm{b} }^{*})(p_{{\mathrm{b}
}^{*}}(\varpi^{n}x)).
\end{equation*}
\end{lem}
\begin{dem}
(i) Si $g\in K_{B}(\hat{x})$, on a 
\begin{equation}\label{eq5.3.1}
g^{-1}xgx^{-1}=(g^{-1}x-x,\frac{1}{2}\beta(x,g^{-1}x))
\end{equation}
avec $g^{-1}x-x\in B^{*}$ et
\begin{equation*}
\frac{1}{2}\beta(x,g^{-1}x)=\frac{1}{2}\beta(x-g^{-1}x,g^{-1}x)
\in\beta(B^{*},\varpi^{-(n+1)}B)
\subset\varpi^{\lambda_{\psi}-n-1}\O.
\end{equation*}

(ii) L'inclusion $K_{B}(x)K'_{B}(\hat{x})\subset K_{B}(\hat{x})$ est
claire. Soit donc $g\in K_{B}(\hat{x})$. Par définition, on a $gx\in
x+B^{*}$. D'autre part, $\varpi^{n+1}x\in B$ et $g\varpi^{n+1}x\in
\varpi^{n+1}x+\varpi^{n+1}B^{*}$. Si $\varpi^{n+1}x\notin\varpi
B^{*}$, l'assertion (i) du lemme \ref{lem2.7.1} montre qu'il existe
$h\in K'_{B}$ tel que $hg\varpi^{n+1}x=\varpi^{n+1}x$. Si
$\varpi^{n+1}x\in\varpi B^{*}$, on a $n\geq1$, $\varpi^{n}x\in
B^{*}\backslash B$ et
$g\varpi^{n}x\in\varpi^{n}x+\varpi^{n}B^{*}\subset\varpi^{n}x+B$, de
sorte que le même argument montre qu'il existe $h\in K'_{B}$ tel que
$hg\varpi^{n}x=\varpi^{n}x$. Dans tous les cas, on a trouvé $h\in
K'_{B}$ tel que $hgx=x$ et il est alors évident que $h\in
K'_{B}(\hat{x})$.

(iii) Soit $g\in Sp({\mathrm{b}
}^{*})$. Désignons par $K_{B}''$ le noyau du
morphisme $p_{Sp(\mathrm{b})}$ de $K_{B}$ sur $Sp(\mathrm{b})$. Il suit du
lemme \ref{lem2.4.1} que
$p_{Sp({\mathrm{b}
}^{*})}(K_{B}'')=Sp({\mathrm{b}
}^{*})$. Il existe donc
$\tilde{g}\in K_{B}''$ tel que
$p_{Sp({\mathrm{b}
}^{*})}(\tilde{g})=g$. Mais, par définition de $K_{B}''$,
on a $\tilde{g}(\varpi^{n+1}x)\in\varpi^{n+1}x+\varpi B^{*}$. D'après
l'assertion (i) du lemme \ref{lem2.7.1}, il existe $h\in K'_{B}$ tel
que $h\tilde{g}\varpi^{n+1}x=\varpi^{n+1}x$. On a alors $h\tilde{g}\in
K_{B}(x)$ et
$p_{Sp({\mathrm{b}
}^{*})}(h\tilde{g})=p_{Sp({\mathrm{b}
}^{*})}(\tilde{g})=g$. D'où
l'assertion (iii).

(iv) L'inclusion $p_{Sp({\mathrm{b}
}^{*})}(K_{B}(x))\subset
Sp({\mathrm{b}
}^{*})(p_{{\mathrm{b}
}^{*}}(\varpi^{n}x))$ est claire. Soit donc
$g$ et $\tilde{g}$ des éléments respectifs de $Sp({\mathrm{b}
}^{*})(p_{{\mathrm{b}
}^{*}}(\varpi^{n}x))$ et 
$K_{B}$ tels que $p_{Sp({\mathrm{b}
}^{*})}(\tilde{g})=g$. On a alors
$\tilde{g}\varpi^{n}x\in\varpi^{n}x+B$. Utilisant l'assertion (i) du
lemme \ref{lem2.7.1}, on conclut comme pour l'assertion (iii).
\end{dem}
Soit $x\in W$. Il suit du lemme précédent que
$x^{-1}K_{B}(\hat{x})x$ est un sous-groupe de
$K_{B}\ltimes\overline{H}(B^{*})$. Rappelons-nous la représentation
$\widetilde{S}_{\overline{\psi}}\widetilde{\rho}_{\overline{\psi}}$ de
$K_{B}\ltimes\overline{H}(B^{*})$ définie au paragraphe \ref{4.2} par
la formule \ref{eq4.2.1}. On définit alors la représentation $\sigma_{\dot{x}}$
du groupe $K_{B}(\hat{x})$ en posant
\begin{equation}\label{eq5.3.2}
\sigma_{\dot{x}}(g)=
\widetilde{S}_{\overline{\psi}}\widetilde{\rho}_{\overline{\psi}}(xgx^{-1})
\mbox{, }g\in K_{B}(\hat{x})
\end{equation}
(on vérifie que la représentation $\sigma_{\dot{x}}$ ne dépend que de la
classe $\dot{x}\in W/B$).

Supposons que $x\notin B^{*}$. Il suit de la formule \ref{eq5.3.1} que
l'on a 
\begin{equation}\label{eq5.3.3}
\sigma_{\dot{x}}(g)=
\psi(\frac{1}{2}\beta(x,g^{-1}x))\widetilde{S}_{\overline{\psi}}(g)
\widetilde{\rho}_{\overline{\psi}}(g^{-1}x-x,0)\mbox{, }g\in K_{B}(\hat{x}).
\end{equation}
De plus, il résulte de l'assertion (ii) du lemme \ref{lem5.3.1} que la
représentation $\sigma_{\dot{x}}$ satisfait les relations suivantes
qui la déterminent entièrement~:
\begin{equation*}
\begin{split}
  \sigma_{\dot{x}}(g)&= \widetilde{S}_{\overline{\psi}}(g)\mbox{, }g\in K_{B}(x),\\
  \sigma_{\dot{x}}(h)& =\psi(\frac{1}{2}\beta(x,h^{-1}x))
\widetilde{\rho}_{\overline{\psi}}(h^{-1}x-x,0)\\
&=\psi(\frac{1}{2}\beta(hx,x))
\widetilde{\rho}_{\overline{\psi}}(x-hx,0)\mbox{, }h\in
K'_{B}(\hat{x}).
\end{split}
\end{equation*}
Dans la suite, on désigne par $\widetilde{K_{B}(\hat{x})}$ le
sous-groupe $\{\pm Id\}K_{B}(\hat{x})$ de $K_{B}$.

Lorsque $x\in B^{*}$,
$\widetilde{K_{B}(\hat{x})}=K_{B}(\hat{x})=K_{B}$ et la représentation
$\sigma_{\dot{x}}$ de $K_{B}$ est équivalente à
$\widetilde{S}_{\overline{\psi}}=S_{0}$ et est réalisée dans le même
espace $\mathcal{H}_{\overline{\psi}}$. On a $-Id\in K_{B}$ et on note
$\mathcal{H}^{\pm}_{\overline{\psi}}$ le sous-espace propre de $-Id$
pour la valeur propre $\pm1$, lequel est invariant par
$\sigma_{\dot{x}}$. On désigne par $\sigma^{\pm}_{\dot{x}}$ la
représentation de $K_{B}$ dans ce sous-espace qui en résulte, lorsque
celui-ci est non nul. On a
$\mathcal{H}_{\overline{\psi}}=\mathcal{H}^{+}_{\overline{\psi}}
\oplus\mathcal{H}^{-}_{\overline{\psi}}$ et donc
$\sigma_{\dot{x}}=\sigma^{+}_{\dot{x}}\oplus \sigma^{-}_{\dot{x}}$, lorsque
$\mathcal{H}^{-}_{\overline{\psi}}$ est non nul~; ceci se produit
exactement lorsque $l\neq 0$. Dans le cas contraire, on a
$\sigma_{\dot{x}}=\sigma^{+}_{\dot{x}}$.

Supposons que $x\notin B^{*}$. Dans ce cas, $-Id$ n'appartient pas à
$K_{B}(\hat{x})$. On étend alors la représentation $\sigma_{\dot{x}}$
en une représentation $\sigma^{\pm}_{\dot{x}}$ du groupe
$\widetilde{K_{B}(\hat{x})}$ en décidant que $\sigma^{\pm}_{\dot{x}}(-Id)=\pm
Id$.
\begin{theo}\label{theo5.3.1}
(i) Si $l=0$, $S^{B,+}_{0}$ est la représentation triviale.

(ii) Si $l>0$, pour tout $n\in\Bbb{N}$, on a
\begin{equation*}
  S^{B,\pm}_{n}=\Ind_{\widetilde{K_{B}(\hat{x})}}^{K_{B}}\sigma^{\pm}_{\dot{x}}\mbox{,
  }x\in(\varpi^{-n}B^{*}\backslash\varpi^{-n}B)+B^{*}.
\end{equation*}

(iii) Si $l<r$, pour tout $n\in\Bbb{N}$, on a
\begin{equation*}
  S^{B,\pm}_{n,1}=\Ind_{\widetilde{K_{B}(\hat{x})}}^{K_{B}}\sigma^{\pm}_{\dot{x}}\mbox{,
  }x\in(\varpi^{-(n+1)}B\backslash\varpi^{-n}B^{*})+B^{*}.
\end{equation*}
\end{theo}
\begin{dem}
On utilise la réalisation $S_{\psi}^{B}$ de la représentation de Weil
dans l'espace $\mathcal{H}_{\psi}^{B}$ donnée au paragraphe \ref{4.2}
et l'opérateur d'entrelacement $\mathcal{I}_{A,B}$ entre cette
dernière et la réalisation $S_{\psi}^{A}$ dans l'espace
$\mathcal{H}_{\psi}^{A}$.

Si $\mathcal{F}$ est un sous-espace de $\mathcal{H}_{\psi}^{B}$, on
désigne par $\mathcal{F}^{+}$ (resp. $\mathcal{F}^{-}$) le sous-espace
de $\mathcal{F}$ constitué des fonctions paires (resp. impaires).

Si $x\in W$, on désigne par $\mathcal{F}_{\hat{x}}$ le sous-espace de
$\mathcal{H}_{\psi}^{B}$ constitué des fonctions dont le support est
contenu dans $K_{B}x$, si $x\notin B^{*}$, et dans $B^{*}$, sinon (la
définition et la notation sont justifiées parce que, d'après le lemme
\ref{lem2.7.1} (i), $x+B^{*}\subset K_{B}x$, si $x\in W\backslash
B^{*}$). Il est clair que les sous-espaces
$\mathcal{F}_{\hat{x}}^{\pm}$ sont invariants sous l'action de la
représentation $S_{\psi}^{B}$ restreinte à $K_{B}$. D'autre part, il
suit du lemme \ref{lem2.7.1} (ii) que les $K_{B}$-orbites dans
$W\backslash B^{*}$ sont les
$\varpi^{-(n+1)}B^{*}\backslash\varpi^{-(n+1)}B$ et
$\varpi^{-(n+1)}B\backslash\varpi^{-n}B^{*}$ $n\in\Bbb{N}$. Comme
l'opérateur d'entrelacement $\mathcal{I}_{A,B}$ conserve les supports
et la parité, on déduit de ceci que l'on a
\begin{equation*}
\begin{split}
  \mathcal{I}_{A,B}(\mathcal{F}_{\hat{x}}^{\pm}) &=
\mathcal{E}^{B\pm}_{0}\mbox{, }
x\in B^{*},\\
\mathcal{I}_{A,B}(\mathcal{F}_{\hat{x}}^{\pm}) & =
\mathcal{E}^{B\pm}_{n+1}\mbox{, }
x\in\varpi^{-(n+1)}B^{*}\backslash\varpi^{-(n+1)}B\mbox{, } n\in\Bbb{N},\\
\mathcal{I}_{A,B}(\mathcal{F}_{\hat{x}}^{\pm}) & =
\mathcal{E}^{B\pm}_{n,1}\mbox{, }
x\in\varpi^{-(n+1)}B\backslash\varpi^{-n}B^{*}\mbox{, } n\in\Bbb{N}.
\end{split}
\end{equation*}
Soit $x\in W$. Désignons par $\mathcal{G}_{\dot{x}}^{\pm}$ l'espace de la
représentation
$\Ind_{\widetilde{K_{B}(\hat{x})}}^{K_{B}}\sigma^{\pm}_{\dot{x}}$. Alors,
l'application $\varphi\longmapsto\tilde{\varphi}$ définie par
$\tilde{\varphi}(g)=
\widetilde{S}_{\overline{\psi}}(g)\varphi(g^{-1}x)$, $g\in K_{B}$,
induit un isomorphisme de $K_{B}$-modules de $\mathcal{F}_{\hat{x}}^{\pm}$
sur $\mathcal{G}_{\dot{x}}^{\pm}$. D'où le théorème.
\end{dem}

\subsection{}\label{5.4}
Dans ce paragraphe nous allons mettre en relation les résultats du
précédent avec ceux du paragraphe \ref{3.5}, dont nous reprenons les
notations. 

Soit $n$ un entier naturel. Le sous-ensemble
$H(\varpi^{-(n+1)}B)=\varpi^{-(n+1)}B\times\varpi^{\lambda_{\psi}-1-2(n+1)}\O$
est un sous-groupe $K_{B}$-invariant de $H(W)$.

On désigne par $\mathcal{H}^{B}_{n}$ le sous-espace de
$\mathcal{H}_{\psi}^{B}$ constitué des fonctions dont le support est
contenu dans $\varpi^{-(n+1)}B$. En fait, $\mathcal{H}^{B}_{n}$ est
l'espace des fonctions $\varphi$ définies sur $\varpi^{-(n+1)}B$ à
valeurs dans l'espace $\mathcal{H}_{\overline{\psi}}$ de la
représentation de Schrödinger de type $\overline{\psi}$ du groupe de
Heisenberg $H(\mathrm{b}^{*})$, vérifiant la relation
\ref{eq4.1.1}. Il est clairement invariant sous la restriction de
$S_{\psi}^{B}$ (resp. $R^{B}_{\psi}$) à $K_{B}$ (resp. $K_{B}\ltimes
H(\varpi^{-(n+1)}B)$)~: on note $S^{B,n}_{\psi}$
(resp. $R^{B,n}_{\psi}$) la représentation de $K_{B}$
(resp. $K_{B}\ltimes H(\varpi^{-(n+1)}B)$) induite par cette dernière
dans $\mathcal{H}^{B}_{n}$.

Soit $\psi_{2n+1}$ le caractère primitif de $O_{2n+1}$ défini par 
\begin{equation*}
  \psi_{2n+1}(p_{O_{2n+1}}(t))=\psi(\varpi^{\lambda_{\psi}-2(n+1)}t)\mbox{, }t\in\O.
\end{equation*}
On considère le $O_{2n+1}$-module symplectique
${\mathrm{b}}_{2n+1}=B/\varpi^{2(n+1)}B^{*}$ et le sous-module isotrope
$Sp({\mathrm{b}}_{2n+1})$-invariant maximal
$\mathrm{u}=\varpi^{n+1}B/\varpi^{2(n+1)}B^{*}$ (voir le lemme
\ref{lem3.5.1}). On rappelle la représentation
$\rho_{\mathrm{u}^{\perp},\psi_{2n+1}}$ de $H(\mathrm{u}^{\perp})$, inflation
de la représentation de Schrödinger
$\rho_{\overline{\mathrm{u}^{\perp}},\psi_{2n+1}}$ de type $\psi_{2n+1}$
de $H(\overline{\mathrm{u}^{\perp}})$.

Soit $\mu:\Bbb{F}_{q}\longrightarrow O_{2n+1}$ l'application définie
par~:
\begin{equation*}
\mu(p_{\Bbb{F}_{q}}(t))=p_{O_{2n+1}}(\varpi^{2n+1}t)\mbox{, }t\in\O.
\end{equation*}
Il est immédiat que $\mu$ induit un isomorphisme de $\Bbb{F}_{q}$-modules de
$\Bbb{F}_{q}$ sur l'idéal minimal $\varpi^{2n+1}\O/\varpi^{2(n+1)}\O$.
De plus, on a
\begin{equation}\label{eq5.4.1}
\psi_{2n+1}\circ\mu=\overline{\psi}.
\end{equation}

L'application
$\mu_{*}:H(\mathrm{b}^{*})\longrightarrow H(\overline{\mathrm{u}^{\perp}})$
définie par~:
\begin{equation}\label{eq5.4.2}
  \mu_{*}(p_{\mathrm{b}^{*}}(x),t)=
(p_{{\mathrm{b}}_{2n+1}}(\varpi^{n+1}x),\mu(t))\mbox{,
  }x\in B^{*}\mbox{, }t\in\Bbb{F}_{q},
\end{equation}
est un morphisme injectif de groupes.
Il suit de la relation \ref{eq5.4.1} que la représentation
$\rho_{\overline{\mathrm{u}^{\perp}},\psi_{2n+1}}\circ\mu_{*}$ est la
représentation de Schrödinger de type $\overline{\psi}$ de
$H(\mathrm{b}^{*})$. On peut donc supposer que
$\mathcal{H}_{\overline{\psi}}$ est également l'espace de la
représentation $\rho_{\overline{\mathrm{u}^{\perp}},\psi_{2n+1}}$ et
écrire alors 
\begin{equation}\label{eq5.4.3}
  \rho_{\overline{\psi}}=\rho_{\overline{\mathrm{u}^{\perp}},\psi_{2n+1}}\circ\mu_{*}.
\end{equation}
Par suite, la représentation de Weil $S_{\overline{\psi}}$ de type
$\overline{\psi}$ de $Sp(\mathrm{b}^{*})$ est également une représentation
de Weil de type $\psi_{2n+1}$ de $Sp(\overline{\mathrm{u}^{\perp}})$. On
choisit alors pour représentation de Weil de type $\psi_{2n+1}$ de
$Sp({\mathrm{b}}_{2n+1})$ relative au morphisme
$r_{\mathrm{u}}:Sp({\mathrm{b}}_{2n+1}) \longrightarrow
Sp(\overline{\mathrm{u}^{\perp}})$ la représentation
$\sigma=S_{\overline{\psi}}\circ r_{\mathrm{u}}$.

On note $\mathcal{H}_{n}$ l'espace de la représentation
\begin{equation}\label{eq5.4.4}
R^{\mathrm{u},\sigma}= \Ind_{Sp({\mathrm{b}}_{2n+1})\ltimes
  H(\mathrm{u}^{\perp})} ^{Sp({\mathrm{b}}_{2n+1})\ltimes
  H({\mathrm{b}}_{2n+1})}\sigma\rho_{\mathrm{u}^{\perp},\psi_{2n+1}}.
\end{equation}
La restriction $S^{\mathrm{u},\sigma}$ de la représentation
$R^{\mathrm{u},\sigma}$ à $Sp({\mathrm{b}}_{2n+1})$ est une
représentation de Weil de type $\psi_{2n+1}$. 

Il suit de la relation \ref{eq3.4.1} que l'espace de la représentation
$R^{\mathrm{u},\sigma}$ s'identifie, via l'application de restriction au
sous-ensemble ${\mathrm{b}}_{2n+1}$ du sous-groupe $H({\mathrm{b}}_{2n+1})$ du
produit semi-direct $Sp({\mathrm{b}}_{2n+1})\ltimes H({\mathrm{b}}_{2n+1})$, à
l'espace des fonctions $\varphi$ définies sur ${\mathrm{b}}_{2n+1}$, à
valeurs dans l'espace $\mathcal{H}_{\overline{\psi}}$ de la
représentation $\rho_{\mathrm{u}^{\perp},\psi_{2n+1}}$, qui vérifient la
relation
\begin{equation}\label{eq5.4.5}
\varphi(x+u)=
\psi_{2n+1}(\frac{1}{2}\beta_{{\mathrm{b}}_{2n+1}}(x,u))\rho_{\mathrm{u}^{\perp},\psi_{2n+1}}(u)
\varphi(x) \mbox{, }x\in {\mathrm{b}}_{2n+1}\mbox{, }u\in \mathrm{u}^{\perp}.
\end{equation}

L'application $p:K_{B}\ltimes H(\varpi^{-(n+1)}B)\longrightarrow
Sp({\mathrm{b}}_{2n+1})\ltimes H({\mathrm{b}}_{2n+1})$
définie par~:
\begin{equation*}
p(g(x,t))=
p_{Sp({\mathrm{b}}_{2n+1})}(g)(p_{{\mathrm{b}}_{2n+1}}(\varpi^{n+1}x),
p_{O_{2n+1}}(\varpi^{2(n+1)-\lambda}t)),
\end{equation*}
$g\in K_{B}$, $(x,t)\in H(\varpi^{-(n+1)}B)$, est un morphisme
surjectif de groupes.

\begin{lem}\label{lem5.4.1}
(i) Si $\varphi$ est un élément de $\mathcal{H}_{n}$, l'application
$\varphi^{B}$ définie sur $\varpi^{-(n+1)}B$ par 
\begin{equation*}
\varphi^{B}(x)=\varphi(p_{{\mathrm{b}}_{2n+1}}(\varpi^{n+1}x))\mbox{,
}x\in\varpi^{-(n+1)}B,
\end{equation*}
est un élément de $\mathcal{H}^{B}_{n}$.

(ii) L'application $\varphi\mapsto\varphi^{B}$ est un isomorphisme de
$\mathcal{H}_{n}$ sur $\mathcal{H}^{B}_{n}$ qui entrelace les
représentations $S^{\mathrm{u},\sigma}\circ p_{Sp({\mathrm{b}}_{2n+1})}$ et
$S^{B,n}_{\psi}$ de $K_{B}$.

(iii) Soit $x\in\varpi^{-(n+1)}B$ et
$\varphi\in\mathcal{H}_{n}$. Alors dire que le support de
$\varphi^{B}$ est contenu dans $K_{B}x$ est équivalent à dire que le
support de $\varphi$ est contenu dans l'orbite de
$p_{{\mathrm{b}}_{2n+1}}(\varpi^{n+1}x)$ modulo $\mathrm{u}^{\perp}$ sous
l'action de $Sp({\mathrm{b}}_{2n+1})$.
\end{lem}
\begin{dem}
(i) Soit $\varphi\in\mathcal{H}_{n}$. Soit $x\in\varpi^{-(n+1)}B$ et
  $b\in B^{*}$. Alors, on a
\begin{equation*}
\begin{split}
\varphi^{B}(x+b) &
=\varphi(p_{{\mathrm{b}}_{2n+1}}(\varpi^{n+1}(x+b)))\\ 
&=\psi_{2n+1}(\frac{1}{2}
\beta_{{\mathrm{b}}_{2n+1}}(p_{{\mathrm{b}}_{2n+1}}(\varpi^{n+1}x),
p_{{\mathrm{b}}_{2n+1}}(\varpi^{n+1}b)))\\ 
&\phantom{=} \quad\rho_{\mathrm{u},\psi_{2n+1}}(p_{{\mathrm{b}}_{2n+1}}(\varpi^{n+1}b))
\varphi^{B}(x)\\
&=\psi_{2n+1}(p_{O_{2n+1}}(\varpi^{2(n+1)-\lambda}\frac{1}{2}\beta(x,b)))
\rho_{\mathrm{u},\psi_{2n+1}}\circ\mu_{*}(p_{\mathrm{b}^{*}}(b))\varphi^{B}(x)\\
&= \psi_{2n+1}\circ\mu(p_{\Bbb{F}_{q}}(\varpi^{1-\lambda}\frac{1}{2}\beta(x,b)))
\rho_{\mathrm{u},\psi_{2n+1}}\circ\mu_{*}(p_{\mathrm{b}^{*}}(b))\varphi^{B}(x)\\
&=\overline{\psi}(p_{\Bbb{F}_{q}}(\varpi^{1-\lambda}\frac{1}{2}\beta(x,b)))
\rho_{\overline{\psi}}(p_{\mathrm{b}^{*}}(b))\varphi^{B}(x)\\
&=\psi(\frac{1}{2}\beta(x,b))\widetilde{\rho}_{\overline{\psi}}(b)\varphi^{B}(x),
\end{split}
\end{equation*}
montrant que $\varphi^{B}$ satisfait la relation \ref{eq4.1.1}.

(ii) Il est clair que l'application $\varphi\mapsto\varphi^{B}$ est
injective. Réciproquement, soit $\varphi\in\mathcal{H}^{B}_{n}$. La relation
\ref{eq4.1.1} satisfaite par $\varphi$, montre qu'il existe une unique
fonction $\varphi^{\mathrm{u}}$ définie sur ${\mathrm{b}}_{2n+1}$ et à valeurs
dans $\mathcal{H}_{\overline{\psi}}$ telle que
\begin{equation*}
\varphi^{\mathrm{u}}(p_{{\mathrm{b}}_{2n+1}}(\varpi^{n+1}x))=\varphi(x)\mbox{,
}x\in\varpi^{-(n+1)}B
\end{equation*}
et que cette fonction satisfait la relation \ref{eq5.4.5}. Ceci montre
que l'application $\varphi\mapsto\varphi^{B}$ est bien une bijection
linéaire de $\mathcal{H}_{n}$ sur $\mathcal{H}^{B}_{n}$. Le fait que
ce soit un opérateur d'entrelacement est facile et est laissé au
lecteur.

(iii) est clair.
\end{dem}

\begin{rems}
(i) Il suit de l'assertion (iii) du lemme précédent que, pour tout
  $x\in\varpi^{-(n+1)}B$, l'application $\varphi\mapsto\varphi^{B}$
  induit un isomorphisme entre les représentations
  $S^{\pm}_{p_{{\mathrm{b}}_{2n+1}}(\varpi^{n+1}x)}\circ
  p_{Sp({\mathrm{b}}_{2n+1})}$ et $S^{\pm}_{\dot{x}}$ de $K_{B}$. En
  particulier, la décomposition en irréductibles de la restriction à
  $K_{B}$ de la représentation de Weil de $Sp(W)$ se ramène à la
  décomposition en irréductibles des représentations de Weil des
  groupes symplectiques $Sp({\mathrm{b}}_{2n+1})$ sur l'anneau local
  fini $O_{2n+1}$, démontrée dans le théorème \ref{theo3.5.1}.

(ii) La décomposition en irréductibles de la restriction à $K_{B}$ de
la représentation de Weil a été obtenue par D. Prasad dans
\cite{prasad-ip-1997} dans le cas où $B=B^{*}$, i.e. $K_{B}$ est le
compact maximal standard. Dans \cite{cliff-mcneilly-szechtman-2003}
G. Cliff, D. McNeilly et F. Szechtman remarquent que leurs résultats
concernant la représentation de Weil de $Sp(W)$ lorsque $W$ est un
module symplectique libre sur un anneau fini local et principal,
permettent d'obtenir la décomposition en irréductibles de la
restriction à $K_{B}$ de la représentation de Weil également lorsque
$B=\varpi B^{*}$, i.e. $K_{B}$ n'est pas dans la classe de conjugaison
du compact maximal standard, mais lui est conjugué par le groupe des
similitudes symplectiques. 
\end{rems}

\section{Restriction de la représentation de Weil à un tore maximal
  elliptique.}\label{6}

\subsection{}\label{6.3}

Pour $Sp(W)$, un tore est elliptique si et seulement s'il
est anisotrope, c'est à dire si et seulement s'il est compact. 
Dans ce paragraphe, nous montrons que la restriction de la
représentation de Weil à un tore maximal elliptique est sans
multiplicité. 

\begin{theo}\label{theo6.3.1}
Soit $T\subset Sp(W)$ un tore maximal elliptique.  La restriction de
la représentation de Weil $S_{\psi}$ à $T$ est somme directe de caractères
intervenant tous avec multiplicité $1$.
\end{theo}
\begin{dem}
Comme $T$ est un groupe commutatif et compact, il est clair que la
restriction de la représentation $S_{\psi}^{B}$ à $T$ est somme
directe de caractères. Il s'agit de montrer qu'ils interviennent tous
avec la multiplicité $1$.

Soit $\mathcal{H}_{\psi}^{\infty}$ l'espace des vecteurs
$\mathcal{C}^{\infty}$ de la représentation $S_{\psi}$, i.e. des
vecteurs $v\in\mathcal{H}_{\psi}$ tels que l'application $x\mapsto
S_{\psi}(x).v$ soit localement constante. C'est un sous-espace
vectoriel dense et $S_{\psi}$-invariant.

Si $G$ est un sous-groupe de $Sp(W)$, on désigne par $\widetilde{G}$
son image réciproque dans $Mp(W)$. Posons $G_{1}=G_{2}=T$. Puisque $T$
est son propre commutant dans $Sp(W)$, $(G_{1},G_{2})$ est une paire
réductive duale dans $Sp(W)$ selon Howe (voir
\cite{howe-ip-1979}). Comme $T$ est compact, il résulte du
théorème \ref{theo4.2.1} que $\widetilde{T}=T\times\{1,\epsilon\}$, où
$\epsilon$ est l'élément non trivial du noyau de la projection de
$Mp(W)$ sur $Sp(W)$. Soit $\chi$ un caractère unitaire de $T$ et soit
$\widetilde{\chi}$ le prolongement de $\chi$ à $\widetilde{T}$ non
trivial sur $\epsilon$. Désignons par $\mathbb{C}_{\chi}$
(resp. $\mathbb{C}_{\widetilde{\chi}}$) le $T$-module
(resp. $\widetilde{T}$-module) irréductible de dimension $1$
correspondant au caractère $\chi$ (resp. $\widetilde{\chi}$). La
conjecture de Howe, démontrée par Waldspurger (voir
\cite{waldspurger-ip-1990}), appliquée à la représentation
$\widetilde{\chi}\otimes\widetilde{\chi}$ de
$\tilde{G}_{1}\tilde{G}_{2}=\tilde{T}$ montre que l'espace des
opérateurs d'entrelacement
$Hom_{\widetilde{T}}(\mathcal{H}_{\psi}^{\infty},\mathbb{C}_{\widetilde{\chi}})$
est de dimension au plus $1$. Comme $S_{\psi}(\epsilon)=-Id$, on a
$\dim
Hom_{\widetilde{T}}(\mathcal{H}_{\psi}^{\infty},\mathbb{C}_{\widetilde{\chi}})=
\dim Hom_{T}(\mathcal{H}_{\psi}^{\infty},\mathbb{C}_{\chi})\geq \dim
Hom_{T}(\mathcal{H}_{\psi},\mathbb{C}_{\chi})$.
\end{dem}


\subsection{}\label{6.1}
Dans les paragraphes qui suivent, nous étudions la décomposition en
irréductibles de la restriction de la représentation de Weil à un tore
maximal elliptique de $Sp(W)$.

Soit $T$ un tore maximal elliptique de $Sp(W)$. D'après
\cite{morris-1991}, $W$ se décompose de manière unique comme une somme
directe orthogonale $W=W_{1}\oplus\cdots\oplus W_{n}$ de
sous-$T$-modules irréductibles sur $k$. Chaque $W_{i}$ est donc un
sous-espace symplectique de $W$ et l'image $T_{i}$ de $T$ par
l'application $x\mapsto x_{\vert W_{i}}$ est un tore maximal
elliptique de $Sp(W_{i})$ sous l'action duquel $W_{i}$ est un module
irréductible sur $k$. De plus, l'application $x\mapsto (x_{\vert
  W_{1}},\ldots,x_{\vert W_{n}})$ est un isomorphisme naturel du tore
$T$ sur le produit direct de tores $T_{1}\times\cdots\times
T_{n}$. Dans la suite, nous identifions $T$ et
$T_{1}\times\cdots\times T_{n}$ au moyen de cet isomorphisme.

Soit $1\leq i\leq n$. Le tore $T_{i}$ étant compact, il est contenu
dans un sous-groupe compact maximal  $K_{B_{i}}$ de $Sp(W_{i})$, avec
$B_{i}$ un bon réseau $T_{i}$-invariant de $W_{i}$. Alors
$B=B_{1}\oplus\cdots\oplus B_{n}$ est un bon réseau $T$-invariant de
$W$, de sorte que $T$ est contenu dans le sous-groupe compact maximal
$K_{B}$ de $Sp(W)$. On utilise la section $s_{B}$ (resp. $s_{B_{i}}$) du
paragraphe \ref{4.6} pour identifier $K_{B}$ (resp. $K_{B_{i}}$) à un
sous-groupe de $Mp(W)$ (resp. $Mp(W_{i})$). On a alors




\begin{pr}\label{pr6.3.1}
La restriction de la représentation de Weil au tore $T$ s'écrit comme
le produit tensoriel extérieur des restrictions pour chaque $i$ de la
représentation de Weil de $Mp(W_{i})$ à $T_{i}$ :
\begin{equation*}
  (S_{W,\psi})_{\vert T}=\otimes_{i=1}^{i=n}(S_{W_{i},\psi})_{\vert T_{i}}
\end{equation*}
\end{pr}
\begin{dem}
C'est une conséquence immédiate du lemme \ref{lem4.7.2} et de la
proposition \ref{pr4.7.1}.
\end{dem}

Il suit de la proposition \ref{pr6.3.1} que pour connaître la
restriction d'une représentation de Weil à un tore elliptique, il
suffit de connaître, pour tout espace symplectique $W$, la restriction
de la représentation de Weil de $Mp(W)$ aux tores maximaux $T$ de $Sp(W)$
pour lesquels le $T$-module $W$ est irréductible sur $k$. C'est ce que
nous allons faire dans les paragraphes \ref{6.1'} à \ref{7.3}.

\subsection{}\label{6.1'}
Dans la suite, si $F'$ est une extension de degré fini d'un corps $F$,
nous désignerons par $\tr_{F'/F}$ (resp. $\norm_{F'/F}$) la trace
(resp. la norme) de $F'$ sur $F$.

Pour simplifier nous dirons qu'un tore $T\subset Sp(W)$ est
irréductible si le $T$-module $W$ est irréductible sur $k$. 

Suivant \cite{howe-1973} et \cite{morris-1991}, nous allons donner une
description des classes de conjugaison de tores maximaux irréductibles
de $Sp(W)$.

Soit $k'$ une extension de degré $r$ de $k$ et $k''$ une extension
quadratique de $k'$. On note $\tau$ l'élément non trivial du groupe de
Galois de $k''$ sur $k'$. Soit $u\in k''$ un élément non nul et de trace
nulle relativement à $k'$. Alors, la formule
\begin{equation}\label{eq6.1.1}
\beta_{u}(x,y)=\frac{1}{2}\sideset{}{_{k''/k}}\tr ux^{\tau}y\mbox{, }x,y\in k'',
\end{equation}
définit une forme symplectique sur le $k$-espace vectoriel
$k''$. L'action de $k''$ sur lui-même par multiplication permet de
l'identifier à un sous-corps de l'algèbre $End_{k}(k'')$ des
endomorphismes de $k$-espace vectoriel de $k''$. On désigne par
$T_{k'',k'}$ le sous-groupe multiplicatif de $k''^{\times}$ constitué
des éléments de norme $1$ relativement à $k'$. Alors $T_{k'',k'}$ est
l'ensemble des $k$-points d'un tore maximal elliptique du groupe
symplectique $Sp(k'',\beta_{u})$.

On désigne par $Aut(k''/k)^{\tau}$ le sous-groupe du groupe des
automorphismes $Aut(k''/k)$ de $k''$ sur $k$ constitué des éléments
qui commutent à $\tau$. Le groupe $Aut(k''/k)^{\tau}$ laisse
invariant le sous-$k$-espace vectoriel de $k''$ constitué des éléments
de trace nulle ainsi que le sous-groupe multiplicatif
$\im(\norm_{k''/k'})$. On peut donc former le groupe produit
semi-direct
$\Gamma_{k'',k'}=Aut(k''/k)^{\tau}\ltimes\im(\norm_{k''/k'})$. L'action
de $Aut(k''/k)^{\tau}$ dans $\ker(\tr_{k''/k'})$ se prolonge à
$\Gamma_{k'',k'}$, le groupe $\im(\norm_{k''/k'})$ agissant par
multiplication. Dans le cas des tores maximaux irréductibles des groupes
symplectiques, les théorème 1.6 et proposition 1.10 de
\cite{morris-1991} s'énoncent ainsi~:

\begin{theo}\label{theo6.1.1}
(i) Soit $T\subset Sp(W)$ un tore maximal irréductible. Il existe une
  extension $k'$ de degré $r$ de $k$, une extension quadratique $k''$
  de $k'$, un élément $u\in k''^{\times}$ de trace nulle relativement
  à $k'$ et $\varphi$ un isomorphisme d'espaces symplectiques de
  $(k'',\beta_{u})$ sur $W$ tels que $T$ soit l'image de $T_{k'',k'}$
  par l'isomorphisme $\tilde{\varphi}:t\mapsto \varphi\circ
  t\circ\varphi^{-1}$. De plus, les tores
  $\tilde{\varphi}(T_{k'',k'})$, $\varphi$ parcourant l'ensemble des
  isomorphismes d'espaces symplectiques de $(k'',\beta_{u})$ sur $W$,
  forment une classe de conjugaison de tores maximaux dans $Sp(W)$,
  notée $\mathcal{C}_{k'',k',u}$.

(ii) Pour $i=1,2$, soit $k'_{i}$ (resp. $k''_{i}$) une extension de degré
  $r$ de $k$ (resp. quadratique de $k'_{i}$) et $u_{i}\in k''_{i}$ un
  élément non nul de trace nulle relativement à $k'_{i}$. Si les
  classes $\mathcal{C}_{k''_{i},k'_{i},u_{i}}$, $i=1,2$ sont
  identiques, il existe un isomorphisme sur $k$ de $k''_{1}$ sur
  $k''_{2}$ qui envoie $k'_{1}$ sur $k'_{2}$.

(iii) Soit $k'$ une extension de degré $r$ de $k$ et $k''$ une extension
quadratique de $k'$. Alors l'application
$u\mapsto\mathcal{C}_{k'',k',u}$ induit une bijection de
$\Gamma_{k'',k'}\backslash\ker(\tr_{k''/k'})$ sur l'ensemble des
classes de conjugaison des tores ma\-ximaux de $Sp(W)$ isomorphes sur $k$
à $T_{k'',k'}$.
\end{theo}

\subsection{}\label{6.2}
Dans ce paragraphe, pour chaque tore maximal irréductible $T$ de $Sp(W)$
nous exhibons un bon réseau $T$-invariant.

Soit $T\subset Sp(W)$ un tore maximal irréductible. D'après le théorème
\ref{theo6.1.1}, il existe une extension $k'$ de degré $r$ de $k$, une
extension quadratique $k''$ de $k'$ et un élément non nul
$u\in\ker(\tr_{k''/k'})$ tels que $W$ s'identifie au $k$-espace
symplectique $(k'',\beta_{u})$ et que $T=T_{k'',k'}$.

On note $\mathcal{O}$ (resp. $\mathcal{O}'$, $\mathcal{O}''$) l'anneau
des entiers, $\varpi$ (resp. $\varpi'$, $\varpi''$) une uniformisante,
$v$ (resp. $v'$, $v''$) la valuation normalisée et $q$ (resp. $q'$, $q''$) le
cardinal du corps résiduel de $k$ (resp. $k'$, $k''$). On prend
$\varpi''$ tel que $\varpi''^{2}=\varpi'$ lorsque $k''$ est ramifié
sur $k'$, et $\varpi''=\varpi'$ dans le cas contraire. On note $e$
l'indice de ramification de $k'$ sur $k$. On a alors
$q'=q^{\frac{r}{e}}$ et $q''=q'$ (resp. $q''=q'^{2}$) si $k''$ est
ramifié (resp. non ramifié) sur $k'$. 

On note $\delta$ l'entier tel que l'idéal
$\varpi'^{\delta}\mathcal{O}'$ soit la différente de $k'$ sur $k$,
i.e. $\delta$ est le plus grand entier tel que
$\tr_{k'/k}(x)\in\mathcal{O}$ pour tout
$x\in\varpi'^{-\delta}\mathcal{O}'$.

Compte tenu du théorème \ref{theo6.1.1} (ii), on peut supposer que
$u=\varpi''$ lorsque $k''$ est ramifié sur $k'$
et $v''(u)\in\{0,1\}$ dans le cas contraire~: c'est ce que nous
faisons désormais.

Comme $T$ est contenu dans $\mathcal{O}''$, les idéaux
$\varpi''^{m}\mathcal{O}''$, $m\in\mathbb{Z}$, sont des réseaux
$T$-invariants de $W$.
\begin{lem}\label{lem6.2.1}
Soit $m\in\mathbb{Z}$ et soit $B=\varpi''^{m}\mathcal{O}''$. Alors le
réseau dual de $B$ est donné par
\begin{equation*}
B^{*}=
\begin{cases}
\varpi''^{2(e\lambda_{\psi}-\delta-(m+1))}B &\mbox{si $k''$ est
    ramifié sur $k'$}\\
\varpi'^{e\lambda_{\psi}-\delta-v''(u)-2m}B &\mbox{si $k''$ n'est
 pas   ramifié sur $k'$}.\\
\end{cases}
\end{equation*}
\end{lem}
\begin{dem} Le réseau dual de $B$ est un $\mathcal{O}''$
  sous-module de $k''$, donc de la forme
  $B^{*}=\varpi''^{l}\mathcal{O}''$, avec $l$ un entier.

Si $k''$ est ramifié sur $k'$, on a $u\in\varpi''\mathcal{O}'^{\times}$
et $\mathcal{O}''=\mathcal{O}'+\mathcal{O}'\varpi''$. Il vient alors
\begin{equation*}
\begin{split}
\beta_{u}(B,\varpi''^{2(e\lambda_{\psi}-\delta-(m+1))}B) &=
\beta_{u}(\varpi''^{m}\mathcal{O}'',
\varpi''^{2(e\lambda_{\psi}-\delta-1)-m}\mathcal{O}'')\\ 
& =
\sideset{}{_{k''/k}}\tr\varpi''^{2(e\lambda_{\psi}-\delta-1)+1}\mathcal{O}''\\
& =
\sideset{}{_{k'/k}}\tr\varpi''^{2(e\lambda_{\psi}-\delta)}\mathcal{O}'\\
&=\varpi^{\lambda_{\psi}}\sideset{}{_{k'/k}}\tr\varpi'^{-\delta}\mathcal{O}'\\
&=\varpi^{\lambda_{\psi}}\mathcal{O},
\end{split}
\end{equation*}
d'où le résultat dans ce cas.

Si $k''$ n'est pas ramifié sur $k'$, on a
$\mathcal{O}''=\mathcal{O}'+\mathcal{O}'\varpi'^{-v''(u)}u$ et il
vient alors
\begin{equation*}
\begin{split}
\beta_{u}(B,\varpi'^{e\lambda_{\psi}-\delta-v''(u)-2m}B) &=
\beta_{u}(\varpi'^{m}\mathcal{O}'',
\varpi'^{e\lambda_{\psi}-\delta-v''(u)-m}\mathcal{O}'')\\
&=\sideset{}{_{k''/k}}\tr\varpi'^{e\lambda_{\psi}-\delta}\mathcal{O}''\\
&=\sideset{}{_{k'/k}}\tr\varpi'^{e\lambda_{\psi}-\delta}\mathcal{O}'\\
&=\varpi^{\lambda_{\psi}}\mathcal{O}, 
\end{split}
\end{equation*}
d'où le lemme.
\end{dem}
On pose $\mu=e\lambda_{\psi}-\delta-v''(u)$. On déduit du lemme
précédent le 
\begin{co}\label{co6.2.1}
On définit le réseau $B$ en posant 
\begin{equation*}
B=
\begin{cases}
\varpi''^{\mu}\mathcal{O}''\mbox{ si $k''$ est ramifié sur
  $k$,}\\ \varpi'^{[\frac{\mu+1}{2}]}\mathcal{O}''\mbox{ si $k''$ est
  non ramifié sur $k'$.}
\end{cases}
\end{equation*}
Alors $B$ est un bon réseau $T$-invariant, autodual si et seulement si
$k''$ est ramifié sur $k'$ ou si $\mu$ est pair.

Lorsque $k''$ est non ramifié sur $k'$ et $\mu$ est impair,
$B^{*}=\varpi'^{-1}B$ et
$\mathrm{b}^{*}=\mathcal{O}''/\varpi'\mathcal{O}''$ est un
$\mathbb{F}_{q}$ espace vectoriel symplectique de dimension
$\frac{2r}{e}$.
\end{co}

\subsection{}\label{7.1}
On se donne un tore maximal irréductible $T$ de $Sp(W)$ et on reprend
les notations du paragraphe précédent. Nous allons décrire les
caractères de $T$.

Lorsque $k''$ n'est pas ramifié sur $k'$, on désigne par $\nu$
l'élément de $\mathcal{O}''^{\times}$ tel que $u=\varpi'^{v''(u)}\nu$
et on pose $d=\nu^{2}$, qui est un élément de $\mathcal{O}'$. Alors,
\begin{equation*}
T=\{\xi+\eta\nu\vert \xi,\eta\in\mathcal{O}'\mbox{ et } 
\xi^{2}-d\eta^{2}=1\}.
\end{equation*}

Lorsque $k''$ est ramifié sur $k'$, on prend $u=\varpi''$ et on a 
\begin{equation*}
T=\{\xi+\eta\varpi''\vert \xi,\eta\in\mathcal{O}'\mbox{ et } 
\xi^{2}-\varpi'\eta^{2}=1\}.
\end{equation*}

Si $j\in\mathbb{N}$, on désigne par $T_{j}$ le $j$-ième sous-groupe de
congruence de $T$~:
\begin{equation*}
T_{j}=\{g\in T\vert g-1\in\varpi''^{j}\mathcal{O}''\}.
\end{equation*}
On a $T_{0}=T$ et, pour $j\geq1$, $T_{j}$ est un sous-groupe strict de
$T$ que nous allons décrire.

On désigne par $\mu_{q'+1}$ le sous-groupe cyclique de $k''^{\times}$
constitué des racines $(q'+1)$-ièmes de l'unité.

\begin{lem}\label{lem7.1.1}
a) On suppose que $k''$ n'est pas ramifié sur $k'$. 

(i) On a $T=\mu_{q'+1}\times T_{1}$. La suite $T_{j}$,
$j\in\mathbb{N}$, est une suite strictement décroissante de
sous-groupes de $T$.

Soit $j$ un entier naturel non nul. 

(ii) Pour tout $\eta\in\mathcal{O}'$, il existe un unique élément
$\theta_{j}(\eta)$ de $T_{j}$ tel que
$\theta_{j}(\eta)-\varpi'^{j}\eta\nu\in k'$. Cet élément vérifie que
$\theta_{j}(\eta)-\varpi'^{j}\eta\nu-1\in\varpi'^{2j}\mathcal{O}'$ et
$\theta_{j}(\eta)^{-1}=\theta_{j}(-\eta)$.

(iii) L'application $\theta_{j}$ est une bijection de $\mathcal{O}'$
sur $T_{j}$ et elle induit, pour tout entier $1\leq s\leq2j$, un
isomorphisme de groupes de $\mathcal{O}'/\varpi'^{s}\mathcal{O}'$ sur
$T_{j}/T_{j+s}$.

b) On suppose que $k''$ est ramifié sur $k'$. 

(i) On a $T=\{\pm1\}\times T_{1}$ et $T_{2j}=T_{2j+1}$, pour tout
entier $j>0$. La suite $(T_{2j+1})_{j\in\mathbb{N}}$ est une suite
strictement décroissante de sous-groupes de $T$.

Soit $j$ un entier naturel.

(ii) Pour tout $\eta\in\mathcal{O}'$, il existe un unique élément
$\theta_{j}(\eta)$ de $T_{2j+1}$ tel que
$\theta_{j}(\eta)-\varpi'^{j}\eta\varpi''\in k'$. Cet élément vérifie
que
$\theta_{j}(\eta)-\varpi'^{j}\eta\varpi''-1\in\varpi'^{2j+1}\mathcal{O}'$
et $\theta_{j}(\eta)^{-1}=\theta_{j}(-\eta)$.

(iii) L'application $\theta_{j}$ est une bijection de $\mathcal{O}'$
sur $T_{2j+1}$ et elle induit, pour tout entier $1\leq s\leq2j+1$, un
isomorphisme de groupes de $\mathcal{O}'/\varpi'^{s}\mathcal{O}'$ sur
$T_{2j+1}/T_{2(j+s)}$.

\end{lem}
\begin{dem}
a) On se place dans le cas où $k''$ n'est pas ramifié sur $k'$. Le
fait que $T=\mu_{q'+1}\times T_{1}$ et que la suite $T_{j}$ est
décroissante est immédiat. Le fait qu'elle l'est strictement résultera
alors du point (iii).

Soit $j\in\mathbb{N}$ non nul et $\eta\in\mathcal{O}'$. Il résulte du
lemme de Hensel qu'il existe un unique $\xi\in\varpi'\mathcal{O}'$ tel
que $1+\xi+\varpi'^{j}\eta\nu$ soit un élément de $T$, c'est à dire
vérifie la relation $(1+\xi)^{2}-d\varpi'^{2j}\eta^{2}=1$. Il suit de
celle-ci que $\xi\in\varpi'^{2j}\mathcal{O}'$. D'où l'assertion (ii)
et le fait que $\theta_{j}$ est une bijection de $\mathcal{O}'$ sur
$T_{j}$. Enfin, pour $\eta\in\mathcal{O}'$, on a
$\theta_{j}(\eta)^{-1}=\theta_{j}(\eta)^{\tau}=\theta_{j}(-\eta)$.

Pour $\eta\in\mathcal{O}'$, on a $\theta_{j}(\eta)\in
1+\varpi'^{j}\eta\nu+\varpi'^{2j}\mathcal{O}'$. On en déduit que pour
$\eta,\eta'\in\mathcal{O}'$,
$\theta_{j}(\eta)\theta_{j}(\eta')\in1+\varpi'^{j}(\eta+\eta')\nu+
\varpi'^{2j}\mathcal{O}'+\varpi'^{3j}\mathcal{O}''$, puis
$\theta_{j}(\eta)\theta_{j}(\eta')\theta_{j}(\eta+\eta')^{-1}\in
1+\varpi'^{2j}\mathcal{O}'+\varpi'^{3j}\mathcal{O}''$. La dernière
affirmation de l'assertion (iii) en découle.

b) On se place dans le cas où $k''$ est ramifié sur $k'$. Il est clair
que la suite $T_{j}$ est décroissante et que $T=\{\pm1\}\times
T_{1}$. Soit $j>0$ un entier et $x\in T_{2j}$. Alors, on a
$x=1+\varpi'^{j}(\xi+\eta\varpi'')$, avec $\xi,\eta\in\mathcal{O}'$
vérifiant la relation $
(1+\varpi'^{j}\xi)^{2}=1+\varpi'^{2j+1}\eta^{2}$. On en déduit
immédiatement que $\xi\in\varpi'^{j+1}\mathcal{O}'$, d'où le fait que
$x\in 1+\varpi'^{j+1}\mathcal{O}'+\varpi'^{j}\varpi''\mathcal{O}'
=1+\varpi''^{2j+1}\mathcal{O}''$, i.e. $x\in T_{2j+1}$. On a bien
$T_{2j}=T_{2j+1}$. A partir de là, la démonstration des assertions
restantes est identique à celle faite dans le cas précédent.
\end{dem}

Soit $j\in\mathbb{N}$ et $\chi$ un caractère de $T_{j}$. On appelle
conducteur de $\chi$, le plus petit entier $\lambda$ tel que
$T_{\lambda}$ soit contenu dans le noyau de $\chi$.

Le seul caractère de conducteur $0$ est le caractère trivial.

Lorsque $k''$ n'est pas ramifié sur $k'$, les caractères de conducteur
au plus 1 de $T$ sont les relèvements à $T$ des caractères de
$T/T_{1}=\mu_{q'+1}$. Ils forment donc un sous-groupe cyclique d'ordre
$q'+1$ de $\hat{T}$.

Soit $(\, ,)_{k'}$ le symbole de Hilbert du corps $k'$. On sait que
tout élément $g$ de $T$ s'écrit $g=z(z^{-1})^{\tau}$ avec
$z\in\mathcal{O}''^{\times}$. On vérifie que le nombre
$(\varpi',zz^{\tau})_{k'}=
\left(\frac{p_{\mathbb{F}_{q'}}(zz^{\tau})}{q'}\right)$ ne dépend
que de $g$ et que la formule
\begin{equation}\label{eq7.1.0}
\eta_{0}(g)=(\varpi',zz^{\tau})_{k'}=
\left(\frac{p_{\mathbb{F}_{q'}}(zz^{\tau})}{q'}\right)
\end{equation}
définit un caractère de conducteur 1 de $T$.

Lorsque $k''$ est ramifié sur $k'$, $T$ possède un unique caractère de
conducteur $1$~: il s'agit du caractère $\epsilon$ trivial sur $T_{1}$
et tel que $\epsilon(-1)=-1$. Tous les autres caractères de $T$ sont
de conducteur pair.
\begin{lem}\label{lem7.1.2}
a) On suppose que $k''$ n'est pas ramifié sur $k'$.

(i) Soit $j\geq1$ un entier. Alors, pour tout entier $\lambda$ tel que
$j<\lambda\leq3j$ et pour tout $b\in\mathcal{O}'$, la formule
\begin{equation}\label{eq7.1.1}
\chi_{b,\lambda,j}(t)=
\psi(-\frac{1}{4}\sideset{}{_{k''/k}}\tr\varpi'^{\mu-\lambda}but)
\mbox{, }t\in T_{j}
\end{equation}
définit un caractère de $T_{j}$ de conducteur inférieur ou égal à
$\lambda$.

Le caractère $\chi_{b,\lambda,j}$ est de conducteur $\lambda$ si et
seulement si $b\in\mathcal{O}'^{\times}$, et tout caractère de $T_{j}$
de conducteur $\lambda$ est de cette forme.

Si $a,b\in\mathcal{O}'^{\times}$, on a 
\begin{equation}\label{eq7.1.2}
\begin{split}
\chi_{a,\lambda,j}=\chi_{b,\lambda,j}&\Longleftrightarrow 
b-a\in\varpi'^{\lambda-j}\mathcal{O}'\\
&\Longleftrightarrow 
b/a\in1+\varpi'^{\lambda-j}\mathcal{O}'.
\end{split}
\end{equation}

(ii) Soit $\chi$ un caractère de $T$ de conducteur $\lambda\geq2$ et
soit $j$ un entier tel que $j<\lambda\leq3j$. Alors $\chi$ est un
prolongement à $T$ d'un caractère  $\chi_{b,\lambda,j}$ pour un
certain $b\in\mathcal{O}'^{\times}$.

b) On suppose que $k''$ est ramifié sur $k'$.

(i) Soit $j\geq0$ un entier. Alors, pour tout entier $\lambda$ tel que
$j<\lambda\leq3j+1$ et pour tout $b\in\mathcal{O}'$, la formule
\begin{equation}\label{eq7.1.3}
\chi_{b,2\lambda,j}(t)=
\psi(-\frac{(-1)^{\lambda-\mu}}{4}
\sideset{}{_{k''/k}}\tr\varpi'^{\mu-\lambda}b\varpi''t)
\mbox{, }t\in T_{2j+1}
\end{equation}
définit un caractère de $T_{2j+1}$ de conducteur inférieur ou égal à
$2\lambda$.

Le caractère $\chi_{b,2\lambda,j}$ est de conducteur $2\lambda$ si et
seulement si $b\in\mathcal{O}'^{\times}$, et tout caractère de $T_{2j+1}$
de conducteur $2\lambda$ est de cette forme.

Si $a,b\in\mathcal{O}'^{\times}$, on a 
\begin{equation}\label{eq7.1.4}
\begin{split}
\chi_{a,2\lambda,j}=\chi_{b,2\lambda,j}&\Longleftrightarrow 
b-a\in\varpi'^{\lambda-j}\mathcal{O}'\\
&\Longleftrightarrow
b/a\in1+\varpi'^{\lambda-j}\mathcal{O}'. 
\end{split}
\end{equation}

(ii) Soit $\chi$ un caractère de $T$ de conducteur $2\lambda\geq2$ et
soit $j$ un entier tel que $j<\lambda\leq3j+1$. Alors $\chi$ est un
prolongement à $T$ d'un caractère $\chi_{b,2\lambda,j}$ pour un certain
$b\in\mathcal{O}'^{\times}$.
\end{lem}
\begin{dem}
a) Si $\eta\in\mathcal{O}'$, on a 
\begin{equation*}
\chi_{b,\lambda,j}(\theta_{j}(\eta))=
\psi(-\frac{1}{4}\sideset{}{_{k''/k}}
\tr\varpi'^{e\lambda_{\psi}-\delta+j-\lambda}bd\eta).
\end{equation*}
On en déduit que $\chi_{b,\lambda,j}\circ\theta_{j}$ est un caractère
du groupe additif $\mathcal{O}'$ dont le noyau contient l'idéal
$\varpi'^{\lambda-j}\mathcal{O}'$. De plus,
$\varpi'^{\lambda-j}\mathcal{O}'$ est le plus grand idéal de
$\mathcal{O}'$ contenu dans $\ker\chi_{b,\lambda,j}\circ\theta_{j}$ si et
seulement si $b\in\mathcal{O}'^{\times}$. Il suit alors du lemme
\ref{lem7.1.1} que $\chi_{b,\lambda,j}$ est un caractère de $T_{j}$ et
il est alors clair que son conducteur est inférieur ou égal à $\lambda$,
l'égalité ayant lieu si et seulement si $b\in\mathcal{O}'^{\times}$.

Réciproquement, si $\chi$ est un caractère de conducteur $\lambda$ tel
que $j<\lambda\leq 3j$, il suit du lemme \ref{lem7.1.1} que
$\chi\circ\theta_{j}$ est un caractère du groupe additif
$\mathcal{O}'$ et que $\varpi'^{\lambda-j}\mathcal{O}'$ est le plus
grand idéal de $\mathcal{O}'$ contenu dans son noyau. Il est alors
clair qu'il existe $b\in\mathcal{O}'^{\times}$ tel que
$\chi=\chi_{b,\lambda,j}$.

Si $a,b\in\mathcal{O}'^{\times}$, on a 
\begin{equation*}
\begin{split}
\chi_{a,\lambda,j}=\chi_{b,\lambda,j}&\Longleftrightarrow
\chi_{a,\lambda,j}\circ\theta_{j}=\chi_{b,\lambda,j}\circ\theta_{j}\\
&\Longleftrightarrow\forall\eta\in\mathcal{O}'\mbox{, }
\psi(\sideset{}{_{k'/k}}\tr\varpi'^{e\lambda_{\psi}-\delta+j-\lambda}(b-a)d\eta)=1\\
&\Longleftrightarrow b-a\in\varpi'^{\lambda-j}\mathcal{O}'\\
&\Longleftrightarrow
b/a\in1+\varpi'^{\lambda-j}\mathcal{O}'. 
\end{split}
\end{equation*}
D'où l'assertion (i).

L'assertion (ii) est une conséquence immédiate de (i).

b) La démonstration est identique à celle du cas a), une fois que l'on a
remarqué que, pour $\eta\in\mathcal{O}'$,
\begin{equation*}
\chi_{b,2\lambda,j}(\theta_{j}(\eta))=
\psi
(\frac{(-1)^{\lambda-\mu}}{2}
\sideset{}{_{k'/k}}\tr\varpi'^{e\lambda_{\psi}-\delta+j-\lambda}b\eta).
\end{equation*}
\end{dem}

\subsection{}\label{7.2}
Soit $T\subset Sp(W)$ un tore maximal irréductible et $B$ le bon
réseau $T$-invariant construit au paragraphe \ref{6.2}. Alors
$T\subset K_{B}$ et on relève $T$ en un sous-groupe de $Mp(W)$ à
l'aide du relèvement $s_{B}$ de $K_{B}$ dans $Mp(W)$ défini dans le paragraphe
\ref{4.6} (voir théorème \ref{theo4.6.1} et la remarque en fin de ce
paragraphe).

Dans ce paragraphe, on suppose que le $\mathcal{O}''$-module $B$ est
un réseau autodual et on considère la représentation de Weil
$S^{B}_{\psi}$ dans l'espace $\mathcal{H}^{B}_{\psi}$. On rappelle les
fonctions $\delta_{v}\in\mathcal{H}^{B}_{\psi}$, $v\in k''$ définies
par la formule \ref{eq4.4.1}.

\begin{theo}\label{theo7.2.1}
a) La fonction $\delta_{0}$ est, à un scalaire multiplicatif près,
  l'unique vecteur de poids trivial de $T$.

b) On suppose que $k''$ n'est pas ramifié sur $k'$. 

(i) Soit $\chi$ un caractère non trivial de $T$. Alors $\chi$
intervient dans la représentation de Weil $S^{B}_{\psi}$ si et
seulement si son conducteur est pair.

(ii) Soit $\chi$ un caractère de $T$ de conducteur pair $2j>0$. Alors,
il existe $a\in\mathcal{O}''^{\times}$ tel que 
\begin{equation*}
\chi_{\vert T_{j}}=\chi_{aa^{\tau},2j,j}
\end{equation*}
et le sous-espace des vecteurs de poids $\chi$ est engendré par la
fonction
\begin{equation}\label{eq7.2.1}
  \varphi_{\chi}=\sum_{g\in T/T_{j}}\chi(g^{-1})\delta_{\varpi'^{\frac{\mu}{2}-j}ga}.
\end{equation}

c) On suppose que $k''$ est ramifié sur $k'$.

Soit $\chi$ un caractère non trivial de $T$. Alors $\chi$
intervient dans la représentation de Weil $S^{B}_{\psi}$ si et
seulement si son conducteur est un entier pair $2j$ et il existe
$a\in\mathcal{O}''^{\times}$ tel que 
\begin{equation*}
\chi_{\vert T_{j}}=\chi_{aa^{\tau},j,[\frac{j}{2}]}.
\end{equation*}
Dans ce cas, le sous-espace des vecteurs
de poids $\chi$ est engendré par la fonction
\begin{equation}\label{eq7.2.2}
\varphi_{\chi}=\sum_{g\in T/T_{j}}\chi(g^{-1})\delta_{\varpi''^{\mu-j}ga}.
\end{equation}
\end{theo}
\begin{dem}
Il résulte du théorème \ref{theo6.3.1} que $\delta_{0}$ est, à un
scalaire multiplicatif près, l'unique vecteur de poids trivial de
$T$.

On pose $\lambda=\frac{\mu}{2}$ si $k''$ n'est pas ramifié sur $k'$ et
$\lambda=\mu$ dans le cas contraire. Soit $x\in k''$ tel que $v''(x)< \lambda$
et $\dot{x}$ son image dans $k''/B$. On écrit $x=\varpi''^{v''(x)}a$
avec $a\in\mathcal{O}''^{\times}$. Il est immédiat que le
stabilisateur $T(\dot{x})$ de $\dot{x}$ dans $T$ est le sous-groupe de
congruence $T_{\lambda-v''(x)}$.

La représentation $\sigma_{\dot{x}}$ de $K_{B}(\dot{x})$ définie par
la formule $\ref{eq5.3.2}$ est un caractère (comme $B$ est un réseau
autodual, $\dot{x}$ et $\hat{x}$ sont le même élément de
$k''/B=k''/B^{*}$). On note $\chi_{\dot{x}}$ sa restriction à
$T(\dot{x})$. Alors, on a
\begin{eqnarray*}
  \chi_{\dot{x}}&=&\chi_{aa^{\tau},2(\lambda-v''(x)),\lambda-v''(x)}
\mbox{, si $k''$ n'est pas ramifié sur $k'$,}\\
\chi_{\dot{x}}&=&\chi_{aa^{\tau},\lambda-v''(x),[\frac{\lambda-v''(x)}{2}]}
\mbox{, si $k''$ est  ramifié sur $k'$.}
\end{eqnarray*}
Soit $S_{T,\dot{x}}$ la représentation de $T$ définie par 
\begin{equation*}
S_{T,\dot{x}}=\Ind_{T(\dot{x})}^{T}\chi_{\dot{x}}.
\end{equation*}
L'espace $\mathcal{H}_{T,\dot{x}}$ de cette représentation est
l'espace des fonctions $\varphi:T\mapsto\mathbb{C}$ qui vérifient
$\varphi(gt)=\chi_{\dot{x}}(t)\varphi(g)$, $g\in T$, $t\in T(\dot{x})$.  Si
$\varphi\in\mathcal{H}_{T,\dot{x}}$, il existe un unique élément de
$\mathcal{H}^{B}_{\psi}$, $\mathcal{I}_{T,\dot{x}}\varphi$, de support
contenu dans la $T$-orbite de $\dot{x}$  vérifiant
\begin{equation*}
  \mathcal{I}_{T,\dot{x}}\varphi(g.x)=\varphi(g^{-1})\mbox{,
  }g\in T.
\end{equation*}

Soit $X\subset k''/B\backslash\{0\}$ un ensemble de représentants des
$T$-orbites dans $k''/B\backslash\{0\}$. Il résulte des lemmes
\ref{lem5.4.1} et \ref{lem3.4.1} que l'on a
\begin{equation*}
S^{B}_{\psi\vert T}=\mathbb{C}\delta_{0}\bigoplus_{\dot{x}\in X}S_{T,\dot{x}}
\end{equation*}
et que, pour $\dot{x}\in k''/B\backslash\{0\}$, l'application
$\mathcal{I}_{T,\dot{x}}$ est un isomorphisme de $T$-modules de
$\mathcal{H}_{T,\dot{x}}$ sur le sous-espace de
$\mathcal{H}^{B}_{\psi}$ constitué des fonctions à support dans la
$T$-orbite de $\dot{x}$. De plus, on a
\begin{equation*}
  S_{T,\dot{x}}=
\bigoplus_{\stackrel{\chi\in\hat{T}}{\chi_{\vert T(\dot{x})}=\chi_{\dot{x}}}}\chi.
\end{equation*}
Il suit de ces considérations que les caractères qui interviennent
dans $S^{B}_{\psi\vert T}$ sont, outre le caractère trivial, ceux
de conducteur $2j>0$ dont la restriction à $T_{j}$ est de la forme
$\chi_{aa^{\tau},2j,j}$ (resp. $\chi_{aa^{\tau},j,[\frac{j}{2}]}$)
avec $a\in\mathcal{O}''^{\times}$ si $k''$ n'est pas ramifié sur $k'$
(resp. $k''$ est ramifié sur $k'$).

Maintenant soit $\chi\in\hat{T}$ un caractère de conducteur $2j>0$ et
$a\in\mathcal{O}''^{\times}$ tels que $\chi_{\vert
  T_{j}}=\chi_{aa^{\tau},2j,j}$ (resp. $\chi_{\vert
  T_{j}}=\chi_{aa^{\tau},j,[\frac{j}{2}]}$) si $k''$ n'est pas ramifié
sur $k'$ (resp. $k''$ est ramifié sur $k'$). Soit
$x=\varpi''^{\lambda-j}a\in k''\backslash B$. Soit
$\varphi\in\mathcal{H}_{T,\dot{x}}$ la fonction de support inclus dans
$T_{j}$ telle que $\varphi(1)=1$. Alors la fonction $\varphi_{\chi}$
définie par
\begin{equation*}
\varphi_{\chi}=\sum_{g\in T/T_{j}}\chi(g^{-1})S_{T,\dot{x}}(g)\varphi,
\end{equation*}
est, à un scalaire multiplicatif près, l'unique vecteur de poids
$\chi$ de la représentation $S_{T,\dot{x}}$. Or, on a clairement
\begin{equation*}
\mathcal{I}_{T,\dot{x}}\varphi=\delta_{x}.
\end{equation*}
Par suite, la fonction 
\begin{equation*}
\begin{split}
  \mathcal{I}_{T,\dot{x}}\varphi_{\chi}&=
\sum_{g\in T/T_{j}}\chi(g^{-1})S^{B}_{\psi}(g)\delta_{x}\\
&=\sum_{g\in T/T_{j}}\chi(g^{-1})\delta_{gx}\\
&=\sum_{g\in T/T_{j}}\chi(g^{-1})\delta_{\varpi''^{\lambda-j}ga}
\end{split}
\end{equation*}
est un vecteur de poids $\chi$ de la représentation $S^{B}_{\psi\vert
  T}$.

Enfin, si $k''$ n'est pas ramifié sur $k'$,
$\sideset{}{_{k''/k'}}\norm$ est un morphisme surjectif de groupes de
$\mathcal{O}''^{\times}$ sur $\mathcal{O}'^{\times}$. On en déduit que
si $\chi\in\hat{T}$ est un caractère de conducteur $2j>0$, il existe
$a\in\mathcal{O}''^{\times}$ tels que $\chi_{\vert
  T_{j}}=\chi_{aa^{\tau},2j,j}$.
\end{dem}

\subsection{}\label{7.3}
On conserve les notations introduites au tout début du paragraphe
précédent, mais on se place dans la situation où $k''$ n'est pas
ramifié sur $k'$ et $\mu$ est impair. Alors
$B=\varpi'^{\frac{\mu+1}{2}}\mathcal{O}''$ est un bon réseau
$T$-invariant dont le réseau dual est
$B^{*}=\varpi'^{\frac{\mu-1}{2}}\mathcal{O}''$. De plus
$A=\varpi'^{\frac{\mu-1}{2}}(\mathcal{O}'+\mathcal{O}'\varpi'\nu)$ est
un réseau autodual tel que $B\subset A\subset B^{*}$. On considère la
représentation de Weil $S_{\psi}^{A}$ dans l'espace
$\mathcal{H}_{\psi}^{A}$. On rappelle le caractère $\eta_{0}$ de $T$
défini au paragraphe \ref{7.1} par la formule \ref{eq7.1.0}.

\begin{theo}\label{theo7.3.1}
a) Un caractère $\chi$ de conducteur au plus $1$ de $T$ intervient
dans la représentation $S_{\psi}^{A}$ si et seulement si
$\chi\neq\eta_{0}$. Dans ce cas, il existe $v\in B^{*}$ tel que la fonction
\begin{equation}\label{eq7.3.1}
  \varphi_{\chi}=\sum_{g\in T/T_{1}}\chi(g^{-1})S_{\psi}^{A}(g)\delta_{v}
\end{equation}
soit non nulle et engendre le sous-espace des vecteurs de poids $\chi$
de $\mathcal{H}_{\psi}^{A}$.

b) (i) Un caractère $\chi$ de conducteur strictement supérieur à $1$
de $T$ intervient dans la représentation $S_{\psi}^{A}$ si et
seulement si son conducteur est impair.

(ii) Soit $\chi$ un caractère de $T$ de conducteur impair
$2j+1>1$. Alors, il existe $a\in\mathcal{O}''^{\times}
\backslash(\varpi'\mathcal{O}'+\mathcal{O}'\nu)$ tel que 
\begin{equation*}
  \chi_{\vert T_{j+1}}=\chi_{aa^{\tau},2j+1,j+1}
\end{equation*}
et le sous-espace des vecteurs de poids $\chi$ de
$\mathcal{H}_{\psi}^{A}$ est engendré par la fonction
\begin{equation}\label{eq7.3.2}
  \varphi_{\chi}=
\sum_{g\in T/T_{j+1}}\chi(g^{-1})S_{\psi}^{A}(g)\delta_{\varpi'^{\frac{\mu-1}{2}-j}a}.
\end{equation}
\end{theo}
\begin{dem}
D'après le théorème \ref{theo6.3.1}, on sait déjà que les caractères
de $T$ interviennent tous sans multiplicité. Il suffit donc de montrer
que les seuls caractères qui interviennent sont ceux indiqués, avec les
fonctions propres correspondantes. 

On réalise la représentation de Weil $S_{\overline{\psi}}$ de
$Sp(\mathrm{b}^{*})$ dans l'espace $\mathcal{H}^{\mathrm{x}}_{\overline{\psi}}$
où $\mathrm{x}=A/B$. Cet espace s'identifie naturellement via l'application
$\varphi\mapsto\varphi\circ p_{\mathrm{b}^{*}}$, au sous-espace
$K_{B}$-invariant $\mathcal{E}^{B}_{0}$ de $\mathcal{H}_{\psi}^{A}$
constitué des fonctions à support dans $B^{*}$ (en particulier, une
base de $\mathcal{H}^{\mathrm{x}}_{\overline{\psi}}$ est constituée
des fonctions $\delta_{v}$, $v$ parcourant un système de représentants
dans $B^{*}$ de $B^{*}/A$). Dans ces conditions la
représentation $\widetilde{S}_{\overline{\psi}}$ est la représentation
de $K_{B}$ dans $\mathcal{E}^{B}_{0}$ induite par $S_{\psi}^{A}$.

On pose $\lambda=\frac{\mu-1}{2}$. Soit $x\in k''$ tel que
$v''(x)\leq\lambda$, $\dot{x}$ son image dans $k''/B$, $\ddot{x}$
celle dans $k''/A$ et $\hat{x}$ celle dans $k''/B^{*}$.  On désigne
par $\sigma_{T,\dot{x}}$ la restriction au stabilisateur $T(\hat{x})$
de $\hat{x}$ dans $T$ de la représentation $\sigma_{\dot{x}}$ de
$K_{B}(\hat{x})$ dans $\mathcal{H}^{\mathrm{x}}_{\overline{\psi}}$
définie par la formule $\ref{eq5.3.2}$.  Soit $S_{T,\dot{x}}$ la
représentation de $T$ définie par
\begin{equation*}
  S_{T,\dot{x}}=\Ind_{T(\hat{x})}^{T}\sigma_{T,\dot{x}}.
\end{equation*}
L'espace $\mathcal{H}_{T,\dot{x}}$ de cette représentation est
l'espace des fonctions $\varphi
:T\longrightarrow\mathcal{H}^{\mathrm{x}}_{\overline{\psi}}$ qui
vérifient $\varphi(gt)=\sigma_{T,\dot{x}}(t)\varphi(g)$, $g\in T$,
$t\in T(\hat{x})$. Si $\varphi\in\mathcal{H}_{T,\dot{x}}$, il existe
un unique élément $\mathcal{I}_{T,\dot{x}}\varphi$ de
$\mathcal{H}_{\psi}^{B}$, de support contenu dans la $T$-orbite de
$\hat{x}$ et vérifiant
\begin{equation}\label{eq7.3.3}
  \mathcal{I}_{T,\dot{x}}\varphi(g.x)=
\widetilde{S}_{\overline{\psi}}(g)\varphi(g^{-1})
  \mbox{, }g\in T.
\end{equation}
Soit $X\subset k''\backslash B^{*}$ un ensemble de représentants
des orbites de $T$ dans $k''/B^{*}\backslash\{0\}$. Il résulte des lemmes
\ref{lem5.4.1} et \ref{lem3.4.1} que l'on a
\begin{equation*}
  S_{\psi\vert
    T}^{B}=\widetilde{S}_{\overline{\psi}\vert T}\bigoplus_{x\in X}S_{T,\dot{x}}
\end{equation*}
et que, pour $x\in X\cup\{0\}$, l'application
$\mathcal{I}_{T,\dot{x}}$ est un isomorphisme de $T$-modules de
$\mathcal{H}_{T,\dot{x}}$ sur le sous-espace $\mathcal{G}_{\hat{x}}$
de $\mathcal{H}_{\psi}^{B}$ constitué des fonctions à support dans la
$T$-orbite de $\hat{x}$. On en déduit évidemment que 
\begin{equation}\label{eq7.3.4}
  S_{\psi\vert
    T}^{A}=\widetilde{S}_{\overline{\psi}\vert T}\bigoplus_{x\in X}S_{T,\dot{x}}
\end{equation}
l'application $\mathcal{I}_{A,B}\circ\mathcal{I}_{T,\dot{x}}$ étant,
pour tout $x\in X\cup\{0\}$, un isomorphisme de $T$-modules de
$\mathcal{H}_{T,\dot{x}}$ sur le sous-espace $\mathcal{E}_{\hat{x}}$
de $\mathcal{H}_{\psi}^{A}$ constitué des fonctions à support dans la
$T$-orbite de $\hat{x}$.

L'image de
$\varpi'\mathcal{O}'+\mathcal{O}'\nu=\varpi'^{-\frac{\mu-1}{2}}A\nu$
dans $\mathbb{F}_{q''}$ étant
$\mathbb{F}_{q'}p_{\mathbb{F}_{q''}}(\nu)$, on voit que la $T$-orbite
de tout élément de $\mathcal{O}''^{\times}$ rencontre
$\mathcal{O}''^{\times}\backslash(\varpi'\mathcal{O}'+\mathcal{O}'\nu)$. On
peut donc supposer que $X\subset k'^{\times}
(\mathcal{O}''^{\times}\backslash(\varpi'\mathcal{O}'+\mathcal{O}'\nu))$.

Soit $x\in X$. \'Ecrivons $x=\varpi'^{v''(x)}a$ avec $v''(x)<\lambda$
et
$a\in\mathcal{O}''^{\times}\backslash(\varpi'\mathcal{O}'+\mathcal{O}'\nu)$.
Alors le stabilisateur $T(\hat{x})$ de $\hat{x}$ dans $T$ est le
sous-groupe de con\-gruen\-ce $T_{j}$ avec $j=\lambda-v''(x)$ tandis que
celui $T(\dot{x})$ de $\dot{x}$ est $T_{j+1}$.

La représentation $\sigma_{T,\dot{x}}$ de $T(\hat{x})$ dans
$\mathcal{H}^{\mathrm{x}}_{\overline{\psi}}$ est donnée par
\begin{equation}\label{eq7.3.5}
  \sigma_{T,\dot{x}}(g)=
\psi(\frac{1}{2}\beta(gx,x))\widetilde{\rho}_{\overline{\psi}}(g^{-1}x-x)
\mbox{, }g\in T(\hat{x}). 
\end{equation}
D'autre part, le groupe $T(\hat{x})$ laisse stable le réseau autodual
$A$ et il opère donc dans le quotient $k''/A$. On désigne par
$T(\ddot{x})$ le stabilisateur de $\ddot{x}$ dans $T(\hat{x})$. Alors,
on a 
\begin{equation*}
  T(\ddot{x})=T(\dot{x})=T_{j+1}.
\end{equation*}
En effet, il est clair que $T(\dot{x})\subset
T(\ddot{x})$. Réciproquement, soit $g\in T(\hat{x})$ tel que
$g(x+A)=x+A$. Alors, on a $g=1+\varpi'^{2j}\xi+\varpi'^{j}\eta\nu$,
avec $\xi$, $\eta\in\mathcal{O}'$. Par hypothèse
$x=\varpi'^{\lambda-j}(u+v\nu)$ avec $u\in\mathcal{O}'^{\times}$ et
$v\in\mathcal{O}'$. La condition $g\ddot{x}=\ddot{x}$ s'écrit
$(g-1)(u+v\nu)\in\varpi'^{j}(\mathcal{O}'+\mathcal{O}'\varpi'\nu)$
qui, après calcul, entraîne
$u\eta+\varpi'^{j}v\xi\in\varpi'\mathcal{O}'$. Comme $u$ est dans
$\mathcal{O}'^{\times}$ et $j>0$, on en déduit que
$\eta\in\varpi'\mathcal{O}'$ et donc que $g\in T_{j+1}$.

Cela dit, si $g\in T(\dot{x})$, on a 
\begin{equation*}
xgx^{-1}=g(g^{-1}x-x,\frac{1}{2}\beta(x,g^{-1}x))\in
K'_{B}B(0,\frac{1}{2}\beta(x,g^{-1}x)).
\end{equation*} 
On en déduit que la
restriction de la représentation $\sigma_{T,\dot{x}}$ à $T(\dot{x})$
est un multiple du caractère
$\chi_{\dot{x}}=\chi_{aa^{\tau},2j+1,j+1}$~; en fait, on a
$(\sigma_{T,\dot{x}})_{\vert T(\dot{x})}=q'\chi_{\dot{x}}$.

Nous allons voir que la représentation $\sigma_{T,\dot{x}}$ est
équivalente à la représentation
$\gamma_{\dot{x}}=\Ind_{T(\dot{x})}^{T(\hat{x})}\chi_{\dot{x}}$.

On définit un opérateur
$T(\hat{x})$-équiva\-riant $\mathcal{I}_{\dot{x}}$ de l'espace
$\mathcal{H}_{\dot{x}}$ de $\gamma_{\dot{x}}$ dans
$\mathcal{H}^{\mathrm{x}}_{\overline{\psi}}$ en posant pour
$\varphi\in\mathcal{H}_{\dot{x}}$
\begin{equation}\label{eq7.3.6}
  \mathcal{I}_{\dot{x}}\varphi= \sum_{t\in
    T(\hat{x})/T(\dot{x})}\varphi(t)\sigma_{T,\dot{x}}(t^{-1})\delta_{0}.
\end{equation}
Comme $T(\hat{x})/T(\dot{x})=T_{j}/T_{j+1}$ est de cardinal $q'$, les
espaces $\mathcal{H}_{\dot{x}}$ et
$\mathcal{H}^{\mathrm{x}}_{\overline{\psi}}$ ont même dimension. Il
nous suffit donc de montrer que $\mathcal{I}_{\dot{x}}$ est
injectif. Soit donc $\varphi\in\mathcal{H}_{\dot{x}}$ tel que
$\mathcal{I}_{\dot{x}}\varphi=0$. Or, il suit des formules
\ref{eq7.3.5} et \ref{eq3.1.1} que, pour $t\in T(\hat{x})$,
$\sigma_{T,\dot{x}}(t^{-1})\delta_{0}$, vu comme un élément de
$\mathcal{H}_{\psi}^{A}$ est à support dans $(t^{-1}x-x)+A$. On voit
donc que, si $t$, $t'\in T(\hat{x})$,
$\sigma_{T,\dot{x}}(t^{-1})\delta_{0}$ et
$\sigma_{T,\dot{x}}(t'^{-1})\delta_{0}$ ont des supports disjoints sauf
si $t^{-1}x+A=t'^{-1}x+A$, c'est à dire si $t\ddot{x}=t'\ddot{x}$. Ce
qui équivaut à $tt'^{-1}\in T(\dot{x})$. Ceci montre que le système de
vecteurs $\sigma_{T,\dot{x}}(t^{-1})\delta_{0}$, $t$ parcourant un
ensemble de représentants dans $T(\hat{x})$ du quotient
$T(\hat{x})/T(\dot{x})$, est libre. Il suit alors de \ref{eq7.3.6} que
$\varphi=0$ comme voulu.

Il suit donc du théorème d'induction par étages que
\begin{equation*}
S_{T,\dot{x}}=\Ind_{T(\dot{x})}^{T}\chi_{\dot{x}}.
\end{equation*}
On en déduit que la représentation $S_{T,\dot{x}}$ se décompose en la
somme directe sans multiplicité des caractères de $T$ dont la
restriction à $T(\dot{x})$ est $\chi_{aa^{\tau},2j+1,j+1}$. Désignons
par $\mathcal{G}_{T,\dot{x}}$ l'espace de la représentation
$\Ind_{T(\dot{x})}^{T}\chi_{\dot{x}}$ et soit $\phi_{0}$ l'élément de
$\mathcal{G}_{T,\dot{x}}$ de support $T(\dot{x})$ et tel que
$\phi_{0}(1)=1$. Si $\chi$ est un caractère de $T$ dont la
restriction à $T(\dot{x})$ est $\chi_{aa^{\tau},2j+1,j+1}$, on désigne
par $\phi_{\chi}$ la fonction définie sur $T$ par 
\begin{equation*}
  \phi_{\chi}(g)=\sum_{t\in
    T/T(\dot{x})}\chi(t^{-1})\phi_{0}(gt)\mbox{, } g\in T.
\end{equation*}
Alors $\phi_{\chi}$ est un élément non nul de $\mathcal{G}_{T,\dot{x}}$
engendrant le sous-espace des vecteurs de poids $\chi$. 

D'autre part, si $\phi\in\mathcal{G}_{T,\dot{x}}$, on définit la fonction
$\widetilde{\phi}$  à valeurs dans
$\mathcal{H}_{\dot{x}}$ en posant
\begin{equation*}
\widetilde{\phi}(g)(t)=\phi(gt)\mbox{, }g\in T\mbox{, }t\in T(\hat{x}).
\end{equation*}
Alors l'application $\phi\mapsto\widetilde{\phi}$, est un isomorphisme
$T$-équivariant de $\mathcal{G}_{T,\dot{x}}$ sur l'espace de la
représentation
$\Ind_{T(\hat{x})}^{T}(\Ind_{T(\dot{x})}^{T(\hat{x})}\chi_{\dot{x}})$.
On définit un isomorphisme $T$-équivariant $\mathcal{J}_{T,\dot{x}}$
de $\mathcal{G}_{T,\dot{x}}$ sur $\mathcal{H}_{T,\dot{x}}$ en posant
pour $\phi\in\mathcal{G}_{T,\dot{x}}$
\begin{equation*}
\mathcal{J}_{T,\dot{x}}\phi(g)=\mathcal{I}_{\dot{x}}(\widetilde{\phi}(g))
\mbox{, }g\in T.
\end{equation*}
On vérifie que $\mathcal{J}_{T,\dot{x}}\phi_{0}$ est l'élément
$\varphi_{0}$ de $\mathcal{H}_{T,\dot{x}}$ de support $T(\hat{x})$ et
tel que $\varphi_{0}(1)=\delta_{0}$.

L'application
$\mathcal{K}_{T,\dot{x}}=
\mathcal{I}_{A,B}\circ\mathcal{I}_{T,\dot{x}}\circ\mathcal{J}_{T,\dot{x}}$
est une injection $T$-équivariante de $\mathcal{G}_{T,\dot{x}}$ dans
$\mathcal{H}^{A}_{\psi}$ dont l'image est le sous-espace
$\mathcal{E}_{\hat{x}}$ des fonctions à support dans l'orbite de
$\hat{x}$ sous $T$. On vérifie que 
\begin{equation*}
  \mathcal{K}_{T,\dot{x}}
\phi_{0}=\delta_{x}.
\end{equation*}
En effet, il suit de la relation \ref{eq7.3.3} que
$\mathcal{I}_{T,\dot{x}}\circ\mathcal{J}_{T,\dot{x}} \phi_{0}$ est
l'élément $\theta$ de $\mathcal{H}^{B}_{\psi}$ supporté par $\hat{x}$
et vérifiant $\theta(x)=\delta_{0}$. On en déduit que
$\mathcal{K}_{T,\dot{x}}
\phi_{0}=\mathcal{I}_{A,B}\theta$ est un élément de
$\mathcal{H}^{A}_{\psi}$ supporté par $\hat{x}$ et vérifiant, pour
tout $b\in B^{*}$,
\begin{equation*}
\begin{split}
\mathcal{I}_{A,B}\theta(x+b)& =\theta(x+b)(0)\\
&=
\psi(\frac{1}{2}\beta(x,b))\widetilde{\rho}_{\overline{\psi}}(b)(\delta_{0})(0)\\
&=\psi(\frac{1}{2}\beta(x,b))\delta_{0}(b).
\end{split}
\end{equation*}
Notre assertion est alors claire.

Maintenant soit $\chi$ un caractère de $T$ tel que
$\chi_{\vert T(\dot{x})}=\chi_{aa^{\tau},2j+1,j+1}$. Désignons par
$\varphi_{\chi}$ l'image de $\phi_{\chi}$ par
$\mathcal{K}_{T,\dot{x}}$. Il est clair que $\varphi_{\chi}$ est un
vecteur non nul de poids $\chi$ dans $\mathcal{H}^{A}_{\psi}$ pour
l'action de $T$ et qu'il vérifie la relation \ref{eq7.3.2}.

Comme $\sideset{}{_{k''/k'}}\norm$ est un morphisme surjectif de
groupes de $\mathcal{O}''^{\times}$ sur $\mathcal{O}'^{\times}$ et
comme toute $T$-orbite dans $\mathcal{O}''^{\times}$ rencontre
$\mathcal{O}''^{\times}\backslash(\varpi'\mathcal{O}'+\mathcal{O}'\nu)$,
on voit que les caractères de conducteur impair $2j+1$, $j>0$, de $T$ sont
exactement les caractères dont la restriction à $T_{j+1}$ est de la
forme $\chi_{aa^{\tau},2j+1,j+1}$ avec
$a\in\mathcal{O}''^{\times}\backslash(\varpi'\mathcal{O}'+\mathcal{O}'\nu)$.

Nous avons donc démontré que chaque caractère d'indice impair
strictement supérieur à $1$ apparaît avec la multiplicité $1$ dans
$\mathcal{H}^{A}_{\psi}$ et qu'une fonction propre correspondante est
celle donnée dans l'énoncé du théorème. Nous avons également montré
que la somme directe des sous-espaces propres correspondants est égale
à
\begin{equation*}
  \bigoplus_{x\in X}S_{T,\dot{x}}.
\end{equation*}

Compte tenu de la relation \ref{eq7.3.4}, il nous suffit donc, pour
achever la démonstration du théorème, de montrer que les seuls
caractères de $T$ qui apparaissent dans $\mathcal{E}^{B}_{0}$ sont les
caractères de conducteur au plus $1$, à l'exclusion du caractère
$\eta_{0}$.

Comme l'action de $K_{B}$ dans $\mathcal{E}^{B}_{0}$ passe au quotient
en la représentation de Weil $S^{\mathrm{x}}_{\overline{\psi}}$ de
$Sp(\mathrm{b}^{*})$, les caractères de $T$ qui apparaissent sont des
caractères qui passent au quotient à $T/T_{1}=\mu_{q'+1}$, donc de
conducteur au plus $1$. Comme ces caractères apparaissent avec au plus
la multiplicité $1$ alors que la représentation de Weil est de
dimension $q'$, on voit qu'ils apparaissent tous sauf l'un d'entre
eux, que l'on note $\eta_{0}$ et que nous allons déterminer. Compte
tenu de ce qui précède, il est clair que l'on a
\begin{equation}\label{eq7.3.7}
\eta_{0}(g)=
\begin{cases}
1 & \mbox{ si }g\in T_{1}\\
-\tr S_{\overline{\psi}}^{\mathrm{x}}(p_{Sp(\mathrm{b}^{*})}(g)) &
\mbox{ si }g\in T\backslash T_{1}.
\end{cases}
\end{equation}

Nous devons donc calculer la restriction du caractère de la
représentation $\widetilde{S}_{\overline{\psi}}^{\mathrm{x}}$ au tore
$T$. Pour ce faire, nous allons commencer par décrire les opérateurs
$S_{\overline{\psi}}^{\mathrm{x}}(p_{Sp(\mathrm{b}^{*})}(g))=
\widetilde{S}_{\overline{\psi}}^{\mathrm{x}}(g)$, $g\in T$. En fait
ces opérateurs sont la restriction des opérateurs $S_{\psi}^{A}(g)$,
$g\in T$ au sous-espace $\mathcal{E}^{B}_{0}$.

Dans la suite, si $x$ est un élément de $k''$ (resp. $\mathcal{O}''$),
on désigne par $\dot{x}$ son image dans $k''/B$
(resp. $\mathbb{F}_{q''}$). On identifie $\mathrm{b}^{*}=B^{*}/B$ avec
$\mathbb{F}_{q''}$ au moyen de l'application $\dot{x}\mapsto
p_{\mathbb{F}_{q''}}(\varpi'^{\frac{1-\mu}{2}}x)$. 

On choisit des bases $e_{1},\ldots,e_{r}$ et $f_{1},\ldots,f_{r}$ de
$\mathcal{O}'$ sur $\mathcal{O}$ telles que
\begin{equation*}
\varpi'^{\frac{\mu-1}{2}}e_{1},\ldots,\varpi'^{\frac{\mu-1}{2}}e_{r},
\varpi'^{\frac{\mu+1}{2}}f_{1}\nu,\ldots,\varpi'^{\frac{\mu+1}{2}}f_{r}\nu
\end{equation*}
soit une base autoduale de $k''$ sur $k$ et que
$\dot{e}_{1},\ldots,\dot{e}_{l}$, où $l=\frac{r}{e}$, soit une base de
$\mathbb{F}_{q'}$ sur $\mathbb{F}_{q}$. On vérifie alors que
$\dot{f}_{1},\ldots,\dot{f}_{l}$ est également une base de
$\mathbb{F}_{q'}$ sur $\mathbb{F}_{q}$ et que
$\dot{e}_{1},\ldots,\dot{e}_{l},\dot{f}_{1}\dot{\nu},\ldots,\dot{f}_{l}\dot{\nu}$
est une base symplectique de $\mathbb{F}_{q''}=\mathrm{b}^{*}$ pour la
forme symplectique $\beta_{\mathrm{b}^{*}}$. De plus,
$\dot{e}_{1},\ldots,\dot{e}_{l}$ est une base du lagrangien
$\mathrm{x}=\mathbb{F}_{q'}$ dont le sous-espace
$\mathrm{y}=\mathbb{F}_{q'}\dot{\nu}$ de $\mathrm{b}^{*}$, engendré
par $\dot{f}_{1}\dot{\nu},\ldots,\dot{f}_{l}\dot{\nu}$, est un
supplémentaire lagrangien.

Dans la suite, on identifie $\mathbb{F}_{q'}^{2}$ avec
$\mathrm{b}^{*}=\mathbb{F}_{q''}$ au moyen de l'application $(x,y)\mapsto
x+y\dot{\nu}$ et on écrit tout endomorphisme $\mathbb{F}_{q}$-linéaire
de $\mathrm{b}^{*}$ sous forme matricielle 
$\left(\begin{smallmatrix}
x & y\\
z & t
\end{smallmatrix}\right)$
avec $x$, $y$, $z$ et $t$ des endomorphismes
$\mathbb{F}_{q}$-linéaires de $\mathbb{F}_{q'}$.

L'application $\gamma:\mathcal{O}'\longrightarrow\mathbb{F}_{q}$
définie par 
\begin{equation*}
\gamma(x)=
p_{\mathbb{F}_{q}}(\varpi^{1-\lambda_{\psi}}
\sideset{}{_{k'/k}}\tr\varpi'^{\mu-1+v''(u)}x)
\end{equation*}
passe au quotient en une forme $\mathbb{F}_{q}$-linéaire non nulle sur
$\mathbb{F}_{q'}$. Il existe donc
$\gamma_{0}\in\mathbb{F}_{q'}^{\times}$ tel que 
\begin{equation}\label{eq7.3.8}
\gamma(x)=\sideset{}{_{\mathbb{F}_{q'}/\mathbb{F}_{q}}}\tr\gamma_{0}\dot{x}\mbox{,
}x\in\mathcal{O}'.
\end{equation}

Cela étant, la forme symplectique $\beta_{\mathrm{b}^{*}}$ induit une
dualité entre les sous-espaces lagrangiens $\mathrm{x}$ et
$\mathrm{y}$ et donc une dualité du $\mathbb{F}_{q}$-espace vectoriel
$\mathbb{F}_{q'}$ sur lui-même, donnée par
\begin{equation*}
  \langle x,y\rangle=\beta_{\mathrm{b}^{*}}(x,y\dot{\nu})\mbox{,
  }x,y\in\mathbb{F}_{q'}.
\end{equation*}
Utilisant les définitions \ref{eq6.1.1} de $\beta$ et \ref{eq2.3.3} de
$\beta_{\mathrm{b}^{*}}$ et compte tenu de la formule \ref{eq7.3.8},
on vérifie que cette dualité est une forme bilinéaire symétrique
donnée par
\begin{equation}\label{eq7.3.9}
\langle x,y\rangle=
\sideset{}{_{\mathbb{F}_{q'}/\mathbb{F}_{q}}}\tr\gamma_{0}\dot{d}xy\mbox{,
}x,y\in\mathbb{F}_{q'}.
\end{equation}

Si $a$ est un $\mathbb{F}_{q}$-endomorphisme de $\mathbb{F}_{q'}$, on
désigne par ${}^ta$ son adjoint relativement à la forme bilinéaire
symétrique $\langle\, ,\rangle$. Alors, un endomorphisme $g=\left(
\begin{smallmatrix}
x & y\\
z & t
\end{smallmatrix}\right)
$
est un élément de $Sp(\mathrm{b}^{*})$ si et seulement si
\begin{equation*}
{}^txz={}^tzx\mbox{, }{}^tyt={}^tty\mbox{ et }{}^txt-{}^tzy=Id.
\end{equation*}

Soit $\varsigma_{\mathrm{b}^{*}}$ l'image de $\varsigma_{B}$ par le
morphisme $p_{Sp(\mathrm{b}^{*})}$. Alors, on a
$\varsigma_{\mathrm{b}^{*}}=\left(
\begin{smallmatrix}
0 & -v^{-1}\\
v & 0
\end{smallmatrix}\right)
$ où $v$ est l'endomorphisme de $\mathbb{F}_{q'}$ tel que
$v(\dot{e}_{i})=-\dot{f}_{i}$, $1\leq i\leq l$. Cet endomorphisme
vérifie ${}^tv=v$, puisque $\dot{f}_{1},\ldots,\dot{f}_{l}$ est la
base duale de la base $\dot{e}_{1},\ldots,\dot{e}_{l}$, relativement à
la forme bilinéaire symétrique $\langle\, ,\rangle$.

De même, si $g=x+y\nu\in T$, on voit que, comme endomorphisme de
$\mathbb{F}_{q''}$ agissant par multiplication,  $\dot{g}=\left(
\begin{smallmatrix}
\dot{x} & \dot{d}\dot{y}\\
\dot{y} & \dot{x}
\end{smallmatrix}\right)
$.

Si $s\in\mathcal{O}'$, la fonction
$\delta_{\varpi'^{\frac{\mu-1}{2}}s\nu}$ ne dépend que de
$\dot{s}\in\mathbb{F}_{q'}$ et on la note $\varphi_{\dot{s}}$. Alors,
les $\varphi_{s}$, $s\in\mathbb{F}_{q'}$, forment une base de l'espace
$\mathcal{E}^{B}_{0}$. On rappelle que $\omega$ désigne l'indice de
Weil relatif au caractère $\frac{1}{2}\overline{\psi}$ du corps
$\mathbb{F}_{q}$. On désigne par $\omega'$ l'indice de Weil relatif au
caractère
$\frac{1}{2}\overline{\psi}\circ\sideset{}{_{\mathbb{F}_{q'}/\mathbb{F}_{q}}}\tr$
du corps $\mathbb{F}_{q'}$.
\begin{lem}\label{lem7.3.1}
a) On a 
\begin{equation}\label{eq7.3.10}
  S_{\psi}^{A}(-1)\varphi_{s}=\left(\frac{-1}{q'}\right)\varphi_{-s}\mbox{,
  }s\in\mathbb{F}_{q'}.
\end{equation}

b) Soit $g\in T$ tel que $\dot{g}\neq\pm1$. Si l'on écrit $g=x+y\nu$,
avec $x$, $y\in\mathcal{O}'$ tels que $x^{2}-dy^{2}=1$, pour tout
$s\in\mathbb{F}_{q'}$, on a
\begin{equation}\label{eq7.3.11}
S_{\psi}^{A}(g)\varphi_{s}=
q^{\frac{-l}{2}}\omega'(1)^{-1}\left(\frac{\gamma_{0}\dot{d}\dot{y}}{q'}\right)
\sum_{t\in\mathbb{F}_{q'}}\overline{\psi}
(\frac{1}{2}\sideset{}{_{\mathbb{F}_{q'}/\mathbb{F}_{q}}}
\tr
\gamma_{0}\dot{d}\dot{y}^{-1}(\dot{x}(s^{2}+t^{2})-2st))\varphi_{t}.
\end{equation}
\end{lem}
\begin{dem}
a) Cette assertion se déduit de la formule \ref{eq4.5.7} du corollaire
\ref{co4.5.1} et de ce que
$\left(\frac{(-1)^{l}}{q}\right)=\left(\frac{-1}{q'}\right)$.

b) On vérifie que $g=p_{1}\varsigma_{B}p_{2}$, avec $p_{1}$,
$p_{2}\in P_{B}$ tels que 
\begin{equation*}
  p_{Sp(\mathrm{b}^{*})}(p_{1})=\left(
\begin{smallmatrix}
\dot{y}^{-1}v & \dot{x}v^{-1}\\
0 & \dot{y}v^{-1}
\end{smallmatrix}\right) \mbox{ et }
p_{Sp(\mathrm{b}^{*})}(p_{2})=\left(
\begin{smallmatrix}
1 & \dot{y}^{-1}\dot{x}\\
0 & 1
\end{smallmatrix}\right).
\end{equation*}
Grâce aux formules \ref{eq4.5.7} et \ref{eq4.5.8}, on en déduit que
\begin{equation}\label{eq7.3.12}
  S_{\psi}^{A}(g)=\left(\frac{-1}{q}\right)^{l}
\omega(1)^{-l}q^{\frac{l}{2}}\left(\frac{\det\dot{y}^{-1}v}{q}\right)
M_{A}(p_{1})M_{A}(\varsigma_{B})M_{A}(p_{2}).
\end{equation}
Commençons par calculer $\left(\frac{\det\dot{y}^{-1}v}{q}\right)$. Soit
$\dot{e}_{1}^{*},\ldots,\dot{e}_{l}^{*}$ la base de $\mathbb{F}_{q'}$
sur $\mathbb{F}_{q}$ duale de la base $\dot{e}_{1},\ldots,\dot{e}_{l}$
pour la forme quadratique
$x\mapsto\tr_{\mathbb{F}_{q'}/\mathbb{F}_{q}}x^{2}$. On a
$\dot{e}_{i}^{*}=\gamma_{0}\dot{d}\dot{f}_{i}$, de sorte que
$v\dot{e}_{i}=-(\gamma_{0}\dot{d})^{-1}\dot{e}_{i}^{*}$, $1\leq i\leq
l$. Soit $\Delta(\mathbb{F}_{q'}/\mathbb{F}_{q})$ le discriminant de
$\mathbb{F}_{q'}$ sur $\mathbb{F}_{q}$. D'après le théorème de parité de
Stickelberger (voir \cite{stickelberger-ip-1897}), on a 
\begin{equation*}
\left(\frac{\Delta(\mathbb{F}_{q'}/\mathbb{F}_{q})}{q}\right)=(-1)^{l-1}.
\end{equation*}
On déduit alors de ce qui précède que 
\begin{equation*}
\det\dot{y}^{-1}v \equiv
(-1)^{l}\sideset{}{_{\mathbb{F}_{q'}/\mathbb{F}_{q}}}\norm
(\gamma_{0}\dot{d}\dot{y})^{-1}\Delta(\mathbb{F}_{q'}/\mathbb{F}_{q})
\mod (\mathbb{F}_{q}^{\times})^{2}.
\end{equation*}
Mais alors, on a 
\begin{equation*}
\begin{split}
\left(\frac{\det\dot{y}^{-1}v}{q}\right)&=
\left(\frac{\sideset{}{_{\mathbb{F}_{q'}/\mathbb{F}_{q}}}\norm
  (\gamma_{0}\dot{d}\dot{y})}{q}\right)\left(\frac{-1}{q}\right)^{l}
\left(\frac{\Delta(\mathbb{F}_{q'}/\mathbb{F}_{q})}{q}\right)\\
&=(-1)^{l-1}\left(\frac{\gamma_{0}\dot{d}\dot{y}}{q'}\right)
\left(\frac{-1}{q}\right)^{l}.
\end{split}
\end{equation*}
Reportant cette égalité dans la formule \ref{eq7.3.12} et tenant
compte tenu de la relation de Hasse-Davenport (voir \cite[paragraphe
  18]{weil-1974} et aussi \cite{perrin-ip-1981})
\begin{equation*}
  -\omega'(1)=(-\omega(1))^{l}
\end{equation*}
il vient 
\begin{equation}\label{eq7.3.13}
  S_{\psi}^{A}(g)=
\omega'(1)^{-1}q^{\frac{l}{2}}\left(\frac{\gamma_{0}\dot{d}\dot{y}}{q'}\right)
M_{A}(p_{1})M_{A}(\varsigma_{B})M_{A}(p_{2}).
\end{equation}

Maintenant, nous allons calculer l'action des opérateurs
$M_{A}(p_{2})$, $M_{A}(\varsigma_{B})$ et $M_{A}(p_{1})$ dans
$\mathcal{E}^{B}_{0}$.

Si $s\in\mathcal{O}'$, on a $p_{2}\varpi'^{\frac{\mu-1}{2}}s\nu\in
\varpi'^{\frac{\mu-1}{2}}(s\nu+y^{-1}xs)+B$. Compte tenu de la formule
\ref{eq4.4.3}, il vient
\begin{equation*}
\begin{split}
M_{A}(p_{2})\varphi_{\dot{s}}&=
\psi(\frac{1}{2}\beta(\varpi'^{\frac{\mu-1}{2}}y^{-1}xs,
\varpi'^{\frac{\mu-1}{2}}s\nu))\varphi_{\dot{s}}\\
&=\psi(\frac{1}{2}\sideset{}{_{k'/k}}\tr\varpi'^{\mu-1+v''(u)}dy^{-1}xs^{2})
\varphi_{\dot{s}}\\
&=\overline{\psi}(p_{\mathbb{F}_{q}}(\frac{1}{2}\varpi^{1-\lambda_{\psi}}
\sideset{}{_{k'/k}}\tr\varpi'^{\mu-1+v''(u)}dy^{-1}xs^{2}))
\varphi_{\dot{s}}.
\end{split}
\end{equation*}
On en déduit immédiatement que, pour tout $s\in\mathbb{F}_{q'}$, on a
\begin{equation}\label{eq7.3.14}
  M_{A}(p_{2})\varphi_{s}=
\overline{\psi}(\frac{1}{2}\sideset{}{_{\mathbb{F}_{q'}/\mathbb{F}_{q}}}\tr
\gamma_{0}\dot{d}\dot{y}^{-1}\dot{x}s^{2})\varphi_{s}.
\end{equation}

Si $s\in\mathbb{F}_{q'}$, on désigne par $\tilde{s}$ un de ses
représentants dans $\mathcal{O}'$. Soit donc
$s\in\mathbb{F}_{q'}$. Alors, on a
$p_{1}\varpi'^{\frac{\mu-1}{2}}\tilde{s}\nu\in
\varpi'^{\frac{\mu-1}{2}}((\dot{x}v^{-1}s)\tilde{\,}+
(\dot{y}v^{-1}s)\tilde{\,}\nu)+B$. Comme pour $M_{A}(p_{2})$, on en
déduit que
\begin{equation}\label{eq7.3.15}
M_{A}(p_{1})\varphi_{s}=
\overline{\psi}(\frac{1}{2}\sideset{}{_{\mathbb{F}_{q'}/\mathbb{F}_{q}}}\tr
\gamma_{0}\dot{d}\dot{x}\dot{y}(v^{-1}s)^{2})\varphi_{\dot{y}v^{-1}s}.
\end{equation}

Pour calculer l'action de $M_{A}(\varsigma_{B})$, on commence par
remarquer que $\varsigma_{B}A=A\nu$ et $A\cap\varsigma_{B}A=B=B\nu$ de
sorte que le groupe quotient $\varsigma_{B}A/A\cap\varsigma_{B}A$ est
isomorphe à $\mathbb{F}_{q'}=A/B$, l'isomorphisme étant donné par
l'application $t\mapsto\varpi'^{\frac{\mu-1}{2}}\tilde{t}\nu+B$. De
plus, pour $t\in\mathbb{F}_{q'}$, on a
\begin{eqnarray*}
\varsigma_{B}(\varpi'^{\frac{\mu-1}{2}}\tilde{t}\nu+B) &=&
-\varpi'^{\frac{\mu-1}{2}}(v^{-1}t)\tilde{\,}+B\\
\varsigma_{B}^{-1}(\varpi'^{\frac{\mu-1}{2}}\tilde{t}\nu+B)&=&
\varpi'^{\frac{\mu-1}{2}}(v^{-1}t)\tilde{\,}+B.
\end{eqnarray*}
Utilisant le lemme \ref{lem4.4.1}, on en déduit que, pour
$s\in\mathbb{F}_{q'}$, 
\begin{equation}\label{eq7.3.16}
M_{A}(\varsigma_{B})\varphi_{s}=q^{-l}\sum_{t\in\mathbb{F}_{q'}}
\overline{\psi}(-\frac{1}{2}\sideset{}{_{\mathbb{F}_{q'}/\mathbb{F}_{q}}}\tr
\gamma_{0}\dot{d}(sv^{-1}t+tv^{-1}s))\varphi_{t}.
\end{equation}
Soit $s\in\mathbb{F}_{q'}$. Comme $v$ est un endomorphisme de
$\mathbb{F}_{q'}$ symétrique relativement à la forme $\langle
x,y\rangle=\sideset{}{_{\mathbb{F}_{q'}/\mathbb{F}_{q}}}\tr\gamma_{0}\dot{d}xy$,
on déduit des relations \ref{eq7.3.14}, \ref{eq7.3.15} et
\ref{eq7.3.16} que
\begin{equation*}
\begin{split}
  M_{A}(p_{1})&M_{A}(\varsigma_{B})M_{A}(p_{2})\varphi_{s}\\
&=q^{-l}\sum_{t\in\mathbb{F}_{q'}}\overline{\psi}
(\frac{1}{2}\sideset{}{_{\mathbb{F}_{q'}/\mathbb{F}_{q}}}\tr
\gamma_{0}\dot{d}
(\dot{y}^{-1}\dot{x}s^{2}-2sv^{-1}t+\dot{x}\dot{y}(v^{-1}t)^{2}))
\varphi_{\dot{y}v^{-1}t}\\
& = q^{-l}\sum_{t\in\mathbb{F}_{q'}}\overline{\psi}
(\frac{1}{2}\sideset{}{_{\mathbb{F}_{q'}/\mathbb{F}_{q}}}\tr
\gamma_{0}\dot{d}\dot{y}^{-1}(\dot{x}(s^{2}+t^{2})-2st))\varphi_{t}.
\end{split}
\end{equation*}
Notre lemme en résulte, compte tenu des formules \ref{eq7.3.12} et
\ref{eq7.3.13}.
\end{dem}

Il résulte des relations \ref{eq7.3.7} et \ref{eq7.3.10} que 
\begin{equation}\label{eq7.3.17}
\eta_{0}(-1)=-\left(\frac{-1}{q'}\right)=\left(\frac{-d}{q'}\right).
\end{equation}

D'autre part, soit $g\in T$ tel que $\dot{g}\neq\pm1$. \'Ecrivons
$g=x+y\nu$ avec $x$, $y\in\mathcal{O}'$ tels que $x^{2}-dy^{2}=1$. Il
suit de la formule \ref{eq7.3.11} que, pour $t\in\mathbb{F}_{q'}$, on
a 
\begin{equation*}
\langle S_{\psi}^{A}(g)\varphi_{t},\varphi_{t}\rangle=
q^{-\frac{l}{2}}\omega'(1)^{-1}\left(\frac{\gamma_{0}\dot{d}\dot{y}}{q'}\right)
\overline{\psi}(\sideset{}{_{\mathbb{F}_{q'}/\mathbb{F}_{q}}}\tr\gamma_{0}
\dot{d}\frac{\dot{x}-1}{\dot{y}}t^{2}).
\end{equation*}
Compte tenu des formules \ref{eq3.3.-1}, \ref{eq3.3.-2} et \ref{eq7.3.7}
et de \cite[Proposition A.9 (1)]{rao-1993}, on en déduit que
\begin{equation*}
\begin{split}
\eta_{0}(g)&=-q^{-\frac{l}{2}}\omega'(1)^{-1}
\left(\frac{\gamma_{0}\dot{d}\dot{y}}{q'}\right)
\sum_{t\in\mathbb{F}_{q'}}
\overline{\psi}(\sideset{}{_{\mathbb{F}_{q'}/\mathbb{F}_{q}}}\tr\gamma_{0}
\dot{d}\frac{\dot{x}-1}{\dot{y}}t^{2})\\
&=-\left(\frac{\gamma_{0}\dot{d}\dot{y}}{q'}\right)
\left(\frac{2\gamma_{0}\dot{d}(\dot{x}-1)\dot{y}^{-1}}{q'}\right).
\end{split}
\end{equation*}
Comme par ailleurs $\left(\frac{\dot{d}}{q'}\right)=-1$, il vient
\begin{equation}\label{eq7.3.18}
\eta_{0}(g)=\left(\frac{2\dot{d}(\dot{x}-1)}{q'}\right).
\end{equation}
On remarque que la formule \ref{eq7.3.18} est encore valable lorsque
$\dot{g}=-1$ (voir la formule \ref{eq7.3.17}) et même lorsque
$\dot{g}=1$, à condition dans ce cas de convenir que
$\left(\frac{0}{q'}\right)=1$.

Maintenant, soit $z=\xi+\zeta\nu\in\mathcal{O}''^{\times}$ tel que
$g=z(z^{-1})^{\tau}$. Alors, on a 
\begin{equation*}
2\dot{d}(\dot{x}-1)=
\frac{4\dot{d}^{2}\dot{\zeta}^{2}}{\dot{\xi}^{2}-\dot{d}\dot{\zeta}^{2}}
\end{equation*}
de sorte que 
\begin{equation*}
\eta_{0}(g)=\left(\frac{\dot{\xi}^{2}-\dot{d}\dot{\zeta}^{2}}{q'}\right)
=\left(\frac{p_{\mathbb{F}_{q'}}(zz^{\tau})}{q'}\right).
\end{equation*}
D'où le théorème, d'après la formule \ref{eq7.1.0}.
\end{dem}

\begin{rem}
Les théorèmes \ref{theo7.2.1} et \ref{theo7.3.1} ont été démontrés par
Yang dans le cas où $W$ est de dimension $2$ (voir
\cite{yang-1998}).
\end{rem}

\subsection{}\label{7.4}
Soit $T$ un tore maximal elliptique de $Sp(W)$. Dans ce paragraphe,
nous décrivons la restriction à $T$ de la représentation de Weil de
$Mp(W)$. 

Rappelons que (voir le paragraphe \ref{6.1}) le $T$-module $W$ se
décompose de manière unique en somme directe orthogonale de
sous-$T$-modules symplectiques et irréductibles sur $k$,
$W=W_{1}\oplus\cdots\oplus W_{n}$, que, pour $1\leq i\leq n$, l'image
$T_{i}$ de $T$ par l'application $x\mapsto x_{\vert W_{i}}$ est un
tore maximal irréductible de $Sp(W_{i})$ et que l'on identifie $T$
avec le produit direct de tores $\prod_{i=1}^{i=n}T_{i}$ au moyen de
l'isomorphisme $x\mapsto(x_{\vert W_{1}},\ldots,x_{\vert W_{n}})$.

Soit $1\leq i\leq n$. 
On désigne par $B_{i}$ le bon réseau
$T_{i}$-invariant de $W_{i}$ introduit dans le corollaire
\ref{co6.2.1} et par $s_{i}$ la restriction à $T_{i}$ du relèvement
$s_{B_{i}}$ de $K_{B_{i}}$ dans $Mp(W_{i})$. On pose $A_{i}=B_{i}$ si
le réseau $B_{i}$ est autodual. Dans le cas contraire, $A_{i}$ désigne
le réseau autodual introduit au début du numéro
\ref{7.3} et vérifiant $B_{i}\subset A_{i}\subset B^{*}_{i}$.

Alors $B=B_{1}\oplus\cdots\oplus B_{n}$ est un bon réseau, 
$A=A_{1}\oplus\cdots\oplus A_{n}$ est un réseau autodual de $W$ et ils
vérifient $B\subset A\subset B^{*}$. De plus, $B$ est stable sous
l'action de $T$ de sorte que $T\subset K_{B}$. On identifie $T$ à un
sous-groupe de $Mp(W)$ à l'aide de la section $s_{B}$ de $K_{B}$ dans
$Mp(W)$. 

On rappelle que si $\varphi_{i}\in\mathcal{H}_{\psi}^{A_{i}}$, $1\leq
i\leq n$, on a défini l'élément
$\varphi_{1}\otimes\cdots\otimes\varphi_{n}$ de
$\mathcal{H}_{\psi}^{A}$ à l'aide de la formule \ref{eq4.7.1}
\begin{theo}\label{theo7.4.1}
(i) Un caractère $\chi$ de $T$ intervient dans la restriction à $T$ de
  la représentation de Weil $S_{W,\psi}$ si et seulement si, pour
  chaque $1\leq i\leq n$, le caractère $\chi_{\vert T_{i}}$ de $T_{i}$
  intervient dans restriction à $T_{i}$ de
  la représentation de Weil $S_{W_{i},\psi}$.

(ii) Soit $\chi$ un caractère de $T$ vérifiant les conditions
  équivalentes du (i).  Pour $1\leq i\leq n$, on désigne par
  $\chi_{i}$ la restriction de $\chi$ à $T_{i}$ et on choisit une
  fonction propre 
  $\varphi_{\chi_{i}}\in\mathcal{H}_{\psi}^{A_{i}}$ de poids
  $\chi_{i}$ relativement à la restriction de la représentation de Weil
  $S_{W_{i},\psi}$ à $T_{i}$. Alors
\begin{equation*}
\varphi_{\chi}:=
\varphi_{\chi_{1}}\otimes\cdots\otimes\varphi_{\chi_{n}}\
\end{equation*}
est une fonction propre de poids $\chi$ dans $\mathcal{H}_{\psi}^{A}$
relativement à la restriction de la représentation de Weil $S_{W,\psi}$
à $T$.
\end{theo}
\begin{dem}
C'est une conséquence immédiate de la proposition \ref{pr4.7.1}.
\end{dem}

\begin{rem}
Le théorème \ref{theo7.4.1}, combiné avec les théorèmes
\ref{theo7.2.1} et \ref{theo7.3.1}, permet de décrire explicitement la
restriction de la représentation de Weil à un tore maximal elliptique
donné de $Sp(W)$.
\end{rem}

\bibliographystyle{smfplain} \bibliography{biblio}

\end{document}